\documentclass[final,1p]{elsarticle}
\usepackage{amssymb}
\usepackage{amsmath}
\usepackage{amsthm}
\usepackage{lineno}
\usepackage{latexsym}
\usepackage{cases}
\usepackage{bm}
\usepackage[ruled,lined]{algorithm2e}
\usepackage{algorithmic}
\usepackage{array}
\usepackage[colorlinks,linkcolor=blue, anchorcolor=blue, citecolor=blue]{hyperref}
\usepackage{color,tabularx}
\usepackage{booktabs}
\usepackage{subfigure}
\usepackage{multirow}
\usepackage{stmaryrd}
\usepackage[T1]{fontenc}
\biboptions{numbers,sort&compress}
\usepackage{graphicx} %use graph format
\usepackage{epstopdf}
\usepackage{adjustbox}
\graphicspath{{./figure/}}
\allowdisplaybreaks

\newtheorem{theorem}{Theorem}[section]

\newtheorem{lemma}[theorem]{Lemma}

\theoremstyle{definition}

\newtheorem{example}[theorem]{Example}

\theoremstyle{remark}
\newtheorem{remark}[theorem]{Remark}

\newcommand{\f}{\frac}

\newcommand{\p}{\partial}

\newcommand{\mal}{\mathcal}
\newcommand{\mbb}{\mathbb}

\newcommand{\mbX}{\mbb{X}}
\newcommand{\mrR}{\mathrm{R}}
\newcommand{\mrE}{\mathrm{E}}

\newcommand{\mrM}{\mathrm{M}}
\numberwithin{equation}{section}

\begin{document}    
	\begin{frontmatter}
		
		% \title{A linear, positivity-preserving, mass-conservation, original energy dissipation and convergent analysis block-centered finite difference method for the multi-species Keller--Segel chemotaxis system}

        \title{A linearly decoupled structure-preserving BCFD method for the two-species Keller--Segel chemotaxis system}
		
		\author[OUC]{Ao Zhang} \ead{zhangao6290@stu.ouc.edu.cn}
		\author[OUC]{Bingyin Zhang} \ead{zhangbingyin@stu.ouc.edu.cn}
		\author[OUC,LAB]{Hongfei Fu\corref{Fu}}\ead{fhf@ouc.edu.cn}
		
		\address[OUC]{School of Mathematical Sciences, Ocean University of China, Qingdao 266100, P.R. China}
		\address[LAB]{Laboratory of Marine Mathematics, Ocean University of China, Qingdao 266100, P.R. China}
		
		\cortext[Fu]{Corresponding author.}
		
	\begin{abstract}
        In this paper, we develop a fully discrete, linearly decoupled block-centered finite difference scheme for the two-species Keller--Segel chemotaxis system. Based on the Slotboom transformation, the cell-density equations are reformulated into a symmetric form, providing a natural framework for the construction of a structure-preserving discretization. The resulting scheme is shown to unconditionally preserve the positivity of both the cell densities and the chemoattractant concentration, while ensuring exact mass conservation of each cell density and the original energy dissipation. Furthermore, optimal-order error estimates for the proposed scheme are rigorously established through a higher-order consistency analysis. In particular, the cell densities exhibit first-order temporal and second-order spatial convergence in the discrete $L^2$-norm, whereas the chemoattractant concentration achieves the same convergence orders in the discrete $H^1$-norm. Finally, a series of numerical experiments is conducted to verify the theoretical results and demonstrate the reliability of the proposed scheme in simulating the finite-time blow-up phenomenon.
	\end{abstract}
		
   \end{frontmatter}
	
%%%%%%%%%%%%%%%%%%%%%%%%%%%%%%%%%%%%%%%%%%%%%%%%%%%%%%%%%%%%%%
	\section{Introduction}
%%%%%%%%%%%%%%%%%%%%%%%%%%%%%%%%%%%%%%%%%%%%%%%%%%%%%%%%%%%%
    Chemotaxis refers to the directional migration of cells or microorganisms in response to spatial variations in chemical signals in their surrounding environment. The mathematical study of this biological process dates back to the seminal works of Keller and Segel \cite{KS'70,KS'71}, who introduced a continuum framework for describing interactions between a mobile population and chemical substances. Their formulation has since become a fundamental prototype for investigating chemotactic aggregation and related collective behaviors. Depending on the biological setting, chemical signals may induce either attractive or repulsive responses, thereby influencing the spatial redistribution of the population. Mathematically, the two-species Keller--Segel chemotaxis model considered here seeks the densities of two interacting cell populations, $u = u(\bm{x}, t)$ and $v = v(\bm{x}, t)$, and the chemoattractant concentration $c = c(\bm{x}, t)$, such that
	\begin{equation}\label{model:KS}
		\left\{
		\begin{aligned}
			&  \p_t u  = \Delta u-\nabla \cdot(u \nabla c),      & \qquad \text{in}~~ \Omega \times(0, T],\\
			&\p_t v  = \Delta v- \nabla \cdot(v \nabla c),      & \qquad \text{in}~~ \Omega \times(0, T],\\
			&\epsilon\, \p_t c  = \beta \Delta c- \alpha c+  u  +  v, & \qquad \text{in}~~ \Omega \times(0, T],
		\end{aligned}
		\right.
	\end{equation}
	subject to periodic boundary conditions and initial conditions
	\begin{equation}\label{model:KS:ic}
		\begin{aligned}
			&u(\bm{x},0)  =u^0(\bm{x}), \quad v(\bm{x},0)  =v^0(\bm{x}), \quad \epsilon\, c(\bm{x},0) = \epsilon\, c^0(\bm{x}),  & & \text { in } \Omega,
		\end{aligned}
	\end{equation}
	where $\bm{n}$ represents the unit outward normal vector on the boundary; when $\epsilon=0$, only the initial conditions $u^0(\bm{x})$ and $v^0(\bm{x})$ are required. Here $\Omega \subset \mrR^d$ ($d=2,3$) is a convex, bounded and open domain. The parameter $\alpha \geq 0$ is the consumption rate of chemoattractant, $\beta$ is a positive diffusion coefficient, and $\epsilon \geq 0$ describes the response rate of the chemoattractant concentration to the cell density. 
	
	Significantly, the Keller--Segel chemotaxis system \eqref{model:KS}--\eqref{model:KS:ic}  satisfies the mass conservation law for the cell densities, i.e.,
	\begin{equation}\label{model:KS:MC}
		\begin{aligned}
			\mrM[u](t) &:= \int_\Omega u(\bm{x}, t)  d\bm{x} = \int_\Omega u^0(\bm{x})  d\bm{x}=\mrM[u](0),\\
			\mrM[v](t) &:= \int_\Omega v(\bm{x}, t)  d\bm{x} = \int_\Omega v^0(\bm{x})  d\bm{x}= \mrM[v](0).
		\end{aligned}
	\end{equation}
	In addition, for positive regular initial data, i.e., $u^0(\bm{x}) > 0$, $v^0(\bm{x}) > 0$, and $c^0(\bm{x}) > 0$, the Keller--Segel chemotaxis system admits unique positive solutions, so that
	\begin{equation}\label{model:KS:positivity}
		\begin{aligned}
			u(\bm{x},t) > 0, \quad v(\bm{x},t) > 0, \quad c(\bm{x},t) > 0, \quad \text { in } \Omega.
		\end{aligned}
	\end{equation}
	Moreover, the Keller--Segel chemotaxis system admits a Wasserstein gradient-flow structure driven by the total free energy
	\begin{equation}\label{def:Energy}
		\begin{aligned}
			\mrE[u,v,c] := \int_\Omega \left[f(u) + f(v) - u c - v c +\f{\beta}{2} |\nabla c|^2 + \f{\alpha}{2} c^2 \right] d\bm{x},
		\end{aligned}
	\end{equation}
	where $f(\rho) := \rho \log(\rho) - \rho$ for $\rho \in (0,+\infty)$, with $\rho=u$ and $v$. A direct calculation shows that the two-species Keller–Segel model is energy dissipative, i.e.,
	\begin{equation}\label{model:KS:Energy}
		\begin{aligned}
			\f{d\mrE[u,v,c]}{dt} = -\int_\Omega \left[ u \Bigl(\nabla \f{\delta \mrE}{\delta u}\Bigr)^2 
            +  v \Bigl(\nabla \f{\delta \mrE}{\delta v} \Bigr)^2 
            + \epsilon\Bigl(\f{\p c}{\p t}\Bigr)^2 \right] d\bm{x} 
            \le 0,
		\end{aligned}
	\end{equation}
    where $\f{\delta \mrE}{\delta u} =  f'(u) -  c$ and	$\f{\delta \mrE}{\delta v} =  f'(v) - c$ are the variational derivatives. In summary, the two-species Keller–Segel model preserves three fundamental physical properties: positivity of the cell densities and chemoattractant concentration, mass conservation for each cell species, and dissipation of the original energy. Therefore, a desirable numerical method should preserve these properties simultaneously at the discrete level.	
	
   In recent years, considerable efforts have been devoted to developing structure-preserving numerical methods for the Keller--Segel model, with particular emphasis on the preservation of some of the fundamental physical properties \eqref{model:KS:MC}, \eqref{model:KS:positivity}, and/or \eqref{model:KS:Energy}. For example, %Xiao et al. \cite{XFH'19} developed a characteristic finite element method for the  surface chemotaxis model that can preserve positivity and mass conservation; however, the optimal-order error analysis is lacking.
    Chertock et al. \cite{AYHA'18} proposed a fourth-order hybrid finite-volume–finite-difference scheme by incorporating a positivity-preserving adaptive flux-rescaling strategy based on local draining time steps, thereby achieving high-order spatial accuracy while maintaining positivity preservation and mass conservation. Cheng and Shen \cite{CS'22} introduced a Lagrange multiplier approach that preserves positivity/bounds and mass conservation without resorting to nonlinear constrained optimization, whereas function transformation-based methods, sometimes combined with the scalar auxiliary variable approach, have been employed to achieve positivity/bound preservation, mass conservation, and modified energy stability \cite{HS'21}. %Kemmochi at al. \cite{KS'22} proposed an energy-decomposition-based scalar auxiliary variable approach to construct a second-order linearly implicit Crank--Nicolson scheme and a fourth-order Runge--Kutta scheme for conservative and dissipative partial differential equations with unbounded energy functionals, which satisfy a modified energy stability property. 
    Recently, Wang et al. \cite{WLF'25} developed a first-order fully discrete finite element scheme that combines an exponential transformation to unconditionally preserve positivity with a normalization procedure to enforce mass conservation. In particular, a rigorous error analysis was carried out for the proposed scheme. However, only an energy bound, rather than energy dissipation, was established. In addition, Yin et al. \cite{YLQ'26} used the direct discontinuous Galerkin method to develop a semi-implicit, decoupled scheme that conserves mass, preserves positivity, and is conditionally energy stable. Optimal-order error estimates were established. However, the resulting scheme requires solving nonlinear algebraic systems at each time step, which may increase the computational cost.  A very recent work by Huang et al. \cite{HGS'26} studied a semi-implicit Euler scheme widely used for the time discretization of Keller–Segel equations in both parabolic-elliptic and parabolic-parabolic forms. They proved that this linearly decoupled, first-order scheme preserves the three key physical properties specified in \eqref{model:KS:MC}, \eqref{model:KS:positivity}, and \eqref{model:KS:Energy} at the semi-discrete level and, more importantly, established optimal-order error estimates in the $L^p$-norm ($1<p <\infty$) for this scheme. They further noted that constructing a fully discrete, linearly decoupled scheme that unconditionally preserves all essential properties while still permitting a rigorous convergence analysis remains highly challenging! %Motivated by this observation, we develop such a scheme, which simultaneously preserves the three physical properties specified in \eqref{model:KS:MC}, \eqref{model:KS:positivity}, and \eqref{model:KS:Energy} at the fully discrete level, and we establish rigorous optimal-order error estimates for the resulting scheme.

    In addition to the aforementioned approaches, Slotboom transformation methods \cite{Slotboom'73, JY'11, HZ'23,LCLL'24} have attracted increasing attention in recent years as an effective framework for developing efficient structure-preserving numerical schemes. For example, Liu et al. \cite{LWZ'18} constructed a second-order central difference scheme based on the Slotboom transformation that can preserve both positivity and mass conservation. Although optimal-order error estimates were not established, the authors provided a stability analysis for the proposed scheme. In addition, a number of Slotboom transformation-based structure-preserving numerical schemes have been developed for the Poisson-Nernst-Planck system (see Refs. \cite{LW'14,HH'20,DWZ'19}). However, to the best of our knowledge, a rigorous convergence analysis for theses schemes is still lacking and challenging, particularly owing to the strong nonlinear coupling between the equations in these coupled systems. 

  % High-order consistency analysis, which occurs severe time-step condition caused in establishing a priori $\ell^{\infty}$ bound for the numerical solution, has been employed as an effective tool for establishing optimal-order error estimates of numerical schemes for nonlinear equations. For example, Liu et al. \cite{LWWQS'21} developed a nonlinear scheme based on the Wasserstein gradient flows for the Poisson--Nernst--Planck system that preserves positivity and energy dissipation, and established rigorous error estimates through a high-order consistency analysis. Furthermore, the high-order consistency analysis technique has been subsequently employed in the development and convergence analysis of second-order structure-preserving schemes for the Keller--Segel and Poisson--Nernst--Planck systems \cite{WWQZ'26,LWWQS'23,DWZ'26}. However, the above numerical schemes are all nonlinear, which presents a great challenge to computational efficiency. 

  Motivated by this observation, we invoke the Slotboom transformation to reformulate the two-species Keller--Segel chemotaxis system \eqref{model:KS} into the following symmetric form:
    \begin{equation}\label{model:KS_transform}
    \left\{
		\begin{aligned}
			&  \p_t u  = \nabla \cdot(e^c \nabla ({e^{-c}}u)), & \qquad \text{in}~~ \Omega \times(0, T],\\
			&  \p_t v  = \nabla \cdot(e^c \nabla ({e^{-c}}v)), & \qquad \text{in}~~ \Omega \times(0, T],\\            
			& \epsilon\,\p_t c  =  \beta\Delta c- \alpha c +  u  +  v, & \qquad \text{in}~~ \Omega \times(0, T].
		\end{aligned}
        \right.
	\end{equation}
Starting from the symmetric form and using the widely used block-centered finite difference (BCFD) method \cite{AWY'97,RP'12,RL'15,XXF'22,WXF'24,LSR'19,XF'25}, which is mass-conservative in nature and flexible in treating Neumann or periodic boundary conditions, we develop a linearly decoupled BCFD scheme that simultaneously preserves all three essential physical properties specified in \eqref{model:KS:MC}, \eqref{model:KS:positivity}, and \eqref{model:KS:Energy} at the fully discrete level. We then establish rigorous optimal-order error estimates via a high-order consistency analysis approach \cite{LWWQS'21,WWQZ'26,LWWQS'23,DWZ'26}, which has been employed as an effective tool for establishing optimal-order error estimates of numerical schemes for nonlinear equations. In summary, the main contributions of the present work are as follows:
	\begin{itemize}
		\item The proposed scheme is fully decoupled and linear, which significantly improves computational efficiency.
        
		\item The proposed scheme unconditionally preserves, at the fully discrete level, the positivity of the cell densities and chemoattractant concentration, as well as the mass conservation of the cell densities. In addition, it ensures the unconditional dissipation of the original energy and the unique solvability of the numerical solutions.  
        
		\item Optimal-order error estimates for the scheme are rigorously established under a mild time-step condition using the high-order consistency analysis.
	\end{itemize}
	
	The rest of the paper is organized as follows. In Section \ref{sec:mc-bcfd}, based on the Slotboom transformation, we propose the structure-preserving BCFD scheme for the two-species Keller--Segel chemotaxis model and prove that it preserves positivity and mass conservation, dissipates the original energy, and is uniquely solvable. Optimal-order error estimates are presented in Section \ref{sec:erro-estimates}. In Section \ref{sec:num}, several numerical experiments are shown to validate the accuracy and physical-property-preserving performance of the proposed scheme. Some conclusions are drawn in Section \ref{sec:conclusion}.  Throughout this paper,  $Q$, with or without subscripts or superscripts, denotes a generic positive constant that is independent of the grid parameters but may have different values in different occurrences.

%%%%%%%%%%%%%%%%%%%%%%%%%%%%%%%%%%%%%%%%%%%%%%%%%%%%%%%%%%%%%%%
\section{Structure-preserving BCFD scheme}\label{sec:mc-bcfd}
%%%%%%%%%%%%%%%%%%%%%%%%%%%%%%%%%%%%%%%%%%%%%%%%%%%%%%%%%%%%%%%
This section is devoted to the construction of a fully discrete, linearly decoupled, and structure-preserving BCFD method for the two-species Keller--Segel chemotaxis system. %The proposed method is positivity-preserving and mass-conserving, and it satisfies an original energy dissipation law. 
For simplicity, we take the domain to be a two-dimensional rectangular $\Omega :=(a, b)^2$. 
    
	%%%%%%%%%%%%%%%%%%%%%%%%%%%%%%%%%%%%%%%%%%%%%%%%%%%%%%%%%%%%%%%
	\subsection{Notations and preliminaries}
	%%%%%%%%%%%%%%%%%%%%%%%%%%%%%%%%%%%%%%%%%%%%%%%%%%%%%%%%%%%%%%%
	Let $N$ and $M$ be positive integers. The time levels are defined by $t_{n}=n\tau$ ($n=0,1,\ldots,N$) with uniform temporal stepsize $\tau=T/N$.  In addition, we introduce two uniform spatial grids along the $x$- and $y$-coordinates, respectively:
    \begin{align*}
        \Pi_z: ~ a=z_{1 / 2}<z_{3 / 2}< \ldots < z_{i-1 / 2}<z_{i+1 / 2} < \ldots <  z_{M+1 / 2}=b,~z=x,y,
    \end{align*}
with equal grid size $h=z_{i+1/2}-z_{i-1/2}=(b-a)/M$ for $i=1, \ldots, M$. Moreover, the staggered middle grid points are denoted by
\begin{align*}
	 \Pi_x^*: ~ x_i= (x_{i-1/2}+x_{i+1/2} ) / 2,  ~~ i=1, \ldots, M,~~
     \Pi_y^*: ~ y_j= (y_{j-1/2}+ y_{j+1/2} ) / 2,  ~~ j=1, \ldots, M.
\end{align*}

Given spatial grid functions $v=\{v_{i, j}\}$, $w^x=\{w^x_{i+1/2, j}\}$, and $w^y=\{w^y_{i, j+1/2}\}$ that defined on $\Pi_x^*\times \Pi_y^*$, $\Pi_x\times \Pi_y^*$, and $\Pi_x^*\times \Pi_y$, respectively, we define
\begin{align*}
	& [d_x v]_{i+1/2, j}=d_x v_{i+1/2, j}=\f{v_{i+1, j}-v_{i, j}}{h}, ~~	
        && [d_y v]_{i, j+1/2}=d_y v_{i, j+1/2}=\f{v_{i, j+1}-v_{i, j}} {h}, \\
	& \bigl[D_x w^x]_{i, j}=D_x w^x_{i, j}=\f{w^x_{i+1/2, j}-w^x_{i-1/2, j}}{h}, ~~
	&&  \bigl[D_y w^y]_{i, j}=D_y w^y_{i, j}=\f{w^y_{i, j+1/2}-w^y_{i, j-1/2}}{h},\\
    & [\mal{A}_x v]_{i+1/2,j} =\mal{A}_x v_{i+1/2,j}=\frac{v_{i,j} + v_{i+1,j}}{2}, ~~
	&&  [\mal{A}_y v]_{i,j+1/2}=\mal{A}_y v_{i,j+1/2}=\frac{v_{i,j} + v_{i,j+1}}{2},
\end{align*}
and denote $\mal{A} v_{i+1/2,j+1/2}:=(\mal{A}_x v_{i+1/2,j},\mal{A}_y v_{i,j+1/2})$, $\bm{d}v_{i+1/2,j+1/2}:=(d_x v_{i+1/2,j},d_y v_{i,j+1/2})$. Besides, we introduce the discrete inner products and norms on $\Pi_x^*\times \Pi_y^*$, $\Pi_x\times \Pi_y^*$, and $\Pi_x^*\times \Pi_y$, respectively, as follows:
\begin{align*}
	& (f, g)_{\rm M}=h^2\sum_{i=1}^{M} \sum_{j=1}^{M} f_{i, j} g_{i, j},  &&\|f\|_{\rm M}^2 = (f, f)_{\rm M}, \\
	& (f, g)_x=h^2\sum_{i=0}^{M-1} \sum_{j=1}^{M} f_{i+1/2, j} g_{i+1/2, j},  &&\|f\|_x^2 = (f, f)_x, \\
	& (f, g)_y=h^2\sum_{i=1}^{M} \sum_{j=0}^{M-1} f_{i, j+1/2} g_{i, j+1/2},  &&\|f\|_y^2 = (f, f)_y,\\
	& (\bm {f}, \bm {g})_{\rm TM}  = (f^x, g^x )_x+ (f^y, g^y)_y,    &&\|\bm {f}\|_{\rm TM}^2 = (\bm{f}, \bm{f})_{\rm TM},\\
    %& \|d_x g\|_\infty =\max _{i, j} |[d_x g]_{i+1/2, j}|,~\|d_y g\|_\infty =\max _{i, j} |[d_y g]_{i, j+1/2}|, &&
    & \|\bm{d}g\|_\infty %= \|d_x g\|_\infty +\|d_y g\|_\infty
    = \max _{i, j} |[d_x g]_{i+1/2, j}|+ \max _{i, j} |[d_y g]_{i, j+1/2}|. &&
\end{align*}

 Next, satrting from the equivalent two-species Keller–Segel chemotaxis system \eqref{model:KS_transform}, we shall propose a linearly decoupled, positivity-preserving, mass-conserving and original energy dissipation Euler-BCFD scheme. 
	
%%%%%%%%%%%%%%%%%%%%%%%%%%%%%%%%%%%%%%%%%%%%%%%%%%%%%%%%%%%%%%%
\subsection{Numerical scheme}
%%%%%%%%%%%%%%%%%%%%%%%%%%%%%%%%%%%%%%%%%%%%%%%%%%%%%%%%%%%%
Let $\mbX_h$ be the set of real-valued periodic grid functions defined on $\Pi_x^*\times \Pi_y^*$, i.e., 
    \begin{equation*}
        \mbX_h:=\{q_h \mid q_h=\{q_{i,j}\},~q_{0,j}= q_{M,j},~q_{1,j}= q_{M+1,j},~q_{i,0}= q_{i,M},~q_{i,1}= q_{i,M+1},~ i,j=1,\ldots,M\}.
    \end{equation*}
Each function $q_h\in \mbX_h$ can be viewed as an $M$-by-$M$ matrix, and its boldface counterpart $\bm{q}_h=vec\,{(q_h)}$ is the corresponding $M^2$-dimensional column vector, with the $x$-label arranged first and the $y$-label last.

Given $\{u_h^0,v_h^0,c_h^0\}\in \mbX_h\times\mbX_h\times\mbX_h$, for $n\geq0$, the positivity-preserving, mass-conserving and original-energy-dissipating BCFD scheme, referred to as PP-MC-ED-BCFD, is proposed as follows: 

    \medskip
	{\bf Step 1}:  Solve the cell density solutions $\{u_h^{n+1},v_h^{n+1}\}\in \mbX_h \times \mbX_h$ via the IMEX-BCFD scheme: 
	\begin{equation}\label{model:KS:shemeI}
		\begin{aligned}
			\f{u_h^{n+1}-u_h^n}{\tau} &= D_x \Big([\mal{A}_x e^{c_h^n}] d_x [e^{-c_h^n}u_h^{n+1}]\Big)
            +D_y \Big([\mal{A}_y e^{c_h^n}] d_y [e^{-c_h^n}u_h^{n+1}]\Big),\\
			\f{v_h^{n+1}-v_h^n}{\tau} &= D_x \Big([\mal{A}_x e^{c_h^n}] d_x [e^{-c_h^n}v_h^{n+1}]\Big) + D_y \Big([\mal{A}_y e^{c_h^n}] d_y [e^{-c_h^n}v_h^{n+1}]\Big).
		\end{aligned}
	\end{equation}

	{\bf Step 2}: Solve the chemoattractant concentration solution $c_h^{n+1}\in \mbX_h$ via the implicit-BCFD scheme:
	\begin{equation}\label{model:KS:shemeII}
	  \epsilon\f{c_h^{n+1}-c_h^n}{\tau} 
      = \beta D_x d_x c_h^{n+1} + \beta D_y d_y c_h^{n+1} - \alpha c_h^{n+1} + u_h^{n+1} 
        + v_h^{n+1}.
	\end{equation}
	
For the purpose of analysis, we review some well-known results in the following lemmas.
\begin{lemma}[\cite{MW'09}] \label{lemma:Dd} 
 Let $\{q_{i, j}\}, \{v_{i+1/2, j}\}$ and $\{w_{i, j+1 / 2}\}$ be any grid functions defined on $\Pi_x^*\times \Pi_y^*$, $\Pi_x\times \Pi_y^*$ and $\Pi_x^*\times \Pi_y$, such that $v_{1/2, j}=v_{M+1/2, j}$, $w_{i, 1/2}=w_{i, M+1 / 2}$ and $q_{i, j}=q_{i+M, j+M}$. Then there holds
$$
\begin{aligned}
 (q, D_x v)_{\rm M}=-(d_x q, v)_x, ~~ (q, D_y w)_{\rm M}=-(d_y q, w)_y .
\end{aligned}
$$
\end{lemma}
\begin{lemma}[\cite{HZ'23}]\label{lem:M-matrix}
	For a real square matrix $\mbb{A}$ with positive diagonal entries and non-positive off-diagonal entries, it is a non-singular M-matrix if and only if there exists a positive diagonal matrix $\mbb{D}$ such that $\mbb{AD}$ has all positive row sums.
\end{lemma}
\begin{lemma}[\cite{Ple'77}]\label{lem:M-matrix:inverse}
	Let $\mbb{A}\in \mbb{R}^{n\times n}$ be a non-singular M-matrix, then $\mbb{A}$ is inverse-positive. That is, $\mbb{A}^{-1}$ exists and $\mbb{A}^{-1} \geq 0.$
\end{lemma}
The following theorem presents the main structure-preserving properties and unique solvability of the proposed PP-MC-ED-BCFD scheme, namely, unconditional positivity preservation, mass conservation, and original energy dissipation.  
	\begin{theorem}\label{thm:structure-pre}
		For the PP-MC-ED-BCFD scheme, there holds
		\begin{itemize}
			\item[(i)] (Positivity preservation) If the initial cell densities $u_h^{0}, v_h^{0} > 0$, then $u_h^{n+1},v_h^{n+1}>0$. Furthermore, if the initial chemoattractant concentration $c_h^{0} > 0$, then we also have $c_h^{n+1} > 0$. 
			\item[(ii)] (Mass conservation) The density solutions satisfy 
			$$
			(u_h^{n+1}, 1)_{\rm M} = ( u_h^0, 1 )_{\rm M}, ~~
			( v_h^{n+1}, 1)_{\rm M} = ( v_h^0, 1)_{\rm M}.
			$$
			\item [(iii)]  (Original energy dissipation) Let $\mal{E}_h^n := ( u_h^n \log \frac{u_h^n}{e^{c_h^n}} +v_h^n \log \frac{v_h^n}{e^{c_h^n}} - u_h^n - v_h^n, 1 )_{\rm M} + \frac{\alpha}{2}\|c_h^n\|^2_{\rm M}+\frac{\beta}{2} \|\bm{d}c_h^n\|^2_{\rm TM}$ be the discrete version of the free energy \eqref{def:Energy}. Then,
			$$
				\mal{E}_h^{n+1} \leq \mal{E}_h^n.
			$$
			\item[(iv)] (Unique solvability) There exist unique solutions $\{u_h^{n+1},v_h^{n+1},c_h^{n+1}\}\in \mbX_h \times \mbX_h \times \mbX_h$ to the scheme \eqref{model:KS:shemeI}--\eqref{model:KS:shemeII}.
		\end{itemize}
	\end{theorem}
	\begin{proof}
		% Since \eqref{model:KS:shemeI}--\eqref{model:KS:shemeII} are linear, thus the solutions at each step are uniquely solvable. 
       First, taking the inner product of \eqref{model:KS:shemeI} with 1 in the sense of $(\cdot, \cdot)_{\rm M}$, and utilizing Lemma \ref{lemma:Dd},  the mass-conservation of the cell density solutions $u_h^{n+1}$ and $v_h^{n+1}$ are easily obtained. 
       
    Second, we first show the positivity of the values $\{u_h^{n+1},v_h^{n+1}\}$ for positive inputs $\{u_h^{n},v_h^{n}\}$. For simplicity, denote $\mal{M}_h^n := e^{c_h^n}$. Suppose that there exist a pair of index $(r,s)$ such that
		\begin{equation*}
			% \min_{\substack{1 \leq i \leq M \\ 1 \leq j \leq M}} \frac{u^{n+1}_{i,j}}{e^{c^n_{i,j}}} 
   %          = \frac{u^{n+1}_{r,s}}{e^{c^n_{r,s}}},
           \min_{\substack{1 \leq i,j \leq M}} \frac{u^{n+1}_{i,j}}{\mal{M}_{i,j}^n} 
            = \frac{u^{n+1}_{r,s}}{\mal{M}^n_{r,s}}, 
		\end{equation*}
		Then, choosing the index $(r,s)$ in \eqref{model:KS:shemeI}, we obtain
\begin{equation*}
    \begin{aligned}
        u_{r,s}^{n+1} &= u_{r,s}^n + \frac{\tau}{h^2}
        \left[ \frac{\mal{M}_{r+1,s}^n + \mal{M}_{r,s}^n}{2}
        \Bigl(\frac{u_{r+1,s}^{n+1}}{\mal{M}_{r+1,s}^n}
        - \frac{u_{r,s}^{n+1}}{\mal{M}_{r,s}^n}\Bigr)
        - \frac{\mal{M}_{r,s}^n + \mal{M}_{r-1,s}^n}{2}
        \Bigl(\frac{u_{r,s}^{n+1}}{\mal{M}_{r,s}^n}
        - \frac{u_{r-1,s}^{n+1}}{\mal{M}_{r-1,s}^n}\Bigr)\right] \\
        &\quad + \frac{\tau}{h^2}
        \left[\frac{\mal{M}_{r,s+1}^n + \mal{M}_{r,s}^n}{2}
        \Bigl(\frac{u_{r,s+1}^{n+1}}{\mal{M}_{r,s+1}^n}
        - \frac{u_{r,s}^{n+1}}{\mal{M}_{r,s}^n}\Bigr)
        - \frac{\mal{M}_{r,s}^n + \mal{M}_{r,s-1}^n}{2}
        \Bigl(\frac{u_{r,s}^{n+1}}{\mal{M}_{r,s}^n}
        - \frac{u_{r,s-1}^{n+1}}{\mal{M}_{r,s-1}^n}\Bigr)\right] \\
        &\geq u_{r,s}^n > 0,
    \end{aligned}
\end{equation*}
        which further implies $\frac{u^{n+1}_{r,s}}{\mal{M}^n_{r,s}} > 0$, and thus, $u_{i,j}^{n+1} > 0$ due to $\mal{M}_{i,j}^n>0$. Thus, the proof of positivity preservation for $u$ is completed. Similarly, the positivity preservation for $v$ can be obtained in an analogous manner.

		Furthermore, assuming $c_h^{n} > 0$, we show the positivity of $c_h^{n+1}$.  Let $\mbb{B}$ of order $M$ be a tridiagonal circulant matrix whose first row is
        $\f{\beta}{h^2} [-2, 1, 0,\ldots, 0, 1]$, and $\mbb{I}$ be the corresponding identity matrix. Then, we can rewrite \eqref{model:KS:shemeII} in a compact matrix form:
		\begin{equation*}
			\mbb{P} \bm{c}_h^{n+1} = \mbb{Q} \bm{c}_h^{n}+ \bm{u}_h^{n+1} + \bm{v}_h^{n+1},
		\end{equation*}
		where $\mbb{P}:=(\f{\epsilon}{\tau} + \alpha) \mbb{I}\otimes \mbb{I} - \mbb{B} \otimes \mbb{I}- \mbb{I}  \otimes \mbb{B}$ and $\mbb{Q}:= \f{\epsilon}{\tau} \mbb{I}\otimes \mbb{I}$.  
		From Lemma \ref{lem:M-matrix}, we can easily check that $\mbb{P}$ is an M-matrix and thus Lemma \ref{lem:M-matrix:inverse} implies that  $\mbb{P}^{-1} \geq 0$. Note that $\mbb{Q}$ is positive and $\bm{c}_h^{n},\, \bm{u}_h^{n+1}\, \bm{v}_h^{n+1} > 0 $. Thus $\bm{c}_h^{n+1} =\mbb{P}^{-1}(\mbb{Q} \bm{c}_h^{n}+ \bm{u}_h^{n+1} + \bm{v}_h^{n+1}) > 0$ is proved.
		
%		The detailed proof of energy dissipation will be provided in \textbf{Appendix}.
		Third, the original energy dissipation property will be proved directly. By definition and conclusion (ii), we have
		\begin{equation}\label{Discrete_Energy:e1}
			\begin{aligned}
				\mal{E}_h^{n+1} - \mal{E}_h^n 
            &= \left[\bigl ( u_h^{n+1}\log h_u^{n+1} ,1\bigr )_{\rm M} 
                - \bigl( u_h^{n}\log h_u^n ,1\bigr )_{\rm M}\right]\\
				& \quad
				+ \left[\bigl( v_h^{n+1}\log h_v^{n+1}, 1\bigr )_{\rm M} 
                - \bigl( v_h^{n}\log h_v^n ,1\bigr)_{\rm M}\right]\\
				& \quad + \frac{\alpha}{2} \bigl[\|c_h^{n+1}\|^2_{\rm M}
                - \|c_h^n\|^2_{\rm M}\bigr]+ \frac{\beta}{2}\bigl[\|\bm{d}c_h^{n+1}\|^2_{\rm TM}
                - \|\bm{d}c_h^n\|^2_{\rm TM}\bigr]
				=: I_1 + I_2 + I_3,
			\end{aligned}
		\end{equation}
where $g_u^{n+1} := \frac{u_h^{n+1}}{\mal{M}_h^n}$, $h_u^n := \frac{u_h^n}{\mal{M}_h^n}$, $g_v^{n+1} := \frac{v_h^{n+1}}{\mal{M}_h^n}$, and $h_v^n := \frac{v_h^n}{\mal{M}_h^n}$.    
In the following, we will analyze each term of \eqref{Discrete_Energy:e1}. 
For the first term $I_1$, a simple calculation yields
	\begin{equation*}
		\begin{aligned}
		I_1
        &= \bigl( \mal{M}_h^n g_u^{n+1}\log g_u^{n+1}  ,1\bigr )_{\rm M} 
                - \bigl(\mal{M}_h^n h_u^{n}\log h_u^n ,1\bigr)_{\rm M} + \bigl(u_h^{n+1},c_h^n-c_h^{n+1}\bigr)_{\rm M}\\
        &= \bigl( \mal{M}_h^n \check{g}_u^{n+1}  ,1\bigr )_{\rm M} 
                - \bigl(\mal{M}_h^n \check{h}_u^{n} ,1\bigr)_{\rm M} 
                + \bigl(u_h^{n+1},c_h^n-c_h^{n+1}\bigr)_{\rm M},
                \end{aligned}
	\end{equation*}
    where $\check{g}_u^{n+1} :=g_u^{n+1}\log g_u^{n+1}$ and $\check{h}_u^n:=h_u^n \log h_u^n$.
Let $\bm{1}$ be an $M^2$-dimensional column vector with all entries being one. By definition of the discrete inner product, we see
	\begin{equation}\label{eq:gn1}
		\bigl (\mal{M}_h^n \check{g}_u^{n+1} , 1\bigr)_{\rm M} 
        = h^2 \sum_{i=1}^{M}\sum_{j=1}^{M} \mal{M}_{i,j}^n (\check{g}_u^{n+1})_{i,j} = h^2\, \bm{1}^{\top} \mbb{M}_d^n \,\check{\bm{g}}_u^{n+1},
	\end{equation}
where $\mbb{M}_d^n:=\operatorname{diag}(\mal{M}_h^n)$ is an $M^2$-by-$M^2$ diagonal matrix. 
Recall that the first equation of \eqref{model:KS:shemeI} is equivalent to
	\begin{equation*}
				\f{\mal{M}_h^n g_u^{n+1}- \mal{M}_h^n h_u^n}{\tau}  = D_x  \big([\mal{A}_x \mal{M}_h^n] [d_x g_u^{n+1}]\big)+D_y  \big([\mal{A}_y \mal{M}_h^n] [d_y g_u^{n+1}]\big),
				% \f{\mal{M}_h^n g_v^{n+1}- \mal{M}_h^n h_v^n}{\tau} & = D_x  \big([\ell_h \mal{M}_h^n] [d_x g_v^{n+1}]\big)+D_y  \big([\ell_h \mal{M}_h^n] [d_y g_v^{n+1}]\big).\label{model:KS:re2}
	\end{equation*}
or reformulated into the compact matrix form
		\begin{equation}\label{model:KS:re1_eq}
			\mbb{A}_u \bm{g}_u^{n+1} = \mbb{M}_d^n \bm{h}_u^n 
            \Longrightarrow 
            \widetilde{\mbb{A}}_u \bm{g}_u^{n+1} = \bm{h}_u^n,
		\end{equation}
where $\mbb{A}_u := \mbb{M}_d^n + \mbb{R}$ with $\mbb{R}=\mbb{R}^{\top}$, $\mbb{R}\bm{1}=\bm{0}$ and $\widetilde{\mbb{A}}_u = \mbb{I} + (\mbb{M}_d^n)^{-1}\mbb{R}$. 

By Lemma \ref{lem:M-matrix}, $\widetilde{\mbb{A}}_u$ is a non-singular M-matrix, and therefore, $\widetilde{\mbb{A}}_u^{-1}\geq 0$ and $\widetilde{\mbb{A}}_u^{-1} \bm{1} = \bm{1}$. Thus, it follows from \eqref{model:KS:re1_eq} that $\bm{g}_u^{n+1}=\widetilde{\mbb{A}}_u^{-1}\bm{h}_u^n$ defines a row convex combination of $\bm{h}_u^n$. Furthermore, note that $x\log x$ is convex in $(0, \infty)$. Then, Jensen's inequality implies that
	\begin{equation}\label{eq:Jensen}
		\check{\bm{g}}_u^{n+1} \leq \widetilde{\mbb{A}}_u^{-1} \check{\bm{h}}_u^n.
	\end{equation}
  Moreover,
	\begin{equation}\label{model:KS:re2_eq}
		 \bm{1}^{\top} \mbb{M}_d^n \widetilde{\mbb{A}}_u =  \bm{1}^{\top} \mbb{M}_d^n\, (\mbb{I} + (\mbb{M}_d^n)^{-1}\mbb{R}) =  \bm{1}^{\top} \mbb{M}_d^n \Longrightarrow
		 \bm{1}^{\top} \mbb{M}_d^n =  \bm{1}^{\top} \mbb{M}_d^n \widetilde{\mbb{A}}_u^{-1}.
	\end{equation}
Therefore, inserting \eqref{eq:Jensen}--\eqref{model:KS:re2_eq} into \eqref{eq:gn1}, we have
	\begin{equation*} 
		\bigl( \mal{M}_h^n \check{g}_u^{n+1}, 1\bigr)_{\rm M} \leq h^2\, \bm{1}^{\top} \mbb{M}_d^n \widetilde{\mbb{A}}_u^{-1} \check{\bm{h}}_u^n = h^2\, \bm{1}^{\top} \mbb{M}_d^n \check{\bm{h}}_u^n = \bigl( \mal{M}_h^n \check{h}_u^{n}  , 1\bigr)_{\rm M},
	\end{equation*}
and thus, the $I_1$-term can be estimated by
	\begin{equation*}
	I_1 \leq \bigl( u_h^{n+1},c_h^n-c_h^{n+1}\bigr)_{\rm M}.
	\end{equation*} 
By the same argument, the estimate of $I_2$ is obtained as
\begin{equation*}
	I_2 \leq \bigl( v_h^{n+1},c_h^n-c_h^{n+1}\bigr)_{\rm M}.
\end{equation*}

Furthermore, recall that $ \frac{1}{2}(a^2 - b^2) \le (a-b) a$,
% using the following identity
% \begin{equation*}
% 	(a-b) a = \frac{1}{2}\left(a^2 - b^2 + (a-b)^2\right),
% \end{equation*}
the $I_3$-term can be bounded as follows:
\begin{equation*}
	\begin{aligned}
		I_3 &=  \frac{\alpha}{2}\|c_h^{n+1}\|^2_{\rm M}-\frac{\alpha}{2}\|c_h^n\|^2_{\rm M}+ \frac{\beta}{2} \|\bm{d}c_h^{n+1}\|^2_{\rm TM}-\frac{\beta}{2} \|\bm{d}c_h^n\|^2_{\rm TM}\\
		&\leq \alpha\bigl( c_h^{n+1} - c_h^n,c_h^{n+1}\bigr)_{\rm M} 
        + \beta\bigl( d_x (c_h^{n+1} - c_h^n),d_x c_h^{n+1}\bigr)_x 
        + \beta\bigl( d_y (c_h^{n+1} - c_h^n),d_y c_h^{n+1}\bigr)_y\\
		&\leq \alpha \bigl( c_h^{n+1} - c_h^n, c_h^{n+1}\bigr)_{\rm M} 
        - \beta \bigl( c_h^{n+1} - c_h^n, D_xd_x c_h^{n+1} +  D_yd_y c_h^{n+1}\bigr)_{\rm M}.
	\end{aligned}
\end{equation*}
Now, inserting the estimates $I_1$--$I_3$ into \eqref{Discrete_Energy:e1} and using \eqref{model:KS:shemeII}, the energy estimate follows
\begin{equation*}
	\begin{aligned}
		\mal{E}_h^{n+1} - \mal{E}_h^n &\leq \bigl( c_h^{n+1} - c_h^n,-\beta D_xd_x c_h^{n+1}- \beta D_yd_y c_h^{n+1} + \alpha c_h^{n+1} - u_h^{n+1} -v_h^{n+1}\bigr)_{\rm M}\\
		& = -\f{\epsilon}{\tau} \|c_h^{n+1}-c_h^n\|^2 \leq 0.
	\end{aligned}
\end{equation*}

Finally, as mentioned earlier, the matrices $\mbb{P}^{-1}$, $\widetilde{\mbb{A}}_u^{-1}$ and $\widetilde{\mbb{A}}_v^{-1}$ are positive and hence non-singular, which directly implies that the existence and uniqueness of the solution to the scheme \eqref{model:KS:shemeI}--\eqref{model:KS:shemeII}. This completes the proof of Theorem \ref{thm:structure-pre}.
	\end{proof}
	
	\begin{remark}
		In addition to the arithmetic mean values $\mal{A}_x e^{c_h^n}= \f{e^{c^n_{i,j}}+ e^{c^n_{i+1,j}}}{2}$ and $\mal{A}_y e^{c_h^n}= \f{e^{c^n_{i,j}}+ e^{c^n_{i,j+1}}}{2}$ used in \eqref{model:KS:shemeI}, the midpoint value of $e^{c_h^n}$ can also be approximated by the geometric or harmonic mean. For example, the  geometric mean in the $x$-direction is defined as $[\mal{G}_x e^{c_h^n}]_{i+1/2,j} := \sqrt{e^{c^n_{i,j}}e^{c^n_{i+1,j}}}$, and the corresponding harmonic mean is $[\mal{H}_x e^{c_h^n}]_{i+1/2,j} := \f{2e^{c^n_{i,j}}e^{c^n_{i+1,j}}}{e^{c^n_{i,j}}+e^{c^n_{i+1,j}}}$. Regardless of which mean is used, the corresponding fully discrete scheme can be shown to unconditionally preserve positivity, mass conservation, and the original energy dissipation.
	\end{remark}
		
		%%%%%%%%%%%%%%%%%%%%%%%%%%%%%%%%%%%%%%%%%%%%%%%%%%%%%%%%%%%%%%%
		\section{Error estimates}\label{sec:erro-estimates}
		%%%%%%%%%%%%%%%%%%%%%%%%%%%%%%%%%%%%%%%%%%%%%%%%%%%%%%%%%%%%%%%
		In this section, we present a rigorous higher-order consistency error analysis for the proposed scheme \eqref{model:KS:shemeI}--\eqref{model:KS:shemeII}, showing that it is first-order accurate in time and second-order accurate in space. 
        
        Let $\{u,v,c\}$ denote the exact solutions to the two-species Keller--Segel chemotaxis system, subject to  the following regularity assumptions:
		\begin{equation}\label{eq:regular}
			u(\bm x,t),v(\bm x,t),c(\bm x,t) \in H^{4}(0,T;C(\Omega)) \cap H^{3}(0,T;C^2(\Omega)) \cap L^{\infty}(0,T;C^6(\Omega)).
		\end{equation}
		Furthermore, for the convenience of analysis, we assume the following separation property:
		\begin{equation}\label{bound:exc1}
			c(\bm x,t) \geq \varepsilon>0,
		\end{equation}
		which holds at the point-wise level. 
        
        Let $\mal{P}_M$ denote the standard Fourier projection operator from $C_{per}(\Omega)$, the space of continuous $\Omega$-periodic functions, onto $\mal{B}^K$, the space of $\Omega$-periodic trigonometric polynomials of degree at most $K$ (with $M = 2K+1$),  and denote
        \begin{equation}\label{def:projection}
        z_M(\bm x,t):=\mal{P}_M z(\bm x,t),\quad \text{for} ~z=u,v,c.
        \end{equation}
  	The following projection approximation holds \cite{LWWQS'21}: 
		\begin{equation}\label{projection:error}
				\|z - z_M\|_{H^{k}} \leq Q h^{\ell-k}\|z\|_{H^{\ell}},
		\end{equation}
        for $z=\{u,v,c\}\in L^{\infty}(0,T;H^{\ell})$ with $0\leq k \leq \ell\le 6$.        
		Clearly, the projection \eqref{def:projection} does not preserve the positivity of $c_M(\bm x,t)$. Nevertheless, provided that $h$ is chosen sufficiently small (equivalently, that $M$ is sufficiently large), it follows from the approximation error estimate \eqref{projection:error} that 
        \begin{equation}\label{bound:exc1a}
         c_M(\bm x,t) \geq \frac{1}{2} \varepsilon.
         \end{equation}
         % Furthermore, let $\Pi_{h}$ denote the canonical grid operator such that $\Pi_{h} z \in \mbb{X}_h$ with $[\Pi_{h} z]_{i,j}=z(x_i,y_j)$ for $i,j=0,1,\ldots,M$. 
 %         For each $0\le n \le N$, define the grid error functions
 %        \begin{equation}\label{error:1}
	% 		%e_z^n := 
 %            z^n -  z_h^n= (z^n-\Pi_{h} z_M^n) + (\Pi_{h} z_M^n-z_h^n) 
 %                    =:\eta_z^n + \theta_z^n, \quad \text{for} ~ z=u,v,c,
 %        \end{equation}
	% with the initial conditions chosen as 
 %        \begin{equation*}
	% 		u^0_h = \Pi_{h}\mal{P}_M u^0,~~v^0_h = \Pi_{h}\mal{P}_M v^0,~~c^0_h = \Pi_{h}\mal{P}_M c^0,
	% 	\end{equation*}
 %        such that $\theta_u^0=\theta_v^0=\theta_c^0=0$.

    Next, we present the main error estimate in Theorem \ref{thm:covergence} only for the parabolic-parabolic case (i.e., $\epsilon>0$); the corresponding error estimate for the parabolic-elliptic case (i.e., $\epsilon=0$) is analogous and can be obtained similarly. Without loss of generality, we take the physical parameters to be $\alpha = \beta = \epsilon \equiv 1$.
%%%%%%%%%%%%%%%%%%%%%%%%%%%%%%%%%%%%%%%%%%%%%%%%%%%%%%%
\begin{theorem}\label{thm:covergence}
	Assume that the regularity assumption \eqref{eq:regular} holds and that there exists two positive constants $Q_{\alpha}$ and $Q^{\alpha}$ such that 
    $Q_{\alpha} h \leq\tau \leq Q^{\alpha} h$ and $\tau \leq \tau^*$. Then there exists a positive constant $Q_*$, independent of $\tau$, $h$ and $n$, such that 
	\begin{equation*}
		\|u^n -  u_h^n\|_{\rm M} + \|v^n -  v_h^n\|_{\rm M} + \|c^n -  c_h^n\|_{\rm M} 
          + \|\nabla c^n-\bm{d} c_h^n\|_{\rm TM} \leq Q_*(\tau + h^2),\quad 0\leq n \leq N.
	\end{equation*}
\end{theorem}
%%%%%%%%%%%%%%%%%%%%%%%%%%%%%%%%%%%%%%%%%%%%%%%%%%%%%%%%%%%%

\subsection{High order consistency analysis}
By consistency, the projection solutions $\{u_M,v_M,c_M\}$ satisfy the proposed scheme \eqref{model:KS:shemeI}--\eqref{model:KS:shemeII} with truncation errors that are first-order accurate in time and second-order accurate in space. However, these errors are insufficient to establish a priori $\ell^{\infty}$ bounds for the numerical solutions, which are essential for guaranteeing the separation property. To overcome this difficulty, we employ a higher-order consistency analysis to construct perturbation solutions of $\{u_M,v_M,c_M\}$ that depend solely on the exact solutions $\{u,v,c\}$ and satisfy a consistency error of $\mal{O}(\tau^3+h^4)$.

In more detail, we shall construct supplementary solutions $\{z_{\tau,1}, z_{\tau,2}, z_{h,1}\}$ for $z=u,v,c$ step by step as follows. 

\paragraph{\indent \bf Part I. Second-order temporal accuracy}
First, taking $M$ large enough and applying Taylor expansion together with the projection approximation \eqref{projection:error}, we see that the projection solutions $\{u_M,v_M,c_M\}$ satisfy
\begin{equation}\label{projection:U_M}
	\begin{aligned}
		\frac{u_M^{n+1}-u_M^n}{\tau} &= \nabla \cdot (e^{c_M^n}\nabla(e^{-c_M^n} u_M^{n+1})) + \tau (G_u^{(0)})^{n+1} + \mal{O}(\tau^2) + \mal{O}(h^{m_0}),\\
		\frac{v_M^{n+1}-v_M^n}{\tau} &= \nabla \cdot (e^{c_M^n}\nabla(e^{-c_M^n} v_M^{n+1})) + \tau (G_v^{(0)})^{n+1} + \mal{O}(\tau^2) + \mal{O}(h^{m_0}),\\ 
		\frac{c_M^{n+1}-c_M^n}{\tau} &= \Delta c_M^{n+1} - c_M^{n+1} + u_M^{n+1} + v_M^{n+1} + \tau (G_c^{(0)})^{n+1} + \mal{O}(\tau^2) + \mal{O}(h^{m_0}),
	\end{aligned}
\end{equation}
where $m_0 \geq 4$, and the regularity assumption \eqref{eq:regular} imply that the leading-term truncation error functions $\{G_u^{(0)},G_v^{(0)},G_c^{(0)}\}$ are sufficiently smooth in the sense that their derivatives are bounded.

Then, we define the temporal correction functions $\{u_{\tau,1},v_{\tau,1},c_{\tau,1}\}$ with trivial initial conditions and periodic boundary conditions as follows:
\begin{equation}\label{partial_tau_1}
	\begin{aligned}
		\partial_t u_{\tau,1} &= \nabla \cdot \bigl(e^{c_M}\nabla(e^{-c_M}(u_{\tau,1}-c_{\tau,1}u_M)) + e^{c_M}c_{\tau,1}\nabla(e^{-c_M}u_M) \bigr) - G_u^{(0)},\\
		\partial_t v_{\tau,1} &= \nabla \cdot \bigl(e^{c_M}\nabla(e^{-c_M}(v_{\tau,1}-c_{\tau,1}v_M)) + e^{c_M}c_{\tau,1}\nabla(e^{-c_M}v_M) \bigr) - G_v^{(0)},\\
		\partial_t c_{\tau,1} &= \Delta c_{\tau,1} - c_{\tau,1} + u_{\tau,1} + v_{\tau,1}- G_c^{(0)}.
	\end{aligned}
\end{equation}
Note that the correction functions depend solely on the projection solutions $\{u_M,v_M,c_M\}$ and thus possess sufficient regularity. %, with their required derivatives uniformly bounded. 
Furthermore, by applying the semi-implicit temporal discretization to \eqref{partial_tau_1}, we obtain
\begin{equation}\label{dis:partial_tau_1}
	\begin{aligned}
		\frac{u_{\tau,1}^{n+1}-u_{\tau,1}^n}{\tau} &= \nabla \cdot \bigl(e^{c_M^n}\nabla(e^{-c_M^n}(u_{\tau,1}^{n+1}-c_{\tau,1}^n u_M^{n+1})) + e^{c_M^n}c_{\tau,1}^n\nabla(e^{-c_M^n}u_M^{n+1}) \bigr)\\
		&\quad\quad
        - (G_u^{(0)})^{n+1} + \mal{O}(\tau),\\
		\frac{v_{\tau,1}^{n+1}-v_{\tau,1}^n}{\tau} &= \nabla \cdot \bigl(e^{c_M^n}\nabla(e^{-c_M^n}(v_{\tau,1}^{n+1}-c_{\tau,1}^n v_M^{n+1})) + e^{c_M^n}c_{\tau,1}^n\nabla(e^{-c_M^n}v_M^{n+1}) \bigr)\\
		&\quad\quad
        - (G_v^{(0)})^{n+1}+ \mal{O}(\tau),\\
		\frac{c_{\tau,1}^{n+1}-c_{\tau,1}^n}{\tau} &= \Delta c_{\tau,1}^{n+1} - c_{\tau,1}^{n+1} + u_{\tau,1}^{n+1} + v_{\tau,1}^{n+1}- (G_c^{(0)})^{n+1}+ \mal{O}(\tau),
	\end{aligned}
\end{equation}
which is only first-order accurate in time.
Define $\check{u}_1:=u_M + \tau \mal{P}_M u_{\tau,1}$, $\check{v}_1:=v_M + \tau \mal{P}_M v_{\tau,1}$, $\check{c}_1:=c_M + \tau \mal{P}_M c_{\tau,1}$. Noting that
\begin{equation*}
	\begin{aligned}
		e^{\check{c}_1} %= e^{c_M + \tau \mal{P}_M c_{\tau,1}}  
        = e^{c_M} + \tau e^{c_M} \mal{P}_M c_{\tau,1} + \mal{O}(\tau^2),\quad
		e^{-\check{c}_1} %= e^{-(c_M + \tau \mal{P}_M c_{\tau,1})}  
        = e^{-c_M} - \tau e^{-c_M} \mal{P}_M c_{\tau,1} + \mal{O}(\tau^2).
	\end{aligned}
\end{equation*}
Then, it follows from \eqref{projection:U_M} and \eqref{dis:partial_tau_1} that $(\check{u}_1, \check{v}_1, \check{c}_1)$ satisfy the semi-discrete approximations 
\begin{equation}\label{dis:U_1}
	\begin{aligned}
		\frac{\check{u}_1^{n+1}-\check{u}_1^n}{\tau} &= \nabla \cdot (e^{\check{c}_1^n}\nabla(e^{-\check{c}_1^n} \check{u}_1^{n+1})) + \tau^2 (G_u^{(1)})^{n+1} + \mal{O}(\tau^3) + \mal{O}(h^{m_0}),\\
		\frac{\check{v}_1^{n+1}-\check{v}_1^n}{\tau} &= \nabla \cdot (e^{\check{c}_1^n}\nabla(e^{-\check{c}_1^n} \check{v}_1^{n+1})) + \tau^2 (G_v^{(1)})^{n+1} + \mal{O}(\tau^3) + \mal{O}(h^{m_0}),\\
		\frac{\check{c}_1^{n+1}-\check{c}_1^n}{\tau} &= \Delta \check{c}_1^{n+1} - \check{c}_1^{n+1} + \check{u}_1^{n+1} + \check{v}_1^{n+1} + \tau^2 (G_c^{(1)})^{n+1} + \mal{O}(\tau^3) + \mal{O}(h^{m_0}).
	\end{aligned}
\end{equation}
which is second-order accurate in time.

\paragraph{\indent \bf Part II. Third-order temporal accuracy}
Following a procedure analogous to that used in \eqref{partial_tau_1}, we further define the temporal correction functions $\{u_{\tau,2},v_{\tau,2},c_{\tau,2}\}$ with trivial initial conditions and periodic boundary conditions as follows:
\begin{equation}\label{partial_tau_2}
	\begin{aligned}
		\partial_t u_{\tau,2} &= \nabla \cdot (e^{\check{c}_1}\nabla(e^{-\check{c}_1}(u_{\tau,2}-c_{\tau,2}\check{u}_1)) + e^{\check{c}_1}c_{\tau,2}\nabla(e^{-\check{c}_1}\check{u}_1)) - G_u^{(1)},\\
		\partial_t v_{\tau,2} &= \nabla \cdot (e^{\check{c}_1}\nabla(e^{-\check{c}_1}(v_{\tau,2}-c_{\tau,2}\check{v}_1))  + e^{\check{c}_1}c_{\tau,2}\nabla(e^{-\check{c}_1}\check{v}_1) )- G_v^{(1)},\\
		\partial_t c_{\tau,2} &= \Delta c_{\tau,2} - c_{\tau,2} + u_{\tau,2} + v_{\tau,2}- G_c^{(1)}.
	\end{aligned}
\end{equation}
Again, the correction functions $\{u_{\tau,2},v_{\tau,2},c_{\tau,2}\}$ depend solely on the projection solutions $\{u_M,v_M,c_M\}$ and are sufficiently smooth. Similar, a semi-implicit temporal discretization to \eqref{partial_tau_2} leads to
\begin{equation}\label{dis:partial_tau_2}
	\begin{aligned}
		\frac{u_{\tau,2}^{n+1}-u_{\tau,2}^n}{\tau} &= \nabla \cdot \bigl(e^{\check{c}_1^n}\nabla(e^{-\check{c}_1^n}(u_{\tau,2}^{n+1}-c_{\tau,2}^n\check{u}_1^{n+1})) +e^{\check{c}_1^n}c_{\tau,2}^n\nabla(e^{-\check{c}_1^n}\check{u}_1^{n+1}) \bigr) \\
		&\quad\quad
        - (G_u^{(1)})^{n+1}+ \mal{O}(\tau),\\
		\frac{v_{\tau,2}^{n+1}-v_{\tau,2}^n}{\tau} &= \nabla \cdot \bigl(e^{\check{c}_1^n}\nabla(e^{-\check{c}_1^n}(v_{\tau,2}^{n+1}-c_{\tau,2}^n\check{v}_1^{n+1})) +e^{\check{c}_1^n}c_{\tau,2}^n\nabla(e^{-\check{c}_1^n}\check{v}_1^{n+1}) \bigr) \\
		&\quad\quad
        - (G_v^{(1)})^{n+1}+ \mal{O}(\tau),\\
		\frac{c_{\tau,2}^{n+1}-c_{\tau,2}^n}{\tau} &= \Delta c_{\tau,2}^{n+1} - c_{\tau,2}^{n+1} + u_{\tau,2}^{n+1} + v_{\tau,2}^{n+1}- (G_c^{(1)})^{n+1}+ \mal{O}(\tau),
	\end{aligned}
\end{equation}
which is again only first-order accurate in time. Define $\check{u}_2:=\check{u}_1 + \tau^2 \mal{P}_M u_{\tau,2}$, $\check{v}_2:=\check{v}_1 + \tau^2 \mal{P}_M v_{\tau,2}$, $\check{c}_2:=\check{c}_1 + \tau^2 \mal{P}_M c_{\tau,2}$. It then follows from  \eqref{dis:U_1} and \eqref{dis:partial_tau_2} that
\begin{equation}\label{dis:U_2}
	\begin{aligned}
		\frac{\check{u}_2^{n+1}-\check{u}_2^n}{\tau} &= \nabla \cdot (e^{\check{c}_2^n}\nabla(e^{-\check{c}_2^n} \check{u}_2^{n+1})) + \tau^3 (G_u^{(2)})^{n+1} + \mal{O}(\tau^4) + \mal{O}(h^{m_0}),\\
		\frac{\check{v}_2^{n+1}-\check{v}_2^n}{\tau} &= \nabla \cdot (e^{\check{c}_2^n}\nabla(e^{-\check{c}_2^n} \check{v}_2^{n+1})) + \tau^3 (G_v^{(2)})^{n+1} + \mal{O}(\tau^4) + \mal{O}(h^{m_0}),\\
		\frac{\check{c}_2^{n+1}-\check{c}_2^n}{\tau} &= \Delta \check{c}_2^{n+1} - \check{c}_2^{n+1} + \check{u}_2^{n+1} + \check{v}_2^{n+1} + \tau^3 (G_c^{(2)})^{n+1} + \mal{O}(\tau^4) + \mal{O}(h^{m_0}),
	\end{aligned}
\end{equation}
which is third-order accurate in time.

\paragraph{\indent \bf Part III. Fourth-order spatial accuracy}
 To further improve the spatial consistency, we introduce spatial correction functions as follows. 
First, a direct spatial BCFD approximation of \eqref{dis:U_2} leads to a fully discrete scheme analogues to \eqref{model:KS:shemeI}--\eqref{model:KS:shemeII}, namely
\begin{equation}\label{dis:U_2_h}
	\begin{aligned}
		\frac{\check{u}_2^{n+1}-\check{u}_2^n}{\tau} &= D_x \bigl([\mal{A}_x e^{\check{c}_2^n}][d_x(e^{-\check{c}_2^n} \check{u}_2^{n+1})]\bigr)+D_y \bigl([\mal{A}_y e^{\check{c}_2^n}][d_y(e^{-\check{c}_2^n} \check{u}_2^{n+1})]\bigr),\\
		&\quad\quad + \tau^3 (G_u^{(2)})^{n+1} + h^2(H_u^{(0)})^{n+1}+ \mal{O}(\tau^4+h^4),\\
		\frac{\check{v}_2^{n+1}-\check{v}_2^n}{\tau} &= D_x \bigl([\mal{A}_x e^{\check{c}_2^n}] [d_x(e^{-\check{c}_2^n} \check{v}_2^{n+1})]\bigr)+D_y \bigl([\mal{A}_y e^{\check{c}_2^n}] [d_y(e^{-\check{c}_2^n} \check{v}_2^{n+1})]\bigr),\\
		&\quad\quad + \tau^3 (G_v^{(2)})^{n+1} + h^2(H_v^{(0)})^{n+1}+ \mal{O}(\tau^4+h^4),\\
		\frac{\check{c}_2^{n+1}-\check{c}_2^n}{\tau} &= D_xd_x \check{c}_2^{n+1} + D_yd_y \check{c}_2^{n+1} - \check{c}_2^{n+1} + \check{u}_2^{n+1} + \check{v}_2^{n+1} \\
		&\quad\quad+ \tau^3 (G_c^{(2)})^{n+1}+ h^2(H_c^{(0)})^{n+1} + \mal{O}(\tau^4+h^4),
	\end{aligned}
\end{equation}
which is second-order accurate in space and third-order accurate in time, and under regularity assumption \eqref{eq:regular}, the leading-term truncation error functions $\{H_u^{(0)},H_v^{(0)},H_c^{(0)}\}$ possess sufficient regularity. 

Next, we define the spatial correction functions $\{u_{h,1},v_{h,1},c_{h,1}\}$  with trivial initial conditions and periodic boundary conditions as follows:
\begin{equation}\label{partial_h_1}
	\begin{aligned}
		\partial_t u_{h,1} &= \nabla \cdot (e^{\check{c}_2}\nabla(e^{-\check{c}_2}(u_{h,1}-c_{h,1}\check{u}_2)) + e^{\check{c}_2}c_{h,1}\nabla(e^{-\check{c}_2}\check{u}_2)) - H_u^{(0)},\\
		\partial_t v_{h,1} &= \nabla \cdot (e^{\check{c}_2}\nabla(e^{-\check{c}_2}(v_{h,1}-c_{h,1}\check{v}_2)) + e^{\check{c}_2}c_{h,1}\nabla(e^{-\check{c}_2}\check{v}_2)) - H_v^{(0)},\\
		\partial_t c_{h,1} &= \Delta c_{h,1} - c_{h,1} + u_{h,1} + v_{h,1}- H_c^{(0)}.
	\end{aligned}
\end{equation}
Note that the correction functions $\{u_{h,1},v_{h,1},c_{h,1}\}$ depend only on the projection solutions $\{u_M,v_M,c_M\}$. Therefore, all their required divided differences are uniformly bounded. Then, the fully discrete approximation of \eqref{partial_h_1} yields
\begin{equation}\label{dis:partial_h_1}
	\begin{aligned}
		\frac{u_{h,1}^{n+1}-u_{h,1}^n}{\tau} &= D_x \bigl ([\mal{A}_x e^{\check{c}_2^n}][d_x(e^{-\check{c}_2^n}(u_{h,1}^{n+1}-c_{h,1}^n\check{u}_2^{n+1}))]+[\mal{A}_x (e^{\check{c}_2^n} c_{h,1}^n)] [d_x(e^{-\check{c}_2^n}\check{u}_2^n)]\bigr)\\
        &+ D_y \bigl ([\mal{A}_y e^{\check{c}_2^n}][d_y(e^{-\check{c}_2^n}(u_{h,1}^{n+1}c_{h,1}^n\check{u}_2^{n+1}))]+[\mal{A}_y (e^{\check{c}_2^n} c_{h,1}^n)] [d_y(e^{-\check{c}_2^n}\check{u}_2^n)]\bigr)\\
	    &- (H_u^{(0)})^{n+1}+ \mal{O}(\tau + h^2) ,\\
		\frac{v_{h,1}^{n+1}-v_{h,1}^n}{\tau} &= D_x\bigl([\mal{A}_x e^{\check{c}_2^n}] [d_x(e^{-\check{c}_2^n}(v_{h,1}^{n+1}-c_{h,1}^n\check{v}_2^{n+1}))]+[\mal{A}_x (e^{\check{c}_2^n} c_{h,1}^n)] [d_x(e^{-\check{c}_2^n}\check{v}_2^n)]\bigr)\\
        &+D_y\bigl([\mal{A}_y e^{\check{c}_2^n}] [d_y(e^{-\check{c}_2^n}(v_{h,1}^{n+1}-c_{h,1}^n\check{v}_2^{n+1}))]+[\mal{A}_y (e^{\check{c}_2^n} c_{h,1}^n)] [d_y(e^{-\check{c}_2^n}\check{v}_2^n)]\bigr)\\
        &- (H_v^{(0)})^{n+1}+ \mal{O}(\tau + h^2),\\
		\frac{c_{h,1}^{n+1}-c_{h,1}^n}{\tau} &= D_xd_x c_{h,1}^{n+1} + D_yd_y c_{h,1}^{n+1}- c_{h,1}^{n+1} + u_{h,1}^{n+1} + v_{h,1}^{n+1}- (H_c^{(0)})^{n+1}+ \mal{O}(\tau + h^2),
	\end{aligned}
\end{equation}
which is first-order accurate in time and second-order accurate in space. 

Now, we define the final corrected solutions $ \check{z}$ for $z=u,v,c$ as follows:
\begin{equation}\label{check_U}
		\check{z}:= \check{z}_2 + h^2 \mal{P}_M z_{h,1}= z_M + \mal{P}_M(\tau z_{\tau,1} + \tau^2 z_{\tau,2} + h^2 z_{h,1}),
\end{equation}
where $\{z_{\tau,1}, z_{\tau,2}, z_{h,1}\}$ are specified as above. Thus, multiplying both sides of \eqref{dis:partial_h_1} by $h^2$ and adding the resulting equation to \eqref{dis:U_2_h}, together with the time-step condition $ Q_{\alpha} h \leq \tau$, it is easy to see that the constructed solutions $\{\check{u}, \check{v}, \check{c}\}$ defined in \eqref{check_U} satisfy the proposed scheme \eqref{model:KS:shemeI}--\eqref{model:KS:shemeII} with a higher-order consistency error of $\mal{O}(\tau^3+h^4)$. This result is summarized as follows. 
\begin{lemma}\label{lem:high-order}
   The constructed correction solutions $\{\check{u}, \check{v}, \check{c}\}$ satisfy the following fully discrete equations with truncation errors $R_u^{n+1}$, $R_v^{n+1}$, and $R_c^{n+1}$:
\begin{align}
	\frac{\check{u}^{n+1}-\check{u}^n}{\tau} &= D_x \bigl([\mal{A}_x e^{\check{c}^n}] [d_x(e^{-\check{c}^n} \check{u}^{n+1})]\bigr)+D_y \bigl([\mal{A}_y e^{\check{c}^n}] [d_y(e^{-\check{c}^n} \check{u}^{n+1})]\bigr)+R_u^{n+1},\label{trun:u}\\
	\frac{\check{v}^{n+1}-\check{v}^n}{\tau} &= D_x \bigl([\mal{A}_x e^{\check{c}^n}] [d_x(e^{-\check{c}^n} \check{v}^{n+1})]\bigr)+D_y \bigl([\mal{A}_y e^{\check{c}^n}] [d_y(e^{-\check{c}^n} \check{v}^{n+1})]\bigr)+R_v^{n+1},\label{trun:v}\\
	\frac{\check{c}^{n+1}-\check{c}^n}{\tau} &= D_xd_x \check{c}^{n+1}+D_yd_y \check{c}^{n+1} - \check{c}^{n+1} + \check{u}^{n+1} + \check{v}^{n+1} + R_c^{n+1}.\label{trun:c}
\end{align}
Moreover, provided that $Q_{\alpha} h \leq \tau$, we have the following estimates:
\begin{equation*}
    \|R_u^{n+1}\|_{\rm M},\,
    \|R_v^{n+1}\|_{\rm M},\,\|R_c^{n+1}\|_{\rm M} \leq \check{Q}(\tau^3+h^4),
\end{equation*}
for some positive constant $\check{Q}$.
\end{lemma}

Note that $\check{c}\in \mal{B}^K$. The projection error estimate \eqref{projection:error}, together with \eqref{bound:exc1a}, ensures that for sufficiently small $\tau$ and $h$, the corrected solution $\check{c}$ satisfies a separation property analogous to \eqref{bound:exc1}, namely,
	\begin{equation}\label{bound:exc2}
		\check{c}(\bm x,t) \geq \varepsilon >0.
	\end{equation}

\begin{remark}	As discussion, the correction functions $\{\check{u},\check{v},\check{c}\}$ depend solely on the exact solutions $\{u,v,c\}$, and therefore, their discrete $W^{1,\infty}$ norms are uniformly bounded. More precisely, there exists a constant $Q^*>0$, independent of $\tau$ and $h$, such that
	\begin{equation}\label{eq:regular:2}
		\begin{aligned}
			\|\check{z}^n\|_{\infty} +\|\bm{d} \check{z}^n\|_{\infty} &\leq Q^*,~~ z=\{u,v,c\},~~\forall n\geq0.
		\end{aligned}
	\end{equation}
 Finally, we point out that the higher-order consistency expansion shall establish an a priori $\ell^{\infty}$ bound for the chemoattractant concentration $c_h$, which is essential for verifying the separation property required in the subsequent error estimates.
\end{remark}

\subsection{Error estimates}
In this subsection, we establish the convergence of the proposed numerical scheme \eqref{model:KS:shemeI}--\eqref{model:KS:shemeII}. To this end, we introduce the following grid error functions:
	\begin{equation}\label{error:2}
		\widetilde{\theta}_u^n := \check{u}^n -  u_h^n,~~\widetilde{\theta}_v^n := \check{v}^n - v_h^n,~~,\widetilde{\theta}_c^n := \check{c}^n - c_h^n,~~0 \leq n \leq N.
	\end{equation}
Then, subtracting \eqref{model:KS:shemeI}--\eqref{model:KS:shemeII} from \eqref{trun:u}--\eqref{trun:c}, we obtain the following error equations:
	\begin{subequations}%\label{error:function}
	\begin{align}
		\frac{\widetilde{\theta}_u^{n+1}-\widetilde{\theta}_u^n}{\tau} &= D_x \bigl([\mal{A}_x e^{\check{c}^n}] [d_x(e^{-\check{c}^n} \check{u}^{n+1})] - [\mal{A}_x e^{c_h^n}] [d_x (e^{-c_h^{n}}u_h^{n+1})]\bigr)\nonumber\\
        &+D_y \bigl([\mal{A}_y e^{\check{c}^n}] [d_y(e^{-\check{c}^n} \check{u}^{n+1})] - [\mal{A}_y e^{c_h^n}] [d_y (e^{-c_h^{n}}u_h^{n+1})]\bigr)+R_u^{n+1}\label{error:function_u},\\
		\frac{\widetilde{\theta}_v^{n+1}-\widetilde{\theta}_v^n}{\tau} &= D_x \bigl([\mal{A}_x e^{\check{c}^n}] [d_x(e^{-\check{c}^n} \check{v}^{n+1})] - [\mal{A}_x e^{c_h^n}] [d_x (e^{-c_h^{n}}v_h^{n+1})]\bigr)\nonumber\\
        &+D_y \bigl([\mal{A}_y e^{\check{c}^n}] [d_y(e^{-\check{c}^n} \check{v}^{n+1})] - [\mal{A}_y e^{c_h^n}] [d_y (e^{-c_h^{n}}v_h^{n+1})]\bigr)+R_v^{n+1}\label{error:function_v},\\
		\frac{\widetilde{\theta}_c^{n+1}-\widetilde{\theta}_c^n}{\tau} &= D_x(d_x \widetilde{\theta}_c^{n+1})+D_y(d_y \widetilde{\theta}_c^{n+1}) - \widetilde{\theta}_c^{n+1} + \widetilde{\theta}_u^{n+1} + \widetilde{\theta}_v^{n+1} + R_c^{n+1}\label{error:function_c},
	\end{align}
\end{subequations}
where, by \eqref{check_U}, the initial errors satisfy 
        % \begin{equation*}
			%\check{u}^0 = \mal{P}_M u^0,~~\check{v}^0_h = \mal{P}_M v^0,~~\check{c}^0_h = \mal{P}_M c^0,
            $\widetilde{\theta}_u^0=\widetilde{\theta}_v^0=\widetilde{\theta}_c^0=0.$
		% \end{equation*}
% Moreover, the introduced error functions benefit from the higher-order consistency estimates \eqref{trun:u}--\eqref{trun:c}, thereby allowing for more accurate error estimates.

Throughout the nonlinear error analysis, we make the following a priori assumption at the previous time level:
\begin{equation}\label{priori_assum}
	 \|\widetilde{\theta}_c^n\|_{\rm M}\leq \tau^{\frac{9}{4}} + h^{\frac{9}{4}}.
\end{equation}
%The validity of \eqref{priori_assum} at the next time level will be established in the subsequent error analysis. 
Under the above assumption and the time-step condition $\tau \leq Q^{\alpha}h$, we can directly obtain the following discrete $W^{1,\infty}$ bounds for the error function:
\begin{align*}
	&\|\widetilde{\theta}_c^n\|_{\infty} 
    \leq Qh^{-1}\|\widetilde{\theta}_c^n\|_{\rm M} 
    \leq Qh^{-1}(\tau^{\frac{9}{4}} + h^{\frac{9}{4}}) 
    \leq Q(\tau^{\frac{5}{4}} + h^{\frac{5}{4}}) 
    \leq 1,\\
	&\|\bm{d}\widetilde{\theta}_c^n\|_{\infty} 
    \leq Qh^{-1}\|\widetilde{\theta}_c^n\|_{\infty} 
    \leq Qh^{-1}(\tau^{\frac{5}{4}} + h^{\frac{5}{4}}) 
    \leq Q(\tau^{\frac{1}{4}} + h^{\frac{1}{4}}) \leq 1. 
\end{align*}
Then, combined with the uniformly boundedness result \eqref{eq:regular:2} and the triangle inequality, the  discrete $W^{1,\infty}$ bound for the numerical solution is obtained as follows:
\begin{equation}\label{infty:cn}
    \begin{aligned}
	\|c_h^n\|_{\infty} \leq \|\check{c}^n\|_{\infty} + \|\widetilde{\theta}_c^n\|_{\infty} \leq \widetilde{Q},\quad
	\|\bm{d} c_h^n\|_{\infty} \leq \|\bm{d}\check{c}^n\|_{\infty} + \|\bm{d}\widetilde{\theta}_c^n\|_{\infty} \leq \widetilde{Q},
\end{aligned}
\end{equation}
where $\widetilde{Q}:=Q^*+1$ is independent of $\tau$ and $h$. 
%Furthermore, for sufficiently small $\tau$ and $h$, the error estimate  implies
% \begin{equation}\label{tilde:error:cn}
% 		\|\widetilde{\theta}_c^n\|_{\infty} \leq Q (\tau^{\frac{5}{4}} + h^{\frac{5}{4}}) \leq \frac{\varepsilon}{2}.
% \end{equation}
Consequently, for sufficiently small $\tau$ and $h$, the separation property \eqref{bound:exc2} yields the following uniform positive lower bound for the numerical solution:
\begin{equation}\label{eq:bound:cn}
	c_h^n \geq \check{c}^n -\|\widetilde{\theta}_c^n\|_{\infty} \ge \check{c}^n - Q (\tau^{\frac{5}{4}} + h^{\frac{5}{4}})\geq \varepsilon-\frac{\varepsilon}{2} \geq \frac{\varepsilon}{2}.
\end{equation}

In the following lemma, we establish an estimate for $\|\bm{d}(e^{-\check{c}^n}-e^{-c_h^n})\|_{\rm TM}$, which plays a crucial role in the main error analysis.
\begin{lemma}\label{lem:grad_c}
	Under the a priori assumption \eqref{priori_assum}, there exist positive constants $Q_1$, $Q_2$, and $Q_3$, independent of $\tau$ and $h$, such that
	\begin{equation*}
		\|\bm{d}(e^{-\check{c}^n}-e^{-c_h^n})\|_{\rm TM} \leq Q_1\|\widetilde{\theta}_c^n\|_{\rm M} + Q_2 h^4 + Q_3 \|\bm{d}\widetilde{\theta}_c^n\|_{\rm TM}.
	\end{equation*}
\end{lemma}
\begin{proof}	First, applying the mean value theorem yields
	\begin{equation}\label{eq:Dxe}
		\begin{aligned}
			 e^{-\check{c}^n_{i+1,j}} -e^{-\check{c}^n_{i,j}} = -h e^{-\xi} [d_x\check{c}^n]_{i+1/2,j}, \quad
          e^{-c^n_{i+1,j}}- e^{-c^n_{i,j}} 	=   he^{-\eta} [d_x c^n]_{i+1/2,j},
		\end{aligned}
	\end{equation}
where $\xi$ lies between $\check{c}_{i,j}^n$ and $\check{c}_{i+1,j}^n$, and $\eta$ lies between $c^n_{i,j}$ and $c^n_{i+1,j}$. Moreover, by Taylor expansion, we have
	\begin{equation*}\label{eq:xi}
		\begin{aligned}
			e^{-\xi} = -\frac{e^{-\check{c}_{i+1,j}^n}-e^{-\check{c}_{i,j}^n}}{\check{c}_{i+1,j}^n-\check{c}_{i,j}^n} 
            &= e^{-\check{c}_{i,j}^n} - \frac{e^{-\check{c}_{i,j}^n}}{2}(\check{c}_{i+1,j}^n-\check{c}_{i,j}^n ) 
            + \frac{e^{-\check{c}_{i,j}^n}}{6}(\check{c}_{i+1,j}^n-\check{c}_{i,j}^n )^2\\
			&\quad -\frac{e^{-\check{c}_{i,j}^n}}{24}(\check{c}_{i+1,j}^n-\check{c}_{i,j}^n )^3 
            +\frac{e^{-\xi_M}}{120} (\check{c}_{i+1,j}^n-\check{c}_{i,j}^n )^4,\\
			e^{-\eta} = -\frac{e^{-c^n_{i+1,j}}-e^{-c^n_{i,j}}}{c^n_{i+1,j}-c^n_{i,j}} 
            &= e^{-c^n_{i,j}} + \frac{e^{-c^n_{i,j}}}{2} (c^n_{i+1,j}-c^n_{i,j}) 
            + \frac{e^{-c^n_{i,j}}}{6} (c^n_{i+1,j}-c^n_{i,j})^2\\
			&\quad -\frac{e^{-c^n_{i,j}}}{24} (c^n_{i+1,j}-c^n_{i,j})^3 
            +\frac{e^{-\eta_M}}{120} (c^n_{i+1,j}-c^n_{i,j})^4,
		\end{aligned}
	\end{equation*}
for some $\xi_M$ between $\check{c}_{i,j}^n$ and $\check{c}_{i+1,j}^n$, and $\eta_M$ between $c^n_{i,j}$ and $c^n_{i+1,j}$.   Thus, the following relation holds:
	\begin{equation}\label{eq:eta-xi}
		\begin{aligned}
			e^{-\eta}-e^{-\xi} &= -( e^{-\check{c}_{i,j}^n}-e^{-c^n_{i,j}}) \\
            &\quad + \frac{1}{2}\bigl[ e^{-\check{c}_{i,j}^n}(\check{c}_{i+1,j}^n-\check{c}_{i,j}^n) - e^{-c^n_{i,j}} (c^n_{i+1,j}-c^n_{i,j}) \bigr]\\
			&\quad - \frac{1}{6}\bigl[ e^{-\check{c}_{i,j}^n}(\check{c}_{i+1,j}^n-\check{c}_{i,j}^n)^2 - e^{-c^n_{i,j}} (c^n_{i+1,j}-c^n_{i,j})^2 \bigr]\\
			&\quad + \frac{1}{24}\bigl[ e^{-\check{c}_{i,j}^n}(\check{c}_{i+1,j}^n-\check{c}_{i,j}^n)^3 - e^{-c^n_{i,j}} (c^n_{i+1,j}-c^n_{i,j})^3 \bigr]\\
			&\quad - \frac{1}{120}\bigl[ e^{-\xi_M}(\check{c}_{i+1,j}^n-\check{c}_{i,j}^n)^4 - e^{-\eta_M} (c^n_{i+1,j}-c^n_{i,j})^4 \bigr] =: \sum_{\ell=1}^5 \mal{L}_\ell.
		\end{aligned}
	\end{equation}
    
Next, we estimate the right-hand side of \eqref{eq:eta-xi} term by term.  For the first term $\mal{L}_1$, an application of the separation estimate \eqref{eq:bound:cn} yields
	\begin{equation*}
		|\mal{L}_1| %= |e^{-\check{c}_{i,j}^n}-e^{-c^n_{i,j}}| 
        = |e^{-\xi_N}[\widetilde{\theta}_c^n]_{i,j}| \leq e^{-\frac{\varepsilon}{2}}|[\widetilde{\theta}_c^n]_{i,j}|=:Q^{(1)}|[\widetilde{\theta}_c^n]_{i,j}|,
	\end{equation*}
where $\xi_N$ lies $\check{c}^n_{i,j}$ and $c^n_{i,j}$.
    
	For the second term $\mal{L}_2$, combining the triangle inequality with the uniform boundedness result \eqref{eq:regular:2} and the lower error bound \eqref{eq:bound:cn}, we obtain
	\begin{equation*}
		\begin{aligned}
			|\mal{L}_2|% &= \frac{1}{2}\left|\bigl( e^{-\check{c}_{i,j}^n}(\check{c}_{i+1,j}^n-\check{c}_{i,j}^n) - e^{-c^n_{i,j}} (c^n_{i+1,j}-c^n_{i,j}) \bigr)\right|\\
			&=\frac{1}{2}\left|e^{-\check{c}_{i,j}^n}([\widetilde{\theta}_c^n]_{i+1,j}-[\widetilde{\theta}_c^n]_{i,j}) - (e^{-\check{c}_{i,j}^n}-e^{-c^n_{i,j}}) (c^n_{i+1,j}-c^n_{i,j}) \right|\\
            &\leq\frac{1}{2}\bigl|e^{-\check{c}_{i,j}^n}([\widetilde{\theta}_c^n]_{i+1,j}|-[\widetilde{\theta}_c^n]_{i,j})\bigr| + \f{1}{2}\bigl|e^{-\frac{\varepsilon}{2}} [\widetilde{\theta}_c^n]_{i.j} (c^n_{i+1,j}-c^n_{i,j}) \bigr|\\
			&\leq \frac{1} {2}e^{-\epsilon}\bigl(|[\widetilde{\theta}_c^n]_{i+1,j}|+|[\widetilde{\theta}_c^n]_{i,j}| \bigr) + \widetilde{Q}e^{-\frac{\varepsilon}{2}}|[\widetilde{\theta}_c^n]_{i,j}|
			=: Q^{(2)}\bigl(|[\widetilde{\theta}_c^n]_{i+1,j}|+|[\widetilde{\theta}_c^n]_{i,j}| \bigr),
		\end{aligned}
	\end{equation*}
    where \eqref{infty:cn} and the estimate result for $\mal{L}_1$ have also been applied.
	 Analogously, we obtain the bounds for the third and fourth terms:
	\begin{equation*} 
		|\mal{L}_3| + |\mal{L}_4| \leq \bigl(Q^{(3)}+Q^{(4)}\bigr)\bigl(|[\widetilde{\theta}_c^n]_{i+1,j}|+|[\widetilde{\theta}_c^n]_{i,j}| \bigr).
	\end{equation*}
	
    Finally, for the remainder term $\mal{L}_5$, the discrete $W^{1,\infty}$ bound \eqref{eq:regular:2} and \eqref{infty:cn} imply that
    \begin{equation}\label{eq:L5:h}
        \begin{aligned}
		\bigl|\check{c}_{i+1,j}^n-\check{c}_{i,j}^n\bigr| &= h\bigl|[d_x\check{c}^n]_{i+1/2,j}\bigr|\leq h\|[d_x\check{c}^n]_{i+1/2,j}\|_{\infty}\leq Q^* h,\\
		\bigl|c^n_{i+1,j}-c^n_{i,j}\bigr| &= h\bigl|[d_xc^n]_{i+1/2,j}\bigr|\leq h\|[d_x c^n]_{i+1/2,j}\|_{\infty} \leq \widetilde{Q}h.
	    \end{aligned}
    \end{equation}
	Thus, by inserting \eqref{eq:L5:h} into $\mal{L}_5$ and using the separation estimates \eqref{bound:exc2} and \eqref{eq:bound:cn}, we have
	\begin{equation*} 
		|\mal{L}_5| \leq  \f{e^{-\f{\varepsilon}{2}}}{120}\bigl( (Q^*)^4+(\widetilde{Q})^4\bigr) h^4 =:Q^{(5)} h^4.
	\end{equation*}
	Now, inserting the estimates for $\mal{L}_1$--$\mal{L}_5$ into \eqref{eq:eta-xi}, we obtain
	\begin{equation}\label{estimate:eta-xi}
		|e^{-\eta}-e^{-\xi}| \leq Q^{(0)} \bigl(|[\widetilde{\theta}_c^n]_{i+1,j}|+|[\widetilde{\theta}_c^n]_{i,j}| \bigr) +  Q^{(5)} h^4,\quad Q^{(0)} := \sum_{k=1}^4Q^{(k)}.
	\end{equation}
%	with $Q^{(0)} := \sum_{k=1}^4Q^{(k)}$. 

Recalling  \eqref{eq:Dxe}, we have
	\begin{equation}\label{eq:Dxe:e1}
		\begin{aligned}
			[d_x(e^{-\check{c}^n}-e^{-c_h^n})]_{i+1/2,j} 
            % &= \frac{1}{h}\bigl( e^{-\check{c}^n_{i+1,j}} -e^{-\check{c}^n_{i,j}} \bigr)  - \frac{1}{h}\bigl( e^{-c^n_{i+1,j}}- e^{-c^n_{i,j}} \bigr)\\
			 &= -e^{-\xi} [d_x\check{c}^n]_{i+1/2,j} + e^{-\eta} [d_x c^n]_{i+1/2,j}\\
			&=  (e^{-\eta}-e^{-\xi})[d_x\check{c}^n]_{i+1/2,j} - e^{-\eta}[d_x \widetilde{\theta}_c^n]_{i+1/2,j}.
		\end{aligned}
	\end{equation}
Thus,	combining the discrete $W^{1,\infty}$ bound \eqref{eq:regular:2} with the separation estimate \eqref{eq:bound:cn}, we conclude from \eqref{estimate:eta-xi} and \eqref{eq:Dxe:e1} that
	\begin{equation}\label{eq:Dxe:point}
		\begin{aligned}
			[d_x(e^{-\check{c}^n}-e^{-c_h^n})]_{i+1/2,j} 
            % &\leq  \bigl|e^{-\eta}-e^{-\xi}\bigr| \bigl|[d_x\check{c}^n]_{i+1/2,j}\bigr| + \bigl|e^{-\eta}\bigr| \bigl|[d_x \widetilde{\theta}_c^n]_{i+1/2,j}\bigr|\\
			\leq Q^*Q^{(0)} \bigl(|[\widetilde{\theta}_c^n]_{i+1,j}|+|[\widetilde{\theta}_c^n]_{i,j}| \bigr) +  Q^*Q^{(5)} h^4 + Q^{(1)} \bigl|[d_x \widetilde{\theta}_c^n]_{i+1/2,j}\bigr|.
		\end{aligned}
	\end{equation}
    Similarly, an analogous argument in the $y$-direction yields
    	\begin{equation}\label{eq:Dye:point}
		\begin{aligned}
			[d_y(e^{-\check{c}^n}-e^{-c_h^n})]_{i,j+1/2} 
            %&\leq  \bigl|e^{-\widetilde{\eta}}-e^{-\widetilde{\xi}}\bigr| \bigl|[d_y\check{c}^n]_{i,j+1/2}\bigr| + \bigl|e^{-\widetilde{\eta}}\bigr| \bigl|[d_y \widetilde{\theta}_c^n]_{i,j+1/2}\bigr|\\
			\leq Q^*Q^{(0)} \bigl(|[\widetilde{\theta}_c^n]_{i,j+1}|+|[\widetilde{\theta}_c^n]_{i,j}| \bigr) +  Q^*Q^{(5)} h^4 + Q^{(1)} \bigl|[d_y \widetilde{\theta}_c^n]_{i,j+1/2}\bigr|.
		\end{aligned}
	\end{equation}
%   where $\widetilde{\xi}$ lies between $\check{c}_{i,j}^n$ and $\check{c}_{i,j+1}^n$, and $\widetilde{\eta}$ lies between $c^n_{i,j}$ and $c^n_{i,j+1}$. 

Therefore, summing \eqref{eq:Dxe:point}--\eqref{eq:Dye:point} over all corresponding grid points, we obtain
	\begin{equation*}
		\begin{aligned}
			\|\bm{d}(e^{-\check{c}^n}-e^{-c_h^n})\|_{\rm TM} \leq Q_1\|\widetilde{\theta}_c^n\|_{\rm M} + Q_2 h^4 + Q_3 \|\bm{d}\widetilde{\theta}_c^n\|_{\rm TM}, 
		\end{aligned}
	\end{equation*}
	where $Q_1:= 2\sqrt{6}Q^*Q^{(0)}$, $Q_2:=\sqrt{6|\Omega|}Q^*Q^{(5)}$, and $Q_3=\sqrt{3}Q^{(1)}$. This completes the proof of the lemma.
\end{proof}

%%%%%%%%%%%%%%%%%%%%%%%%%%%%%%%%%%%%%%%%%%%%%%%%%%%%%%%%%%%%%%%%%%%%%%%%%%%%
Finally, with these preparations in place, we turn to the proof of Theorem \ref{thm:covergence}, which is divided into three parts.
\paragraph{\indent \bf Part I. Error estimates for  $\|\widetilde{\theta}_u^{n+1}\|_{\rm M}$, $\|\widetilde{\theta}_v^{n+1}\|_{\rm M}$, and $\|\widetilde{\theta}_c^{n+1}\|_{\rm M}$}
We first establish a partial error estimate for $\|\widetilde{\theta}_u^{n+1}\|_{\rm M}$.
Taking the discrete inner product of \eqref{error:function_u} with $\widetilde{\theta}_u^{n+1}$ and using Lemma \ref{lemma:Dd}, we obtain
\begin{equation}\label{error:u}
	\begin{aligned}
		&\frac{1}{2\tau} \bigl( \|\widetilde{\theta}_u^{n+1}\|_{\rm M}^2 - \|\widetilde{\theta}_u^n\|_{\rm M}^2 +\|\widetilde{\theta}_u^{n+1}-\widetilde{\theta}_u^n\|_{\rm M}^2  \bigr) \\
		&\qquad + \bigl( [\mal{A} e^{\check{c}^n}]\odot[\bm{d}(e^{-\check{c}^n} \check{u}^{n+1})] - [\mal{A} e^{c_h^n}]\odot [\bm{d} (e^{-c_h^{n}}u_h^{n+1})],\bm{d} \widetilde{\theta}_u^{n+1} \bigr)_{\rm TM} \\
        & \quad = \bigl ( R_u^{n+1},\widetilde{\theta}_u^{n+1} \bigr )_{\rm M} 
        \le \frac{1}{2} \bigl(\|\widetilde{\theta}_u^{n+1}\|_{\rm M}^2 + \|R_u^{n+1}\|_{\rm M}^2\bigr),
	\end{aligned}
\end{equation}
where $\odot$ denotes the Hadamard product of two vectors. In the following, we analyze in detail the second term on the left-hand side of \eqref{error:u}, which can be decomposed as follows:
\begin{equation*}
	\begin{aligned}
		&\bigl ( [\mal{A} e^{\check{c}^n}]\odot[\bm{d}(e^{-\check{c}^n} \check{u}^{n+1})] - [\mal{A} e^{c_h^n}]\odot [\bm{d} (e^{-c_h^{n}}u_h^{n+1})],\bm{d} \widetilde{\theta}_u^{n+1}\bigr )_{\rm TM}\\
		&= \bigl ( [\mal{A} e^{\check{c}^n}] \odot [\bm{d}((e^{-\check{c}^n} - e^{-c_h^n})\check{u}^{n+1})],\bm{d} \widetilde{\theta}_u^{n+1} \bigr )_{\rm TM} 
        + \bigl ( [\mal{A} e^{c_h^n}]\odot[\bm{d}(e^{-c_h^n}\widetilde{\theta}_u^{n+1})],\bm{d} \widetilde{\theta}_u^{n+1} \bigr)_{\rm TM}\\
		&\quad + \bigl ( [\mal{A} (e^{\check{c}^n}- e^{c_h^n})]\odot [\bm{d}(e^{-c_h^n}\check{u}^{n+1})],\bm{d} \widetilde{\theta}_u^{n+1} \bigr )_{\rm TM} =: \sum_{\ell=1}^3 \mal{J}_\ell.
	\end{aligned}
\end{equation*}	

For the first term $\mal{J}_1$, applying the separation estimate \eqref{eq:bound:cn} and the estimate established in Lemma \ref{lem:grad_c}, we have
\begin{equation*}
	\begin{aligned}
		\mal{J}_1 %&= \bigl ( [\mal{A} e^{\check{c}^n}][\bm{d}((e^{-\check{c}^n} - e^{-c_h^n})\check{u}^{n+1}],\bm{d} \widetilde{\theta}_u^{n+1} \bigr )_{\rm TM}\\
		& \leq \|\mal{A} e^{\check{c}^n}\|_{\infty} \|\bm{d}((e^{-\check{c}^n} - e^{-c_h^n})\check{u}^{n+1})\|_{\rm TM} \|\bm{d} \widetilde{\theta}_u^{n+1}\|_{\rm TM}\\
		& \leq \|\mal{A} e^{\check{c}^n}\|_{\infty}\, 
        \bigl( \|\check{u}^{n+1}\|_{\infty}\|\bm{d}(e^{-\check{c}^n} - e^{-c_h^n})\|_{\rm TM} +\|\bm{d}\check{u}^{n+1}\|_{\infty}\|e^{-\check{c}^n} - e^{-c_h^n}\|_{\rm M} \bigr) \, \|\bm{d} \widetilde{\theta}_u^{n+1}\|_{\rm TM}\\
		&\leq \frac{Q^{(6)}}{4\delta} \|\widetilde{\theta}_c^n\|_{\rm M}^2 + \frac{Q^{(7)}}{4\delta} \|\bm{d}\widetilde{\theta}_c^n\|_{\rm TM}^2 + \frac{Q^{(8)}}{4\delta} h^8 + 3\delta\|\bm{d} \widetilde{\theta}_u^{n+1}\|_{\rm TM}^2,
	\end{aligned}
\end{equation*}
where the constants $Q^{(6)}:= \left(Q^*e^{Q^*}Q_1+e^{Q^*-\epsilon}\right)^2$, $Q^{(7)}:=\left(Q_3Q^*e^{Q^*}\right)^2$ and $Q^{(8)}:=\left(Q_2Q^*e^{Q^*}\right)^2$. Here, we have also used the following estimate for $\|e^{-\check{c}^n} - e^{-c_h^n}\|_{\rm M}$:
\begin{equation}\label{eq:ech}
    \|e^{-\check{c}^n} - e^{-c_h^n}\|_{\rm M}^2=h^2\sum_{i=1}^{M} \sum_{j=1}^{M} |e^{-\check{c}^n_{i,j}} - e^{-c^n_{i,j}}|^2= h^2\sum_{i=1}^{M} \sum_{j=1}^{M} |e^{-\hat{\eta}}|^2 |\widetilde{\theta}_c^n|^2  \leq e^{-\varepsilon}\|\widetilde{\theta}_c^n\|_{\rm M}^2,
\end{equation}
where $\hat{\eta}$ lies between $\check{c}^n_{i,j}$ and $c^n_{i,j}$.

To estimate the second term $\mal{J}_2$, we first note that
\begin{equation*} 
\begin{aligned}
   	[d_x(e^{-c_h^n}\widetilde{\theta}_u^{n+1})]_{i+1/2,j} %&= \frac{1}{h}(e^{-c^n_{i+1,j}}[\widetilde{\theta}_u^{n+1}]_{i+1,j} - e^{-c^n_{i,j}}[\widetilde{\theta}_u^{n+1}]_{i,j})\\
    &= e^{-c^n_{i+1,j}}[d_x\widetilde{\theta}_u^{n+1}]_{i+1/2,j} + [d_x e^{-c_h^n}]_{i+1/2,j}[\widetilde{\theta}_u^{n+1}]_{i,j}. 
\end{aligned}
\end{equation*}
Thus, using the established discrete $W^{1,\infty}$ bound \eqref{infty:cn}, the separation estimate \eqref{eq:bound:cn} and the mean value theorem \eqref{eq:Dxe}, we have
\begin{equation}\label{eq:ec:point:x}
	\begin{aligned}
		&[\mal{A}_x e^{c_h^n}]_{i+1/2,j}\, [d_x(e^{-c_h^n}\widetilde{\theta}_u^{n+1})]_{i+1/2,j}\, [d_x\widetilde{\theta}_u^{n+1}]_{i+1/2,j}\\
		&\quad= [\mal{A}_x e^{c_h^n}]_{i+1/2,j}\,e^{-c^n_{i+1,j}}\,\bigl|[d_x\widetilde{\theta}_u^{n+1}]_{i+1/2,j}\bigr|^2 \\
        &\qquad+ [\mal{A}_x e^{c_h^n}]_{i+1/2,j}\,[d_xe^{-c_h^n}]_{i+1/2,j}\,[\widetilde{\theta}_u^{n+1}]_{i,j}\,[d_x\widetilde{\theta}_u^{n+1}]_{i+1/2,j}\\
		&\quad \geq e^{\frac{\varepsilon}{2}-\widetilde{Q}} 
          \bigl|[d_x\widetilde{\theta}_u^{n+1}]_{i+1/2,j}\bigr|^2 
        + [\mal{A}_x e^{c_h^n}]_{i+1/2,j}\, e^{-\eta} \,[d_x c_h^n]_{i+1/2,j}\,[\widetilde{\theta}_u^{n+1}]_{i,j}\,[d_x\widetilde{\theta}_u^{n+1}]_{i+1/2,j}\\
		% &\geq e^{\frac{\varepsilon}{2}-\widetilde{Q}} \bigl|[d_x\widetilde{\theta}_u^{n+1}]_{i+1/2,j}\bigr|^2 - \frac{1}{4\delta}\bigl|[\mal{A}_x e^{c_h^n}]_{i+1/2,j}\,e^{-\eta} [d_x c_h^n]_{i+1/2,j}\,[\widetilde{\theta}_u^{n+1}]_{i,j}\bigr|^2 - \delta\, \bigl|[d_x\widetilde{\theta}_u^{n+1}]_{i+1/2,j}\bigr|^2\\
		&\quad \geq \bigl(e^{\frac{\varepsilon}{2}-\widetilde{Q}} - \delta \bigr) \bigl|[d_x\widetilde{\theta}_u^{n+1}]_{i+1/2,j}\bigr|^2 - \frac{Q^{(9)}}{4\delta}\bigl|[\widetilde{\theta}_u^{n+1}]_{i,j}\bigr|^2,
	\end{aligned}
\end{equation}
where $Q^{(9)}:=(\widetilde{Q})^2e^{2\widetilde{Q}-\varepsilon}$.
Similarly, an analogous argument in the $y$-direction yields
\begin{equation}\label{eq:ec:point:y}
    \begin{aligned}
    &[\mal{A}_y e^{c_h^n}]_{i,j+1/2}\,[d_y(e^{-c_h^n}\widetilde{\theta}_u^{n+1})]_{i,j+1/2}\,[d_y\widetilde{\theta}_u^{n+1}]_{i,j+1/2} \\
    &\quad\geq 
    \bigl(e^{\frac{\varepsilon}{2}-\widetilde{Q}} - \delta \bigr) \bigl|[d_y\widetilde{\theta}_u^{n+1}]_{i,j+1/2}\bigr|^2 - \frac{Q^{(9)}}{4\delta}\bigl|[\widetilde{\theta}_u^{n+1}]_{i,j}\bigr|^2.
\end{aligned}
\end{equation}
Therefore, \eqref{eq:ec:point:x}--\eqref{eq:ec:point:y} directly imply that
\begin{equation*}%\label{eq:J2}
	\begin{aligned}
		\mal{J}_2 \geq \bigl(e^{\frac{\varepsilon}{2}-\widetilde{Q}} - \delta \bigr)  \|\bm{d}\widetilde{\theta}_u^{n+1}\|_{\rm TM}^2 - \frac{Q^{(9)}}{2\delta}\|\widetilde{\theta}_u^{n+1}\|_{\rm M}^2.
	\end{aligned}
\end{equation*}

Finally, for the remaining term $\mal{J}_3$, applying the established discrete $W^{1,\infty}$ bounds \eqref{eq:regular:2} and \eqref{infty:cn}, the separation estimate \eqref{eq:bound:cn} and an estimate analogous to \eqref{eq:ech}, we have
\begin{equation*}
	\begin{aligned}
		\mal{J}_3 %&= \bigl ( [\mal{A} e^{\check{c}^n}-\mal{A} e^{c_h^n}] [\bm{d}(e^{-c_h^n}\check{u}^{n+1})],\bm{d} \widetilde{\theta}_u^{n+1} \bigr )_{\rm TM}\\
		&\leq \|\bm{d}(e^{-c_h^n}\check{u}^{n+1})\|_{\infty} \|\mal{A} e^{\check{c}^n}-\mal{A} e^{c_h^n}\|_{\rm TM} \|\bm{d} \widetilde{\theta}_u^{n+1}\|_{\rm TM}\\
		&\leq \bigl( \|e^{-c_h^n}\|_{\infty} \|\bm{d}\check{u}^{n+1}\|_{\infty} + \|\bm{d} e^{-c_h^n}\|_{\infty} \|\check{u}^{n+1}\|_{\infty}\bigr)\| e^{\check{c}^n}- e^{c_h^n}\|_{\rm M} \|\bm{d} \widetilde{\theta}_u^{n+1}\|_{\rm TM}\\
		&\leq (Q^*+\widetilde{Q})e^{-\f{\varepsilon}{2}+\max\{Q^*,\widetilde{Q}\}}\|\widetilde{\theta}_c^{n}\|_{\rm M} \|\bm{d} \widetilde{\theta}_u^{n+1}\|_{\rm TM}
        \leq \frac{Q^{(10)}}{\delta}\|\widetilde{\theta}_c^{n}\|_{\rm M}^2 + \delta\|\bm{d} \widetilde{\theta}_u^{n+1}\|_{\rm TM}^2,
	\end{aligned}
\end{equation*}
where $Q^{(10)}:=(Q^*+\widetilde{Q})e^{-\f{\varepsilon}{2}+\max\{Q^*,\widetilde{Q}\}}$.
 
Therefore, substituting the estimates for $\mal{J}_1$--$\mal{J}_3$ into \eqref{error:u}, we obtain 
\begin{equation}\label{error:u:estimate}
	\begin{aligned}
		&\frac{1}{2\tau} \bigl( \|\widetilde{\theta}_u^{n+1}\|_{\rm M}^2 - \|\widetilde{\theta}_u^n\|_{\rm M}^2 +\|\widetilde{\theta}_u^{n+1}-\widetilde{\theta}_u^n\|_{\rm M}^2  \bigr) + \bigl(e^{\frac{\varepsilon}{2}-\widetilde{Q}} -5\delta\bigr)\|\bm{d}\widetilde{\theta}_u^{n+1}\|_{\rm TM}^2\\
		 &\quad\leq  \Bigl(\frac{Q^{(9)}}{2\delta}+\frac{1}{2}\Bigr)\|\widetilde{\theta}_u^{n+1}\|_{\rm M}^2 
         + \Bigl(\frac{Q^{(6)}}{4\delta}+\frac{Q^{(10)}}{\delta}\Bigr)\|\widetilde{\theta}_c^{n}\|_{\rm M}^2
         + \frac{Q^{(7)}}{4\delta}\|\bm{d}\widetilde{\theta}_c^{n}\|_{\rm TM}^2 + \frac{1}{2} \|R_u^{n+1}\|_{\rm M}^2 + \frac{Q^{(8)}}{4\delta}h^8.
	\end{aligned}
\end{equation}
Proceeding analogously to the analysis of $\widetilde{\theta}_u^{n+1}$,  we obtain the estimate for $\widetilde{\theta}_v^{n+1}$ as follows:
\begin{equation}\label{error:v:estimate}
	\begin{aligned}
		&\frac{1}{2\tau} \bigl( \|\widetilde{\theta}_v^{n+1}\|_{\rm M}^2 - \|\widetilde{\theta}_v^n\|_{\rm M}^2 +\|\widetilde{\theta}_v^{n+1}-\widetilde{\theta}_v^n\|_{\rm M}^2  \bigr) + \bigl(e^{\frac{\varepsilon}{2}-\widetilde{Q}} -5\delta\bigr) \|\bm{d}\widetilde{\theta}_v^{n+1}\|_{\rm TM}^2\\
		&\quad\leq \Bigl(\frac{Q^{(9)}}{2\delta}+\frac{1}{2}\Bigr)\|\widetilde{\theta}_v^{n+1}\|_{\rm M}^2 
        + \Bigl(\frac{Q^{(6)}}{4\delta}+\frac{Q^{(10)}}{\delta}\Bigr)\|\widetilde{\theta}_c^{n}\|_{\rm M}^2
        + \frac{Q^{(7)}}{4\delta}\|\bm{d}\widetilde{\theta}_c^{n}\|_{\rm TM}^2+ \frac{1}{2} \|R_v^{n+1}\|_{\rm M}^2 + \frac{Q^{(8)}}{4\delta}h^8.
	\end{aligned}
\end{equation}

In addition, taking the discrete inner product of \eqref{error:function_c} with $\widetilde{\theta}_c^{n+1}$ and with $-D_x(d_x\widetilde{\theta}_c^{n+1})-D_y(d_y\widetilde{\theta}_c^{n+1})$, respectively, and applying Lemma \ref{lemma:Dd}, we obtain the estimate for $\widetilde{\theta}_c^{n+1}$ as follows:
\begin{equation}\label{error:c:estimate}
\begin{aligned}
		&\frac{1}{2\tau} \bigl( \|\widetilde{\theta}_c^{n+1}\|_{\rm M}^2 - \|\widetilde{\theta}_c^n\|_{\rm M}^2 +\|\widetilde{\theta}_c^{n+1}-\widetilde{\theta}_c^n\|_{\rm M}^2  \bigr) \\
		&\quad= -\|\bm{d}\widetilde{\theta}_c^{n+1}\|_{\rm M}^2 - \|\widetilde{\theta}_c^{n+1}\|_{\rm M}^2 + \bigl( \widetilde{\theta}_u^{n+1},\widetilde{\theta}_c^{n+1} \bigr)_{\rm M}+ \bigl( \widetilde{\theta}_v^{n+1},\widetilde{\theta}_c^{n+1} \bigr)_{\rm M} + \bigl( R_c^{n+1},\widetilde{\theta}_c^{n+1} \bigr)_{\rm M}\\
        &\quad \leq -\|\bm{d}\widetilde{\theta}_c^{n+1}\|_{\rm M}^2 +\|\widetilde{\theta}_u^{n+1}\|_{\rm M}^2 + \|\widetilde{\theta}_v^{n+1}\|_{\rm M}^2 + \frac{1}{2}\|R_c^{n+1}\|_{\rm M}^2,
	\end{aligned}
\end{equation}
and
\begin{equation}\label{error:gradc:estimate}
	\begin{aligned}
		&\frac{1}{2\tau} \bigl( \|\bm{d}\widetilde{\theta}_c^{n+1}\|_{\rm TM}^2 - \|\bm{d}\widetilde{\theta}_c^n\|_{\rm TM}^2 
        +\|\bm{d} (\widetilde{\theta}_c^{n+1}-\widetilde{\theta}_c^n)\|_{\rm TM}^2  \bigr) \\
		% &\quad= -\|D_x(d_x\widetilde{\theta}_c^{n+1})+D_y(d_y\widetilde{\theta}_c^{n+1})\|_{\rm M}^2 - \|\bm{d}\widetilde{\theta}_c^{n+1}\|_{\rm TM}^2+\bigl( \widetilde{\theta}_u^{n+1},-D_x(d_x\widetilde{\theta}_c^{n+1})-D_y(d_y\widetilde{\theta}_c^{n+1})\bigr)_{\rm M}\\
  %       &\qquad+ \bigl( \widetilde{\theta}_v^{n+1},-D_x(d_x\widetilde{\theta}_c^{n+1})-D_y(d_y\widetilde{\theta}_c^{n+1}) \bigr)_{\rm M} + \bigl( R_c^{n+1},-D_x(d_x\widetilde{\theta}_c^{n+1})-D_y(d_y\widetilde{\theta}_c^{n+1}) \bigr)_{\rm M}\\
		&\quad \leq -\|\bm{d} \widetilde{\theta}_c^{n+1}\|_{\rm TM}^2 + \|\widetilde{\theta}_u^{n+1}\|_{\rm M}^2 + \|\widetilde{\theta}_v^{n+1}\|_{\rm M}^2 
        + \frac{1}{2}\|R_c^{n+1}\|_{\rm M}^2.
	\end{aligned}
\end{equation}

Denote
\[
\Theta^{n}:=  \|\widetilde{\theta}_u^n\|_{\rm M}^2 + \|\widetilde{\theta}_v^n\|_{\rm M}^2 +  \|\widetilde{\theta}_c^n\|_{\rm M}^2 + \|\bm{d}\widetilde{\theta}_c^n\|_{\rm TM}^2.
\]
Then, combining \eqref{error:u:estimate}--\eqref{error:gradc:estimate} and choosing $\delta \le \delta^*:=\frac{e^{\frac{\varepsilon}{2}-\widetilde{Q}}}{5}$ so as to discard the nonnegative terms on the left-hand side, we obtain
\begin{equation}\label{sum:uvc}
	\begin{aligned}
		\Theta^{n+1} -\Theta^{n} 
     %   + 2\tau\bigl(e^{\frac{\varepsilon}{2}-\widetilde{Q}} -5\delta\bigr)
   %     [\|\bm{d}\widetilde{\theta}_u^{n+1}\|_{\rm TM}^2 + \|\bm{d}\widetilde{\theta}_v^{n+1}\|_{\rm TM}^2]\\
		& \leq 2\tau\left( \frac{Q^{(9)}}{2\delta}+\frac{5}{2}\right)\left(\|\widetilde{\theta}_u^{n+1}\|_{\rm M}^2+\|\widetilde{\theta}_v^{n+1}\|_{\rm M}^2\right) + \frac{Q^{(7)}}{\delta}\tau\|\bm{d}\widetilde{\theta}_c^{n+1}\|_{\rm TM}^2+ \frac{Q^{(8)}}{\delta}\tau h^8\\
		&\qquad + 4\tau\left( \frac{Q^{(6)}}{4\delta}+\frac{Q^{(10)}}{\delta}\right)\|\widetilde{\theta}_c^{n+1}\|_{\rm M}^2 
        + \tau \left(\|R_u^{n+1}\|_{\rm M}^2 + \|R_v^{n+1}\|_{\rm M}^2+ 2\|R_c^{n+1}\|_{\rm M}^2\right)\\
       &\leq \hat{Q} \tau \Theta^{n+1}+ \tau \left(\f{Q^{(8)}}{\delta}+4\check{Q}^2\right) (\tau^3+h^4)^2,
	\end{aligned}
\end{equation}
where %$\delta=\frac{e^{\frac{\varepsilon}{2}-\widetilde{Q}}}{6}$ is chosen to ensure the nonnegative of the coefficient associated with the discrete gradient term on the left-hand side. 
$\hat{Q}:=\max\{ \f{Q^{(9)}}{\delta}+5,\f{Q^{(6)}+4Q^{(10)}}{\delta},\f{Q^{(7)}}{\delta}\}$ and the higher-order consistency errors in Lemma \ref{lem:high-order} are applied in the last step.
Therefore, an application of the discrete Gr\"{o}nwall's inequality to \eqref{sum:uvc} directly yields 
\begin{equation}\label{result}
	\|\widetilde{\theta}_u^{n+1}\|_{\rm M} + \|\widetilde{\theta}_v^{n+1}\|_{\rm M} + \|\widetilde{\theta}_c^{n+1}\|_{\rm M} + \|\bm{d}\widetilde{\theta}_c^{n+1}\|_{\rm TM} 
    \leq Q_4 \bigl(\tau^3+h^4\bigr),\quad \text{for} ~\tau \leq \tau^*:=1/(2\hat{Q}),
\end{equation}
where $Q_4:=2e^{\hat{Q}T}\sqrt{2T(\f{Q^{(8)}}{\delta}+4\check{Q}^2)}$.

\paragraph{\indent \bf Part III.  Completion of the inductive assumption \eqref{priori_assum} and final conclusions}
Invoking the higher-order error estimate \eqref{result}, for sufficiently small $\tau$ and $h$, we obtain
\begin{equation*}
	\|\widetilde{\theta}_c^{n+1}\|_{\rm M}\leq \tau^{\frac{9}{4}} + h^{\frac{9}{4}}.
\end{equation*}
This completes the a priori assumption \eqref{priori_assum} at the next time level $t_{n+1}$, thereby closing the induction and establishing the higher-order error estimate for all time levels under consideration.

To establish the desired final estimates in Theorem \ref{thm:covergence}, it suffices to derive bounds for the following terms:
\begin{equation}\label{eq:error:final}
    \begin{aligned}
        \|z^n-z_h^n\|_{\rm M} &\leq \|z^n- z_M^n\|_{\rm M} + \|\check{z}^n-z_M^n\|_{\rm M} + \|\widetilde{\theta}_z^n\|_{\rm M}, \quad z=\{u,v,c\},\\
        \|\nabla c^n-\bm{d}c_h^n\|_{\rm TM} 
        &\leq \|\nabla c^n-\bm{d} c^n\|_{\rm TM} 
        + \|\bm{d} (c^n-c_M^n)\|_{\rm TM}
        + \|\bm{d}(\check{c}^n-c_M^n)\|_{\rm TM}
        +\|\bm{d} \widetilde{\theta}_c^n\|_{\rm TM}.
    \end{aligned}
\end{equation}
Note that the projection approximation \eqref{projection:error}, the definition of the corrected solutions in \eqref{check_U}, and the regularity assumption \eqref{eq:regular} together imply that
\begin{equation}\label{eq:error:pre}
    \begin{aligned}
       \|z^n- z_M^n\|_{\rm M} + \|\check{z}^n-z_M^n\|_{\rm M} &\leq Q(\tau+h^2), \quad z=\{u,v,c\},\\
        \|\nabla c^n-\bm{d} c^n\|_{\rm TM} + \|\bm{d} (c^n-c_M^n)\|_{\rm TM}
        + \|\bm{d}(\check{c}^n-c_M^n)\|_{\rm TM} &\leq Q(\tau+h^2).
    \end{aligned}
\end{equation}
Therefore, combining the estimates \eqref{result}--\eqref{eq:error:pre}, we  obtain the final convergence result:
\begin{equation*}
		\|u^n -  u_h^n\|_{\rm M} + \|v^n -  v_h^n\|_{\rm M} + \|c^n -  c_h^n\|_{\rm M} 
          + \|\nabla c^n-\bm{d} c_h^n\|_{\rm TM} \leq Q_*(\tau + h^2),\quad 1 \leq n \leq N.
\end{equation*}
This completes the proof of Theorem \ref{thm:covergence}. \qed

% \begin{remark}
%     The terms involving $\|\check{z}^n-z_M^n\|_{\rm M}$ and $\|\bm{d}(\check{c}^n-c_M^n)\|_{\rm TM}$ on the left-hand side of \eqref{eq:error:pre} can be directly estimated using the definition of the corrected solutions. More precisely, for $z=\{u,v,c\}$, we have 
%     \begin{align*}
%         \|\check{z}^n-z_M^n\|_{\rm M} &\leq \|\check{z}^n-z_M^n\|_{\infty}
%         =\|\tau\mal{P}_Mz_{\tau,1}^n + \tau^2\mal{P}_Mz_{\tau,2}^n + h^2 \mal{P}_M z_{h,1}^n\|_{\infty} \leq Q(\tau+h^2),\\
%         \|\bm{d}\check{c}^n-\bm{d}c_M^n\|_{\rm TM} &\leq \|\bm{d}\check{c}^n-\bm{d}c_M^n\|_{\infty}=\|\tau\bm{d}(\mal{P}_Mc_{\tau,1}^n) + \tau^2\bm{d}(\mal{P}_Mc_{\tau,2}^n) + h^2 \bm{d}(\mal{P}_M c_{h,1}^n)\|_{\infty} \leq Q(\tau+h^2),
%     \end{align*}
%     due to the regularity assumption \eqref{eq:regular}.
% \end{remark}

\begin{remark}
   The higher-order convergence estimate $\mathcal{O}(\tau^3+h^4)$ in \eqref{result} is not a superconvergence result for the proposed numerical scheme \eqref{model:KS:shemeI}--\eqref{model:KS:shemeII}; rather, it is an auxiliary error estimate established directly for the corrected solutions $\{\check{u},\check{v},\check{c}\}$. Combining this estimate with \eqref{eq:error:final}--\eqref{eq:error:pre}, we then recover the optimal-order error estimates stated in Theorem \ref{thm:covergence}. 
   
   The error analysis can be carried out along the same lines for the parabolic-elliptic case (i.e., $\epsilon=0$, $\alpha>0$), with the main modification being that the a priori assumption \eqref{priori_assum} is replaced by  $\|\tilde{\theta}_u^{n}\|_{\rm M}+\|\tilde{\theta}_v^{n}\|_{\rm M} \leq \tau^{\f{9}{4}}+h^{\f{9}{4}}$. Under this modified assumption, the convergence estimates stated in Theorem \ref{thm:covergence} follow by adapting the preceding analysis.
\end{remark}

%%%%%%%%%%%%%%%%%%%%%%%%%%%%%%%%%%%%%%%%%%%%%%%%%%%%
\section{Numerical results}\label{sec:num}
In this section, we present several numerical experiments to assess the accuracy and structure-preserving properties of the proposed scheme \eqref{model:KS:shemeI}--\eqref{model:KS:shemeII} for both the parabolic-parabolic case (i.e., $\epsilon=1$) and the parabolic-elliptic case (i.e., $\epsilon=0$) of model \eqref{model:KS}--\eqref{model:KS:ic}.

\subsection{Accuracy test}
\begin{example}\label{exam:s1}
In this example, we set the computational domain $\Omega=(-1,1)^2$ and assess the convergence rate of the proposed scheme \eqref{model:KS:shemeI}--\eqref{model:KS:shemeII} for the general two-species Keller–Segel chemotaxis system: 
\begin{equation*}
\left\{
	\begin{aligned}
		\p_t u  &= \Delta u- \nabla \cdot(u \nabla c) + f_1,\\
		\p_t v  &= \Delta v- \nabla \cdot(v \nabla c) + f_2,\\
		\epsilon\, \p_t c  &=  \Delta c -  c+ u  +  v + g,
	\end{aligned}
    \right.
\end{equation*}
where the source terms $(f_1,f_2,g)$ are given such that the exact solutions are prescribed as
\begin{equation*}
	\begin{aligned}
		u(\bm x, t) &= \sin(\pi x)\sin(\pi y)\sin(t)+1.1,\\
		v(\bm x, t) &= \f{1}{2\pi^2+1} \sin(\pi x)\sin(\pi y)\sin(t)+1.1,\\
		c(\bm x,t) &= \sin(\pi x)\sin(\pi y)\sin(t) + 1.1.
	\end{aligned}
\end{equation*}
\end{example}
%The periodic boundary conditions are assumed, and the source terms $(f_1,f_2,g)$ can be computed from the above solution.

The temporal and spatial convergence orders of the proposed PP-MC-ED-BCFD scheme are investigated with $\tau = h^2$. Table \ref{tab:accuracy:I} presents the error results for the parabolic-parabolic case (i.e., $\epsilon=1$), while  the parabolic-elliptic case (i.e., $\epsilon=0$) is reported in Table \ref{tab:accuracy:II}. Both cases confirm that the proposed scheme achieves first-order accuracy in time and second-order accuracy in space in the relevant norms, in agreement with the theoretical conclusions established in Theorem \ref{thm:covergence}.
\begin{table}[t]
	\centering \caption{ $L^2$-norm errors of $u$, $v$, $c$ and $\nabla c$ for Example \ref{exam:s1} with $\epsilon=1$. } \label{tab:accuracy:I}
	\setlength{\tabcolsep}{1.5mm}
	\begin{tabular}{c c c c c c c c c }
		\toprule
		$M$ &  $\|u-u_h\|_{\rm M}$&  Order&  $\|v-v_h\|_{\rm M}$&  Order& $\|c-c_h\|_{\rm M}$& Order &  $\|\nabla c-\bm{d} c_h\|_{\rm TM}$& Order  \\
		\midrule
		10 & 3.09e-02 & --    & 1.27e-02 & --    & 2.36e-02 &--&4.20e-02 & --\\
		20 & 7.02e-03 & 2.14  & 2.61e-03 & 2.28  & 5.87e-03 & 2.00 &1.07e-02 & 1.98\\
		40 & 1.71e-03 & 2.03  & 6.18e-04 & 2.08  & 1.47e-03 & 2.00 &2.68e-03 & 1.99\\
		80 & 4.26e-04 & 2.01  & 1.52e-04 & 2.02  & 3.67e-04 & 2.00 &6.71e-04 & 2.00\\
		\bottomrule
	\end{tabular}
\end{table}
\begin{table}[t]
	\centering \caption{ $L^2$-norm errors of $u$, $v$, $c$ and $\nabla c$ for Example \ref{exam:s1} with $\epsilon=0$. } \label{tab:accuracy:II}
	\setlength{\tabcolsep}{1.5mm}
	\begin{tabular}{c c c c c c c c c }
		\toprule
		$M$ &  $\|u-u_h\|_{\rm M}$&  Order&  $\|v-v_h\|_{\rm M}$&  Order& $\|c-c_h\|_{\rm M}$& Order &  $\|\nabla c-\bm{d} c_h\|_{\rm TM}$& Order  \\
		\midrule
		10 & 2.88e-02 & --    & 1.01e-02 & --    & 2.59e-02 &--&5.21e-02 & --\\
		20 & 6.53e-03 & 2.14  & 2.01e-03 & 2.32  & 6.43e-03 & 2.01 &1.31e-02 & 1.99\\
		40 & 1.59e-03 & 2.03  & 4.73e-04 & 2.09  & 1.60e-03 & 2.00 &3.29e-03 & 2.00\\
		80 & 3.96e-04 & 2.01  & 1.16e-04 & 2.02  & 4.01e-04 & 2.00 &8.22e-04 & 2.00\\
		\bottomrule
	\end{tabular}
\end{table}

\subsection{Structure-preserving test}
    In this subsection, we examine the preservation of physical properties by the proposed scheme in simulations of both the parabolic–parabolic and parabolic–elliptic cases of the 2D two-species Keller-Segel model. As shown in Refs. \cite{HS'24,EVC'12,L'24}, global solutions to the two-species Keller-Segel chemotaxis system exist, provided that the total initial mass satisfies the critical relation $M[u^0]+M[v^0]<8\pi$.
        
	\begin{example}[2D parabolic-parabolic case]\label{exam:eg2:pp}
	   In this example, we consider the computational domain $\Omega=(0,1)^2$, the final time $T=10$, and prescribe the initial  densities and concentration as follows:
			\begin{align*}
				u^0(\bm x)&=12 \exp \left(-2\big((x-0.5)^2+(y-0.5)^2\big)\right),\\
				v^0(\bm x)&=10 \exp \left(-2\big((x-0.5)^2+(y-0.5)^2\big)\right),\\
				c^0(\bm x)&=8 \exp \left(-2\big((x-0.5)^2+(y-0.5)^2\big)\right) ,
			\end{align*}
		so that the total mass satisfies $M[u^0]+M[v^0]\approx  8.79+ 7.32=16.11<8\pi$, and hence global solutions exist. 
	\end{example}		

\begin{figure} [t]
	\centering
	\includegraphics[width=0.325\linewidth]{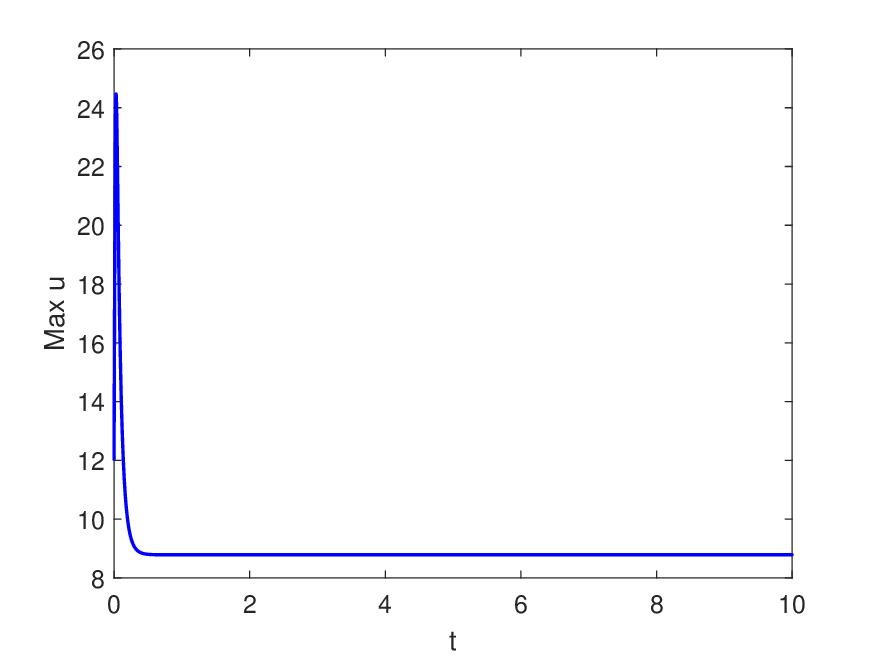}
	\includegraphics[width=0.325\linewidth]{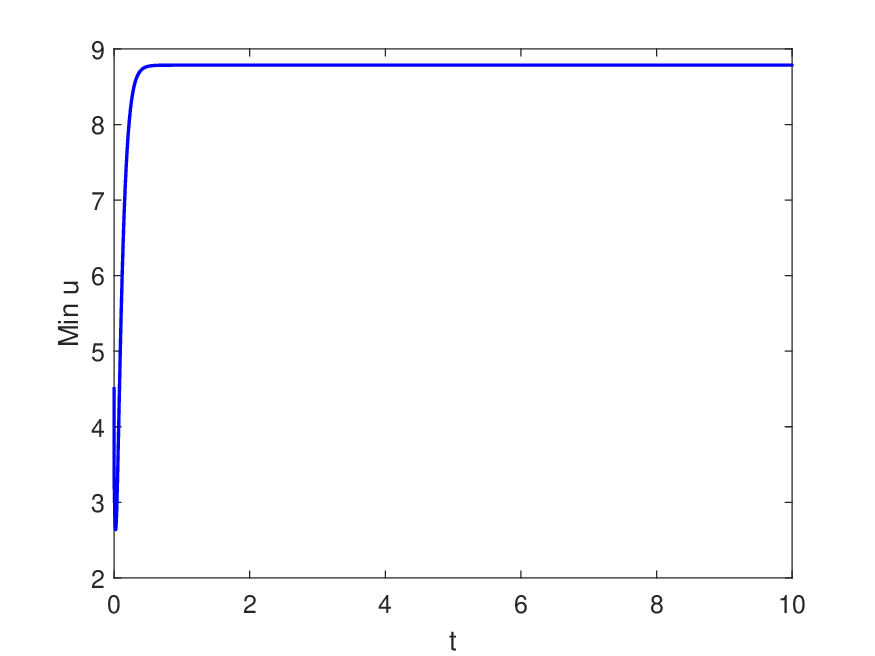}
	\includegraphics[width=0.325\linewidth]{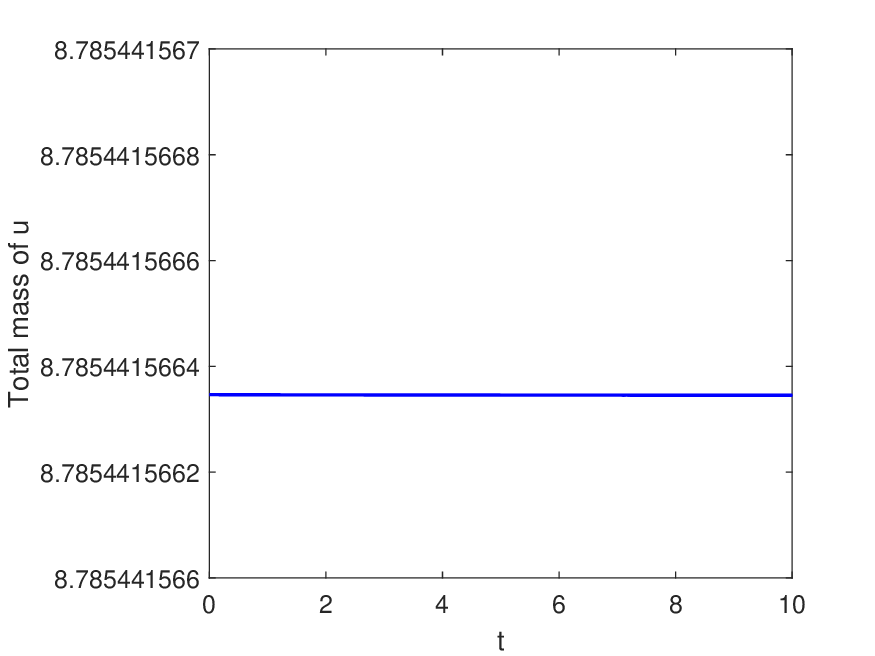}\\
	\includegraphics[width=0.325\linewidth]{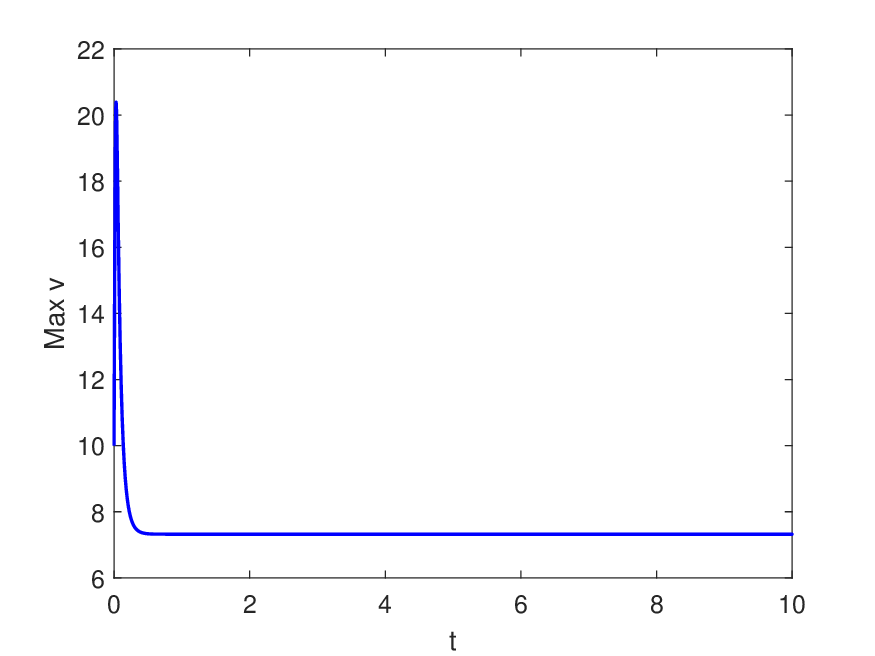}
	\includegraphics[width=0.325\linewidth]{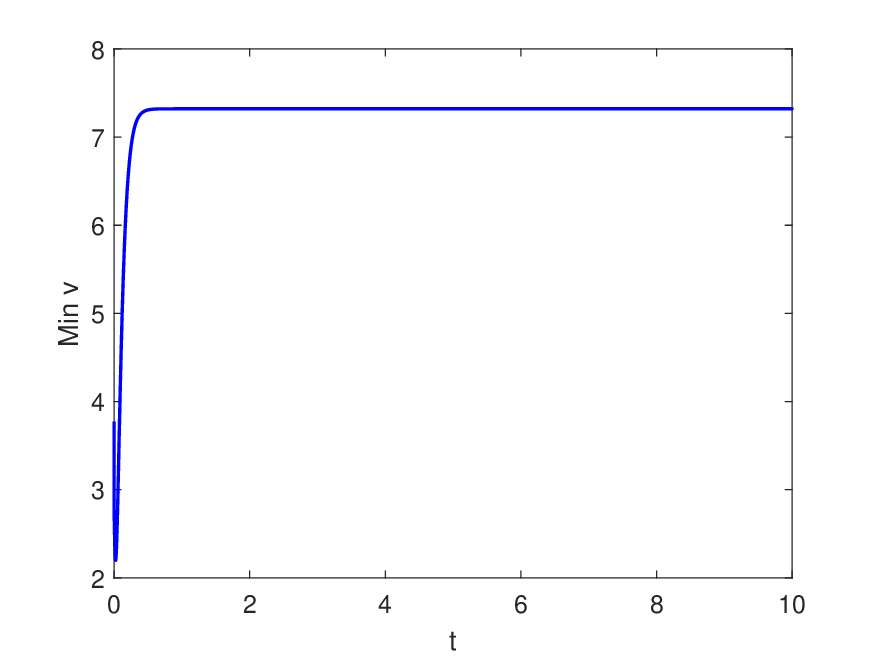}
	\includegraphics[width=0.325\linewidth]{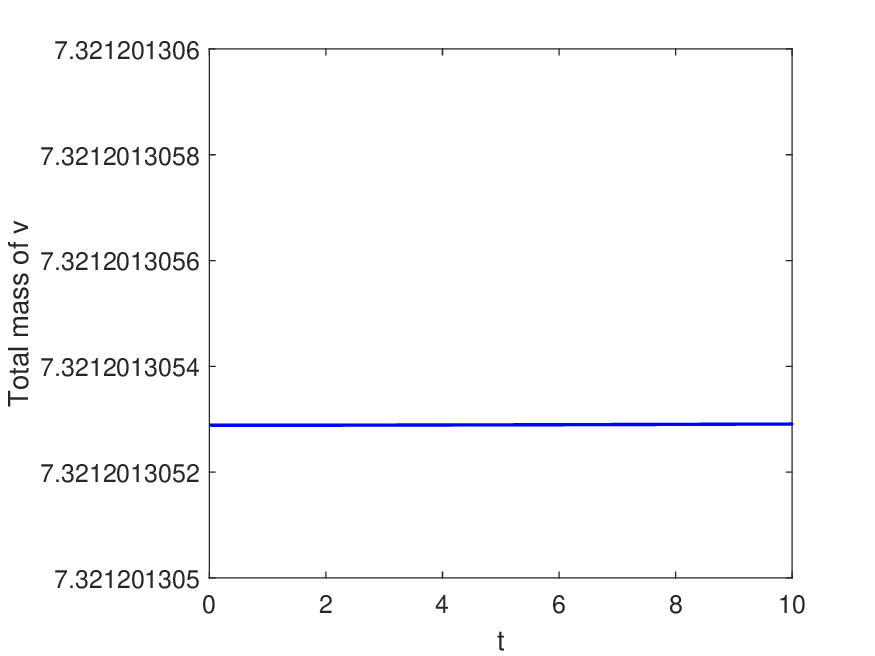}\\
	\includegraphics[width=0.325\linewidth]{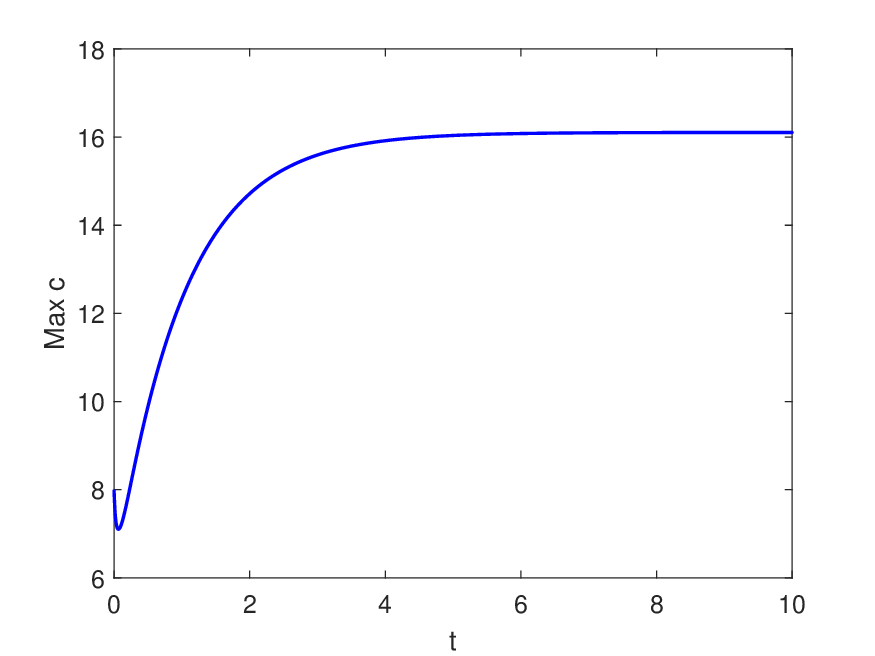}
	\includegraphics[width=0.325\linewidth]{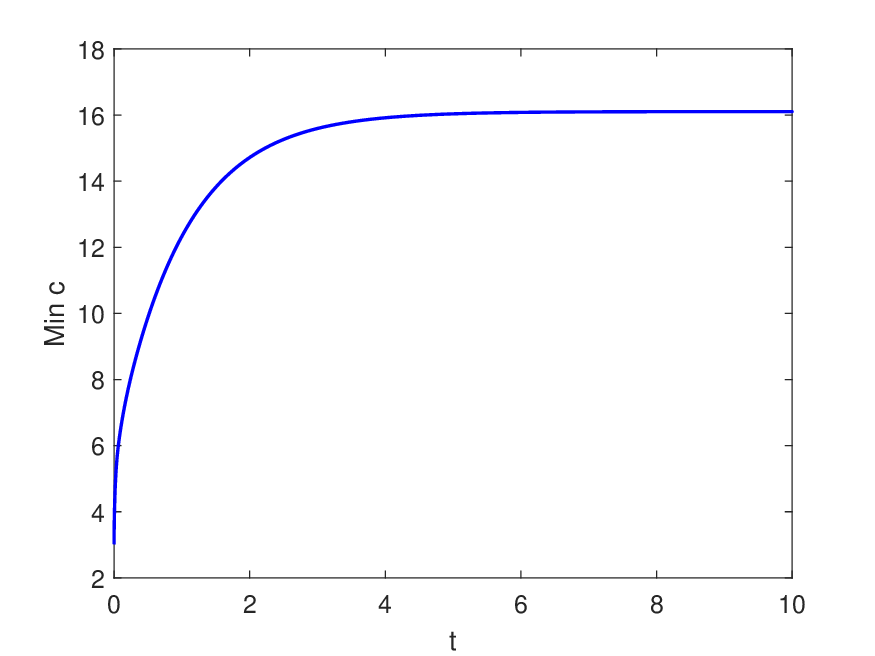}
	\includegraphics[width=0.325\linewidth]{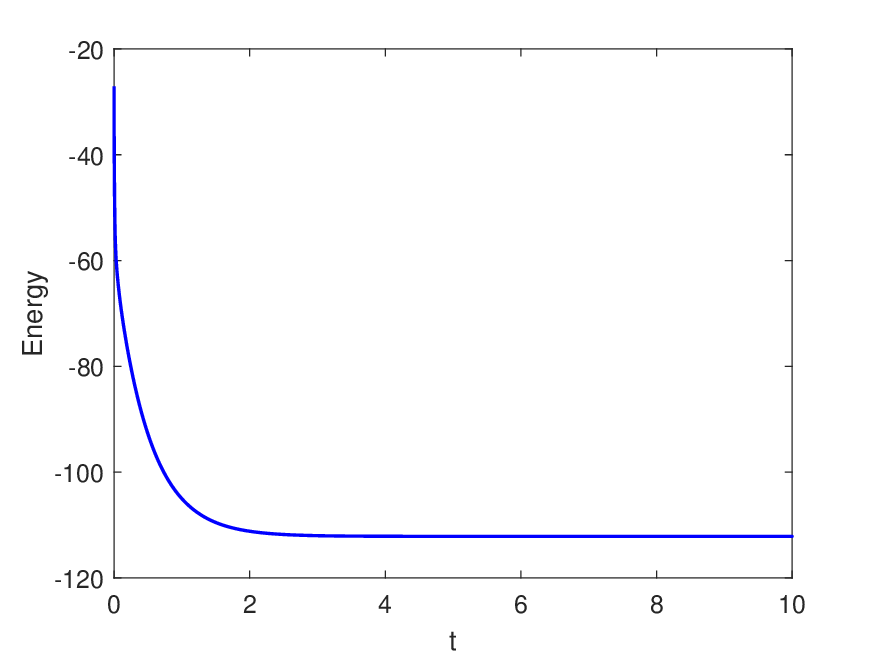}
	\caption{ Evolutions of the extrema of $u$, $v$, $c$, the total mass of $u,v$ and the original energy (from left to right) for Example \ref{exam:eg2:pp}. }
	\label{fig:eg2:mass}
\end{figure}

In the simulation, the computational domain is discretized with $M=80$ grid points in each spatial direction, and the time stepsize is set to $\tau=2.0\times10^{-3}$. The evolution of the maximum and minimum values of $u$, $v$, and $c$, together with the total masses of $u$ and $v$ and the discrete original energy, is presented in Fig. \ref{fig:eg2:mass}. The numerical observations are summarized as follows: (i) the maximum values of both $u$ and $v$ initially decrease, while the minimum values initially increase; thereafter, both quantities gradually approach their steady states; (ii) the densities $u$ and $v$ and the chemoattractant concentration $c$ remain positive, and the total masses of $u$ and $v$ are conserved to at least 11 significant figures throughout the the time evolution; and (iii) the discrete original energy defined in Theorem \ref{thm:structure-pre} is dissipative. These observations are fully consistent with the structure-preserving properties proven in Theorem \ref{thm:structure-pre}. 

\begin{example}[2D parabolic-elliptic case]\label{exam:eg3:pe}
    In this example, we investigate the structure-preserving properties of the proposed scheme for  the parabolic-elliptic case with $\Omega=(-2,2)^2$ and $T=10$. The initial conditions are prescribed as
	\begin{align*}
		u^0(\bm x)=\f{30}{1+40(x^2+y^2)},\quad v^0(\bm x)=\f{20}{ 1+40(x^2+y^2)},
	\end{align*}
	so that  $M[u^0]+M[v^0]\approx 12.49+ 8.33=20.82<8\pi$, and hence global solutions exist.
\end{example}

\begin{figure}[t]
	\centering
	\includegraphics[width=0.325\linewidth]{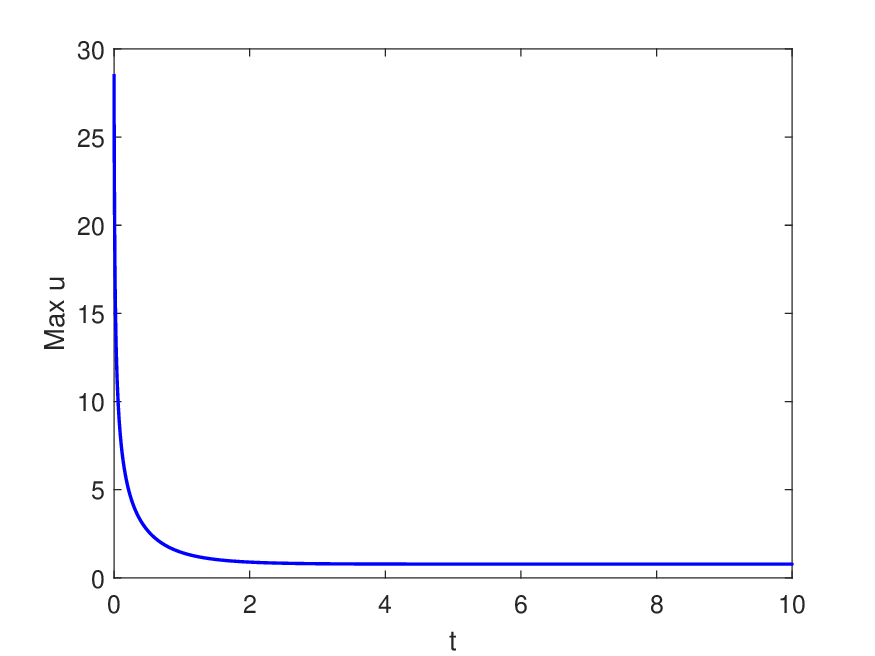}
	\includegraphics[width=0.325\linewidth]{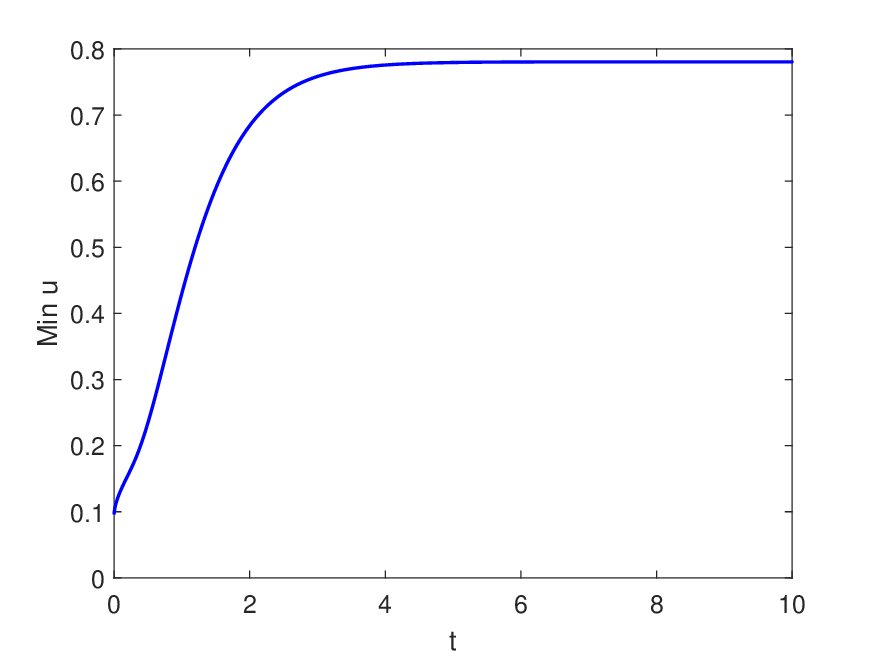}
	\includegraphics[width=0.325\linewidth]{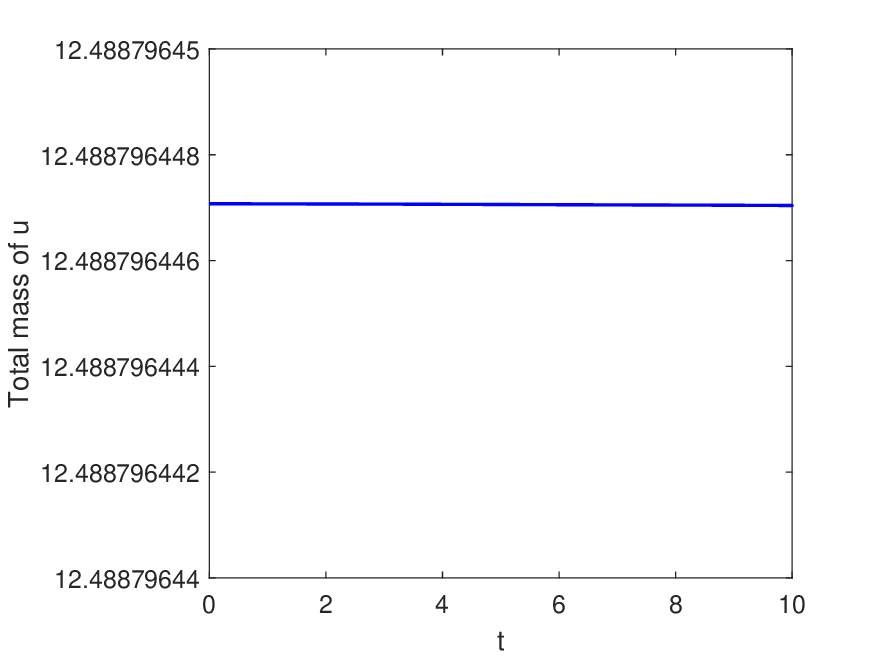}
	\includegraphics[width=0.325\linewidth]{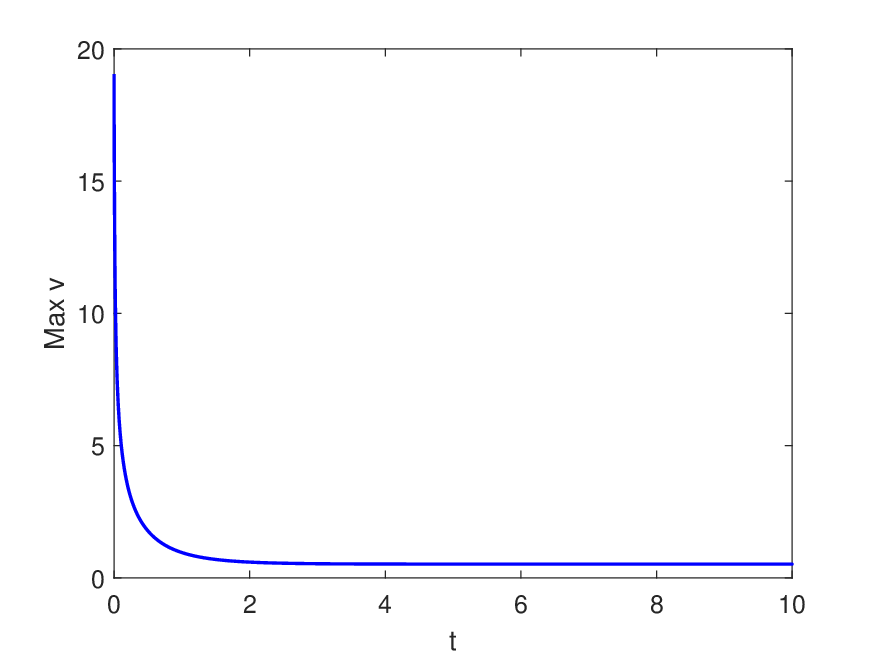}
	\includegraphics[width=0.325\linewidth]{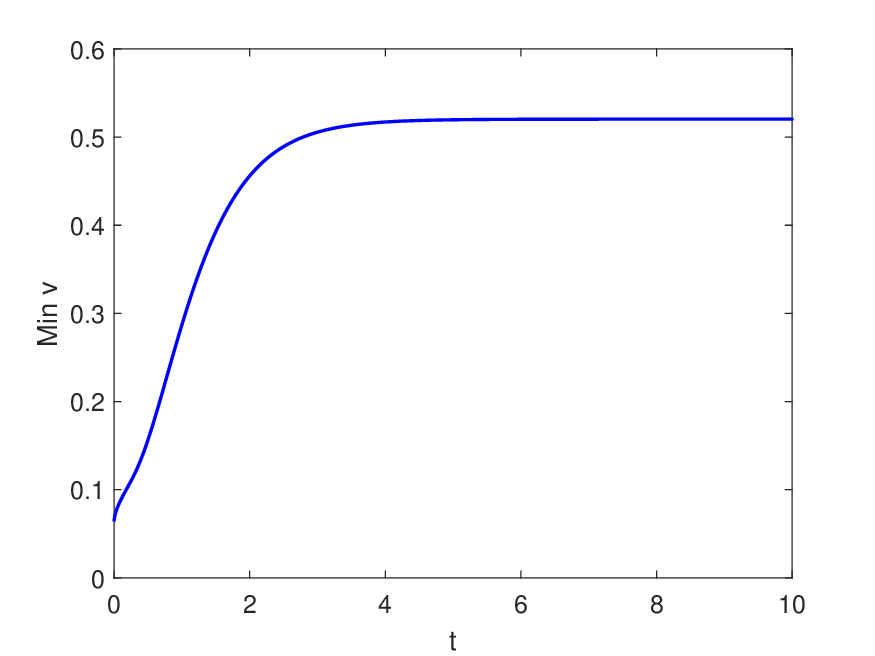}
	\includegraphics[width=0.325\linewidth]{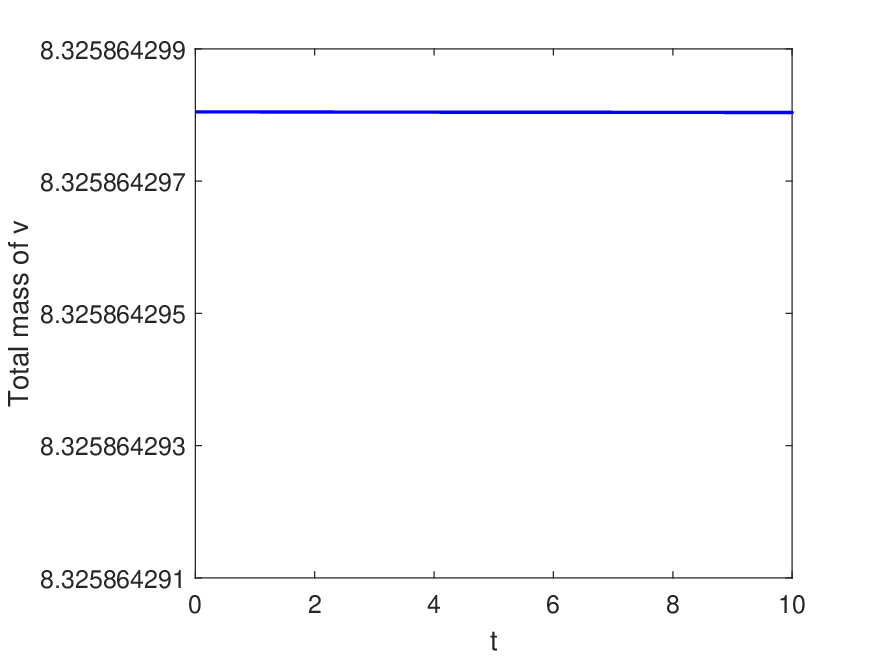}
	\includegraphics[width=0.325\linewidth]{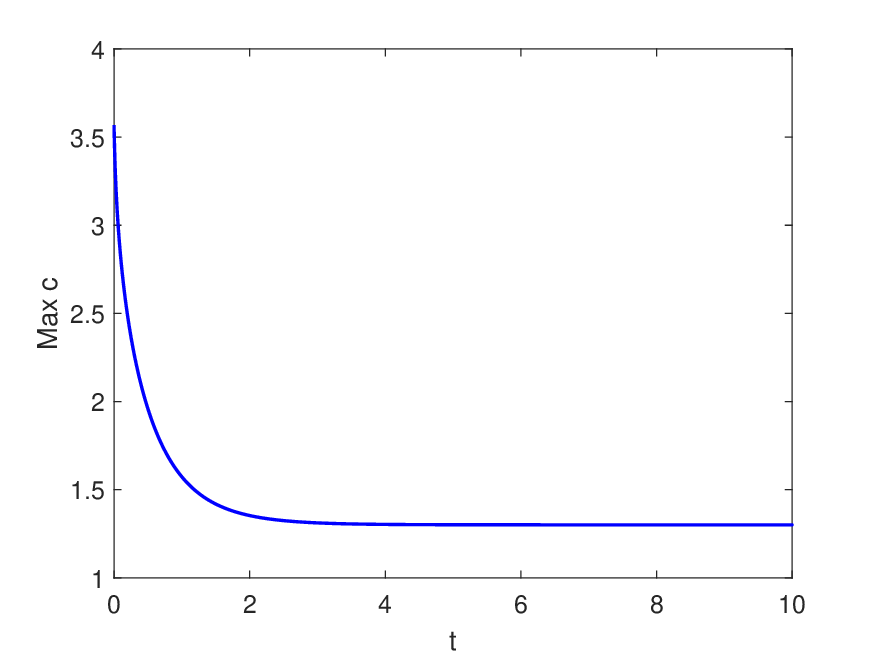}
	\includegraphics[width=0.325\linewidth]{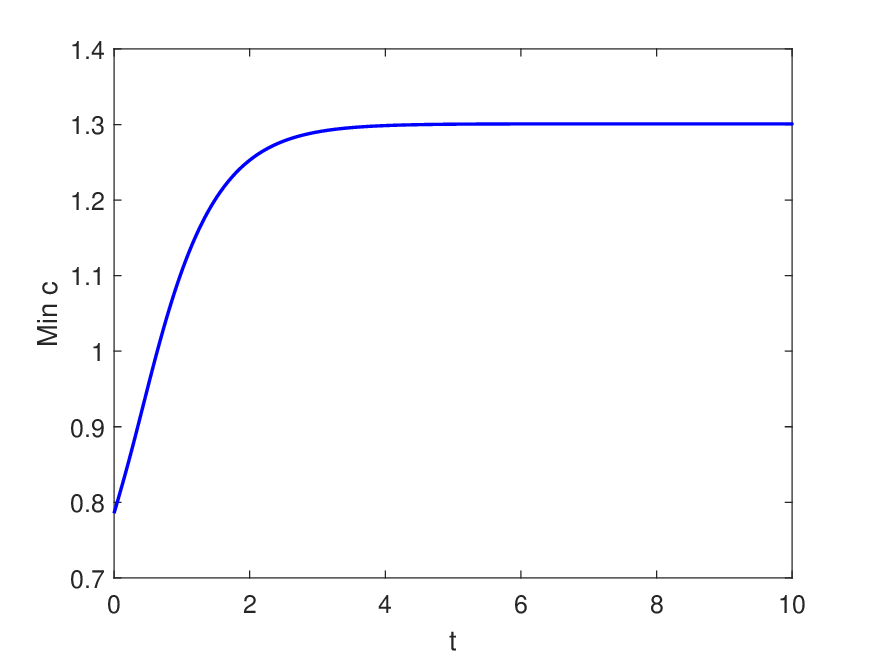}
	\includegraphics[width=0.325\linewidth]{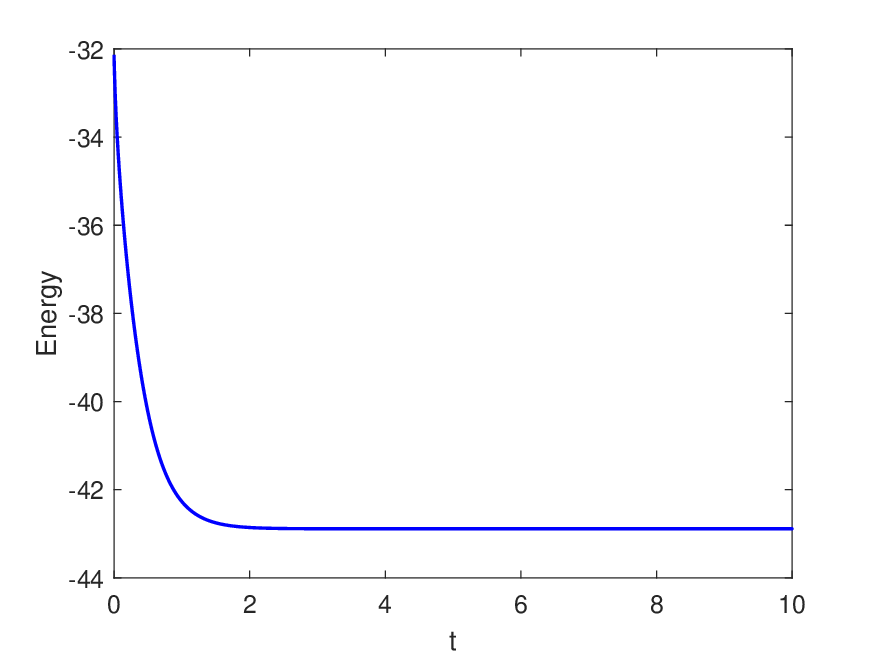}
	\caption{ Evolutions of the extrema of $u$, $v$, $c$, and the total mass of $u,v$ and original energy (from left to right) for Example \ref{exam:eg3:pe}. }
	\label{fig:eg2:mass1}
\end{figure}

In the simulation, the computational domain is discretized with $M=80$ grid points in each spatial direction, and the time stepsize is set to $\tau=2.0\times10^{-3}$. As illustrated in Fig. \ref{fig:eg2:mass1}, the numerical observations for the parabolic-elliptic case are analogous to those for the parabolic-parabolic case in Example \ref{exam:eg2:pp}.
%%%%%%%%%%%%%%%%%%%%%%%%%%%%%%%%%%%%%%%%%%%%%%%%%%%%%%%%%%%%%

%%%%%%%%%%%%%%%%%%%%%%%%%%%%%%%%%%%%%%%%%%%%%%%%%%%%%%%%%%%
\subsection{Blow up of two-species}
    In this subsection, we present numerical simulations to examine the structure-preserving properties and blow-up phenomena of the proposed scheme for both the parabolic-parabolic and parabolic-elliptic cases of the two-species Keller-Segel system. These simulations are inspired by the numerical examples for the single-species Keller--Segel system in Refs. \cite{AGR'23,XF'25}.

  \begin{example}[Blow up in the parabolic-parabolic case]\label{exam:eg4:pp}
	   In this example, the computational space-time domain is chosen as $\Omega=(-0.5,0.5)^2$ and $T=1.8\times10^{-2}$, and we modify the initial conditions in Example \ref{exam:eg2:pp} as follows: 
			\begin{equation*}
				\begin{aligned}
					u^0(\bm x)&=1000 \exp\left( -100 ((x-0.25)^2+(y-0.25)^2 )\right), \\
					v^0(\bm x)&=500 \exp\left( -100 ((x-0.25)^2+(y-0.25)^2)\right), \\
					c^0(\bm x)&=0.
				\end{aligned}
			\end{equation*}	
		The corresponding total initial mass is $M[u^0]+M[v^0]\approx 31.40+ 15.70=41.7>8\pi$, and the solutions exhibit finite-time blow-up behavior.
		\end{example}

\begin{figure}[htbp]
			\centering
			\includegraphics[width=0.325\linewidth]{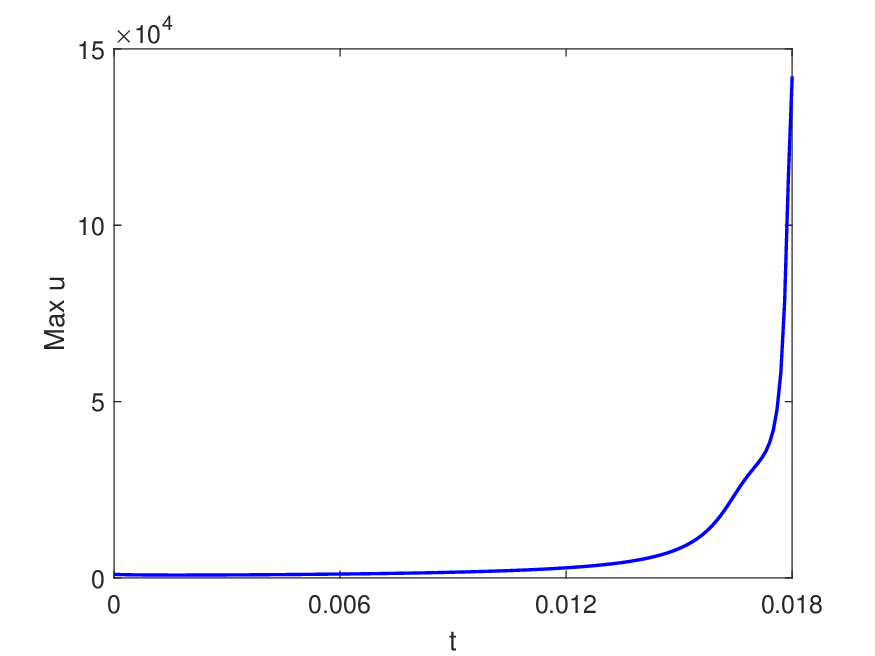}
			\includegraphics[width=0.325\linewidth]{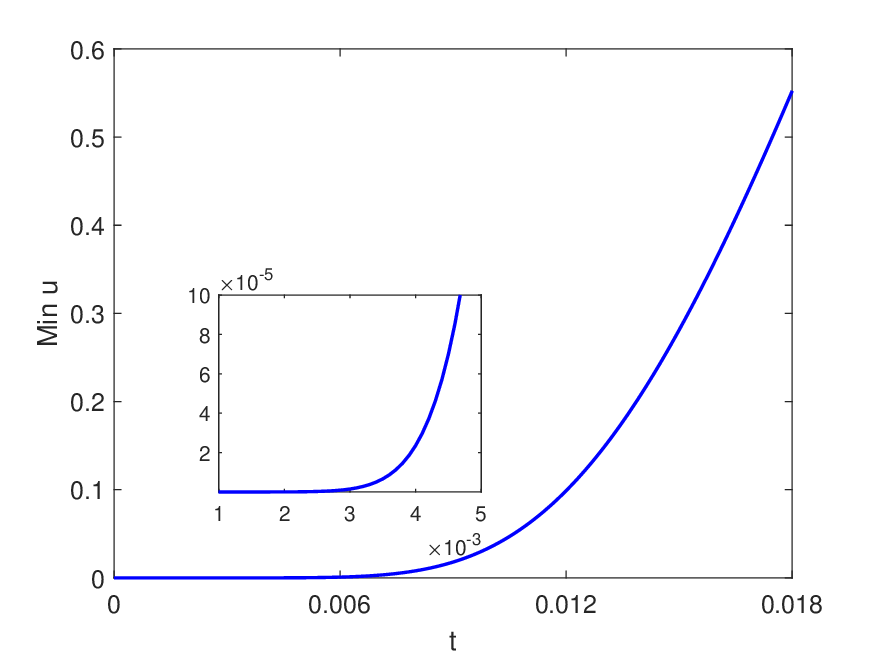}
			\includegraphics[width=0.325\linewidth]{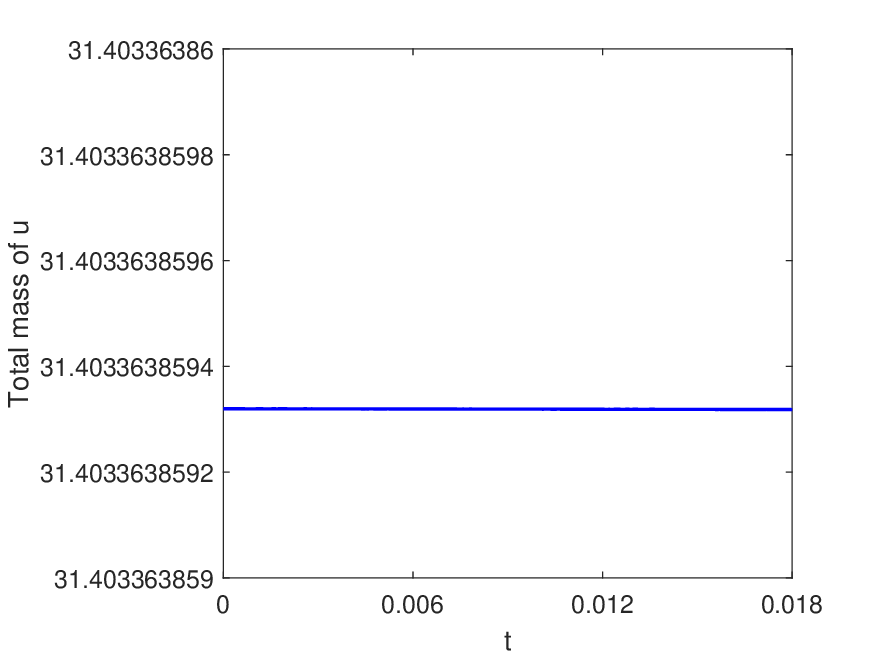}
			\includegraphics[width=0.325\linewidth]{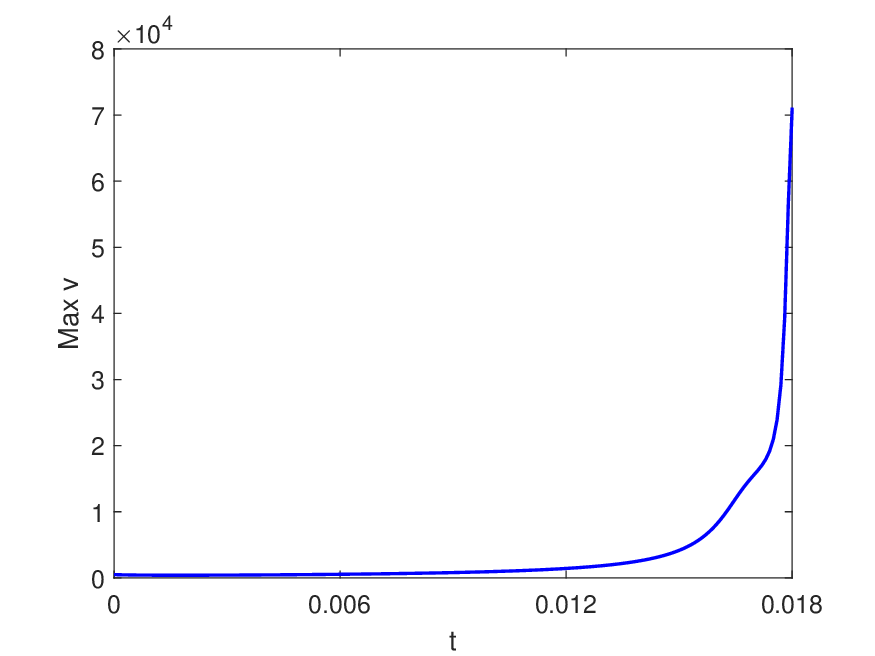}
			\includegraphics[width=0.325\linewidth]{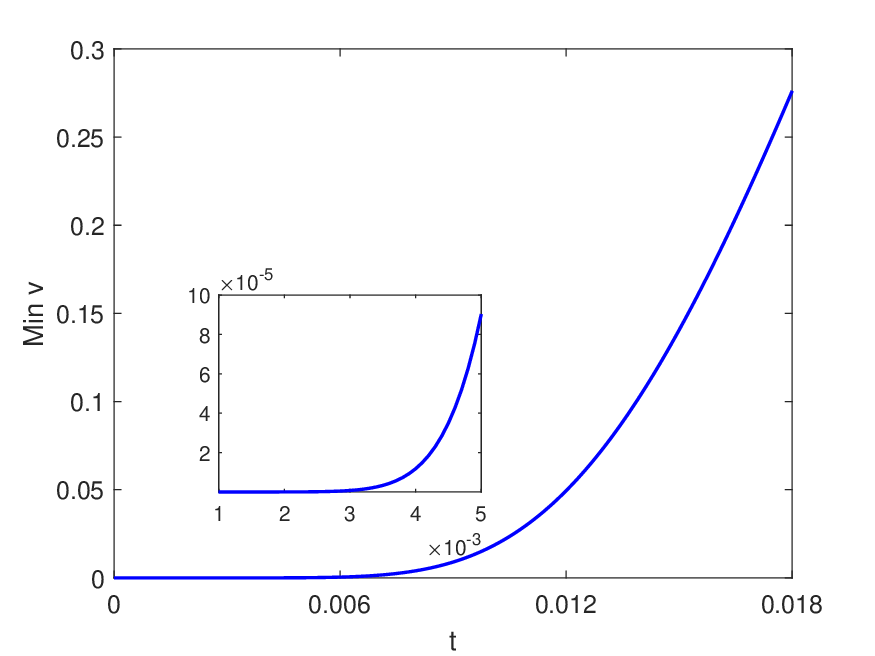}
			\includegraphics[width=0.325\linewidth]{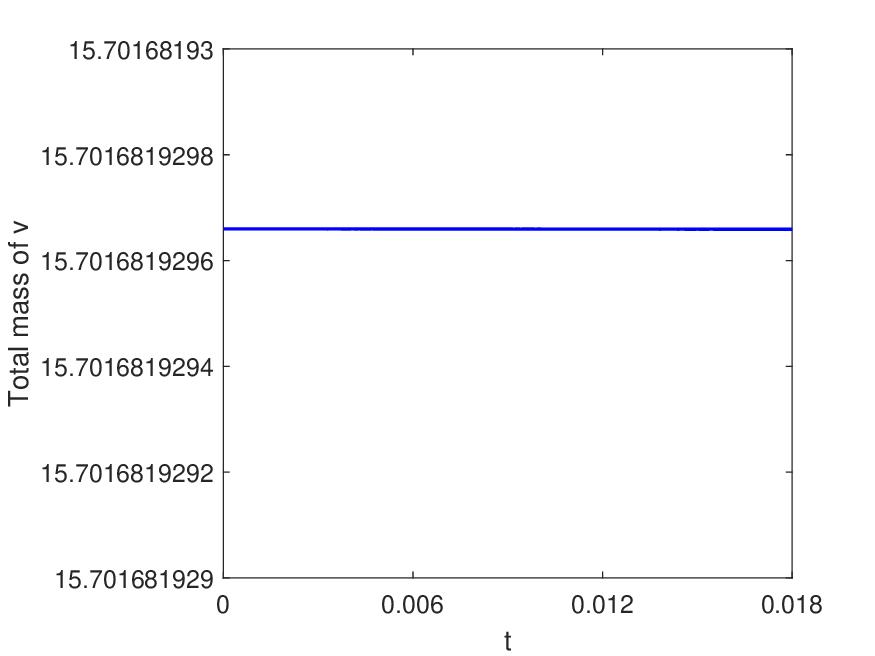}
			\includegraphics[width=0.325\linewidth]{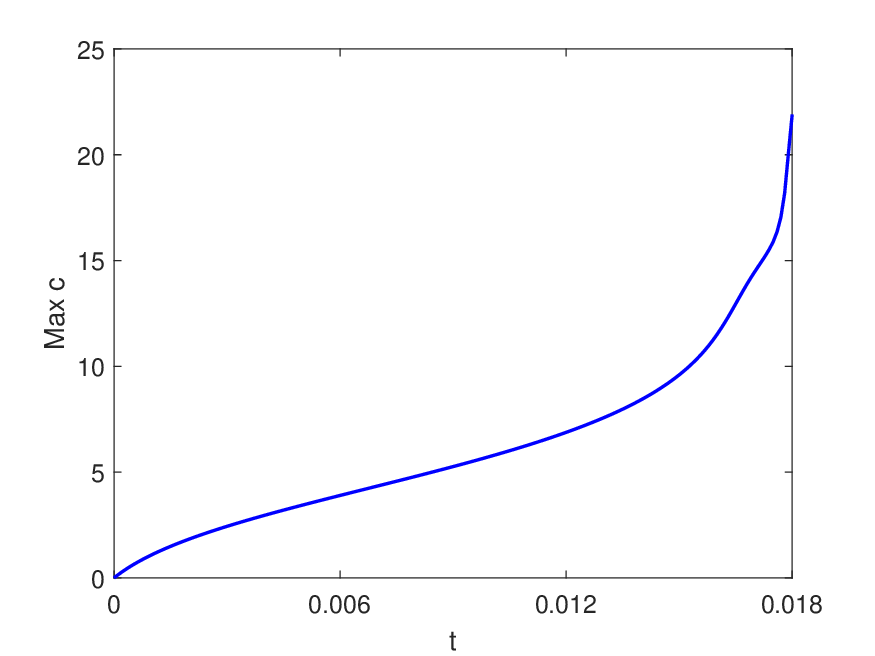}
			\includegraphics[width=0.325\linewidth]{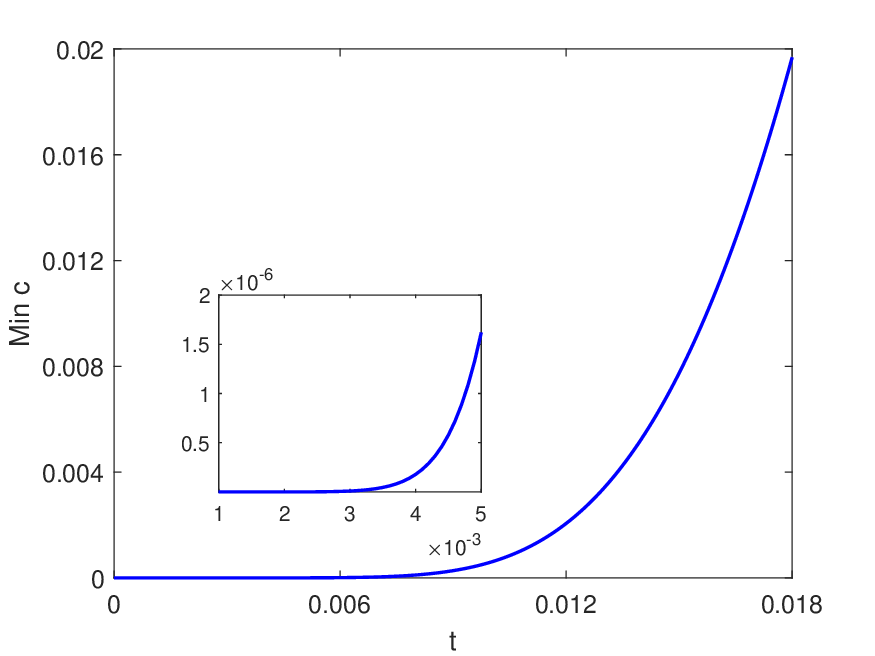}
			\includegraphics[width=0.325\linewidth]{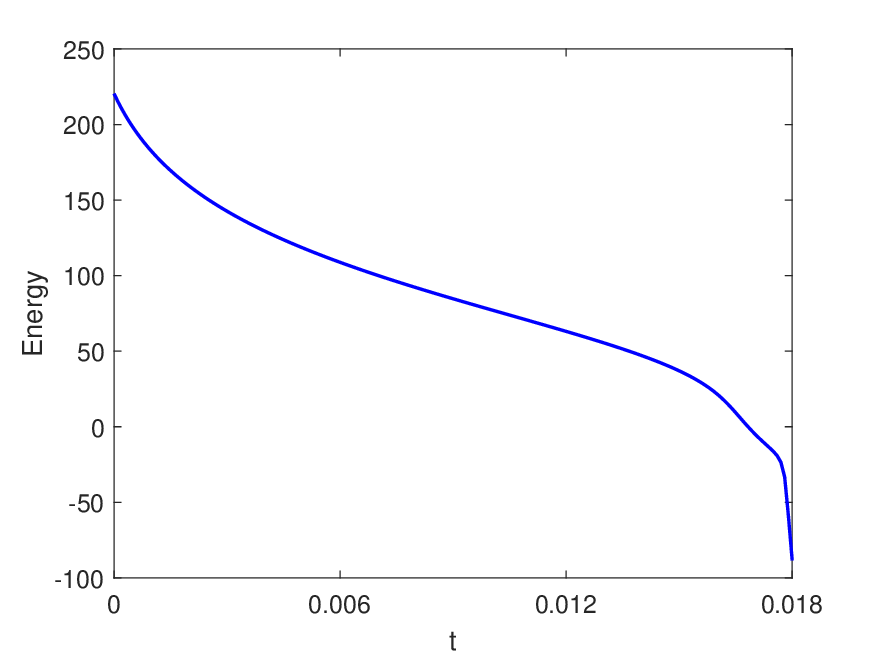}
			\caption{Evolutions of the extrema of $u$, $v$, $c$, and the total mass of $u,v$ and original energy (from left to right) for Example \ref{exam:eg4:pp}.}
			\label{fig:eg4:mass}
		\end{figure}

\begin{figure}[htbp]
			\centering
			\includegraphics[width=0.24\linewidth]{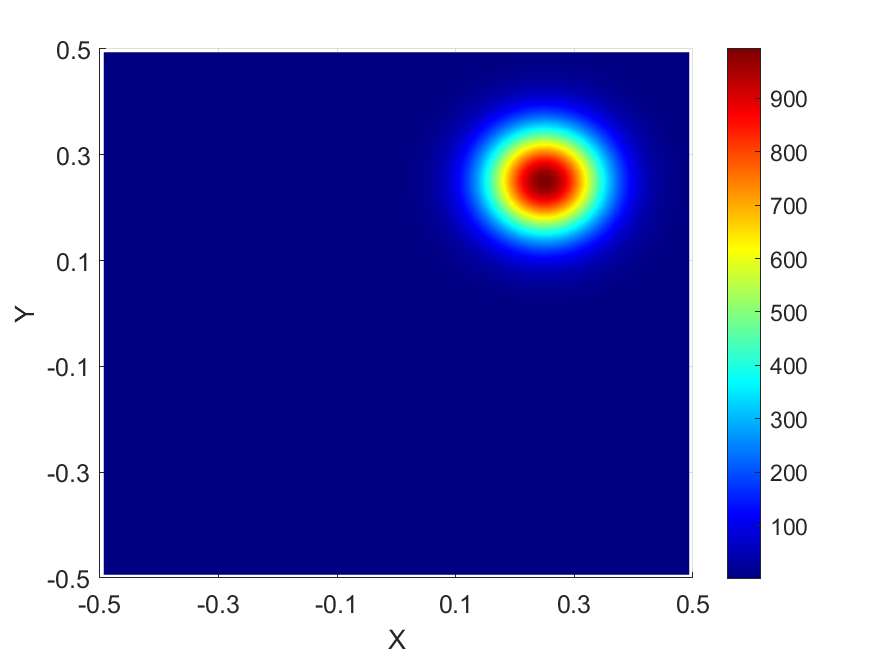}
			\includegraphics[width=0.24\linewidth]{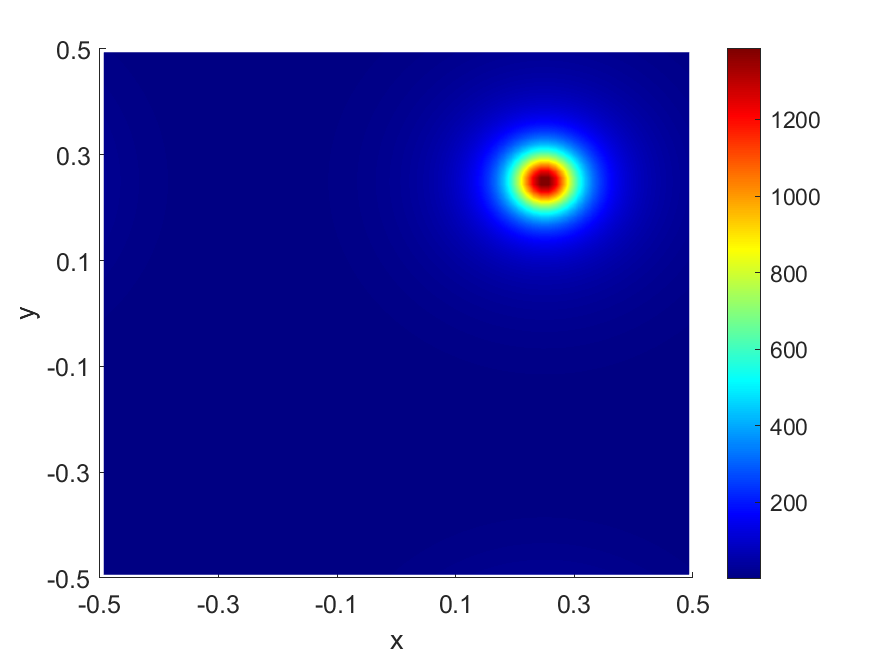}
			\includegraphics[width=0.24\linewidth]{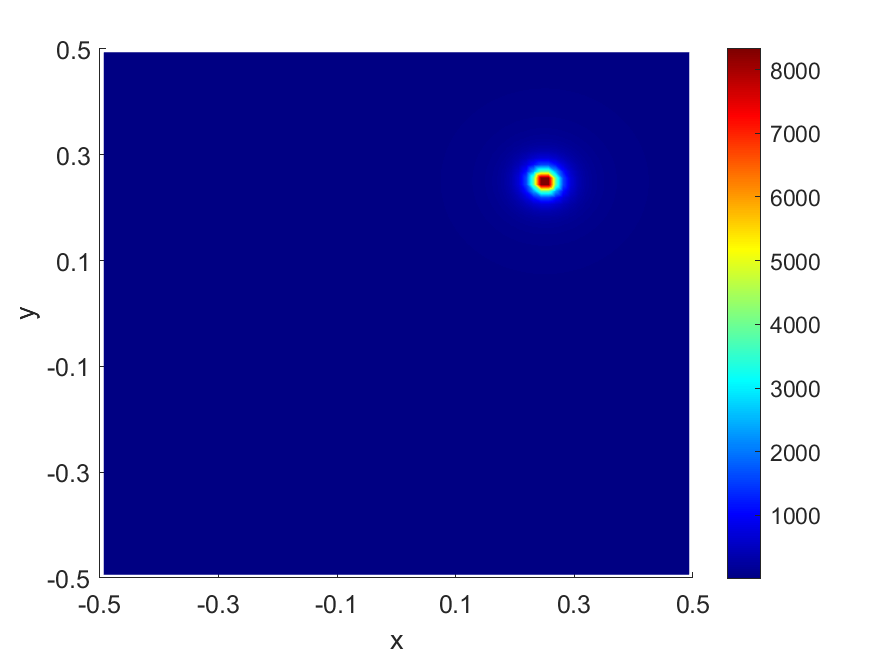}
			\includegraphics[width=0.24\linewidth]{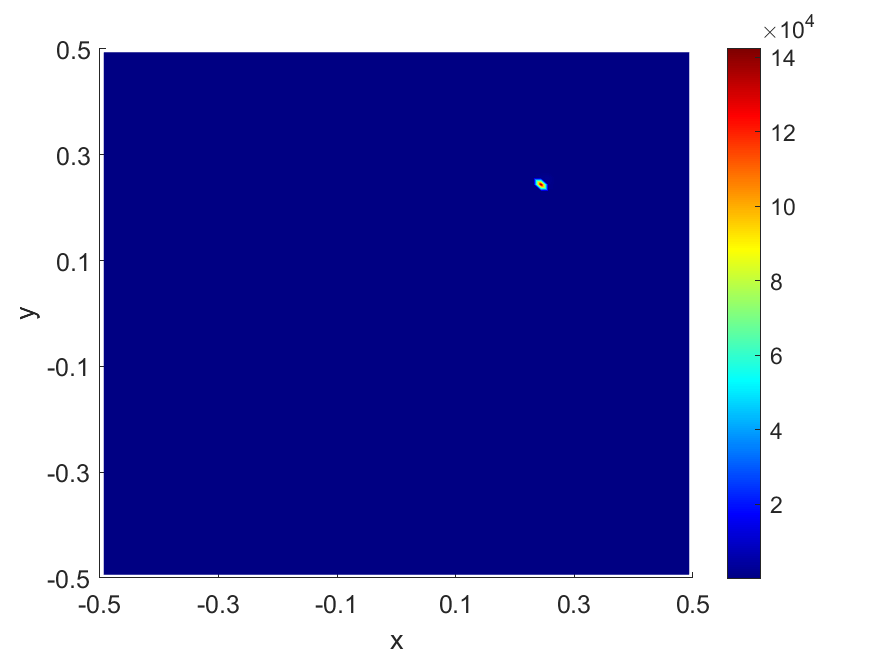}
			\includegraphics[width=0.24\linewidth]{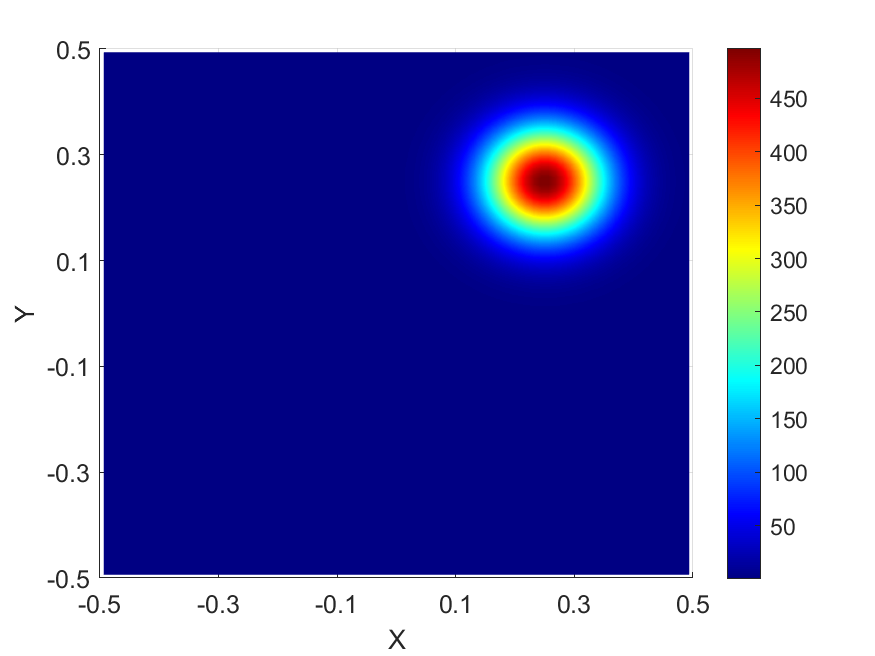}
			\includegraphics[width=0.24\linewidth]{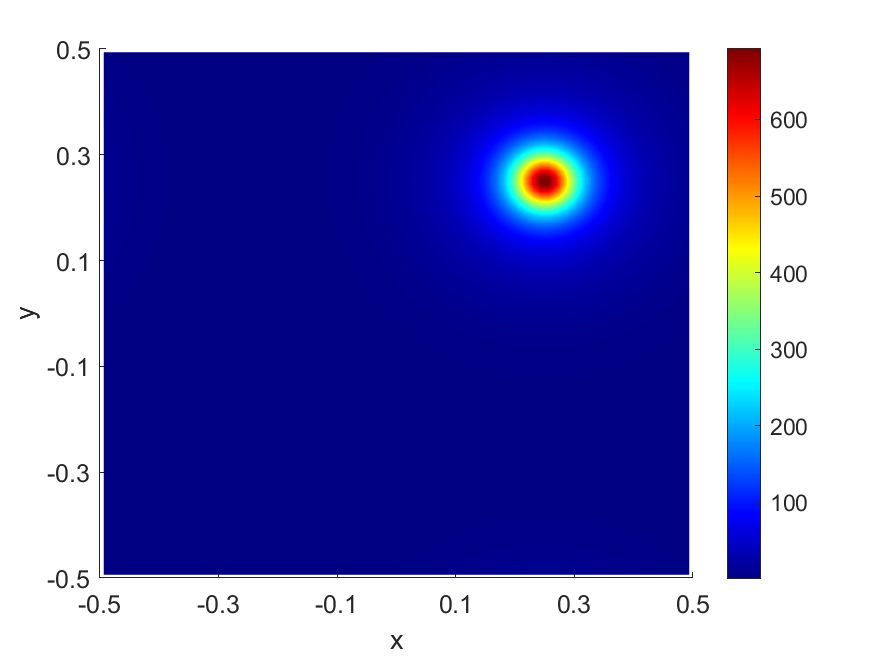}
			\includegraphics[width=0.24\linewidth]{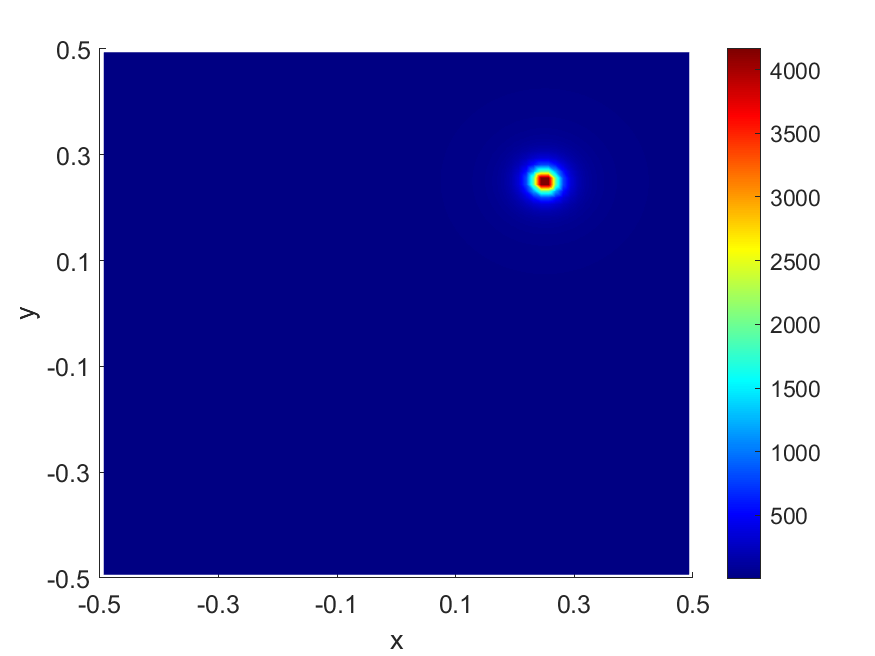}
			\includegraphics[width=0.24\linewidth]{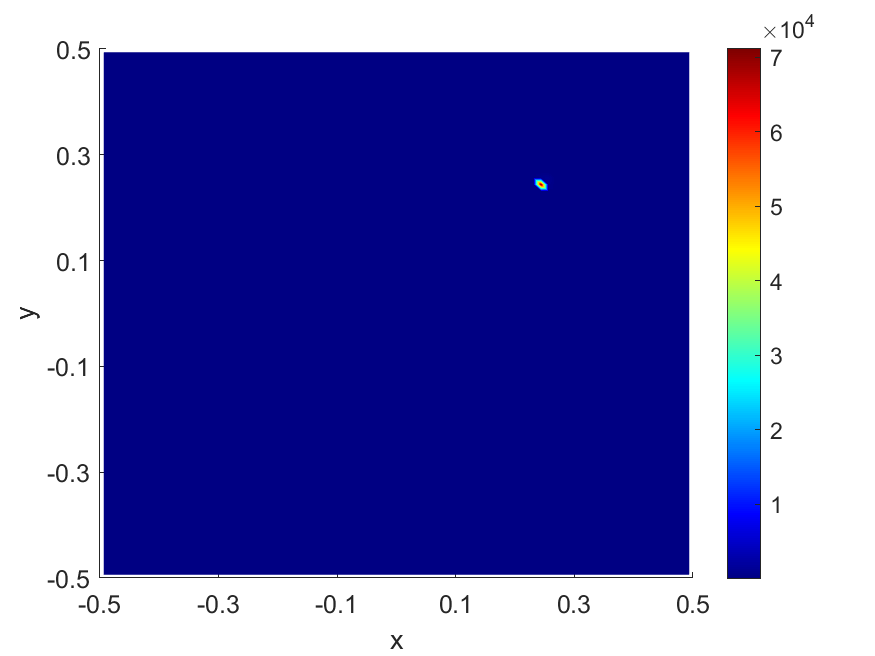}
			\includegraphics[width=0.24\linewidth]{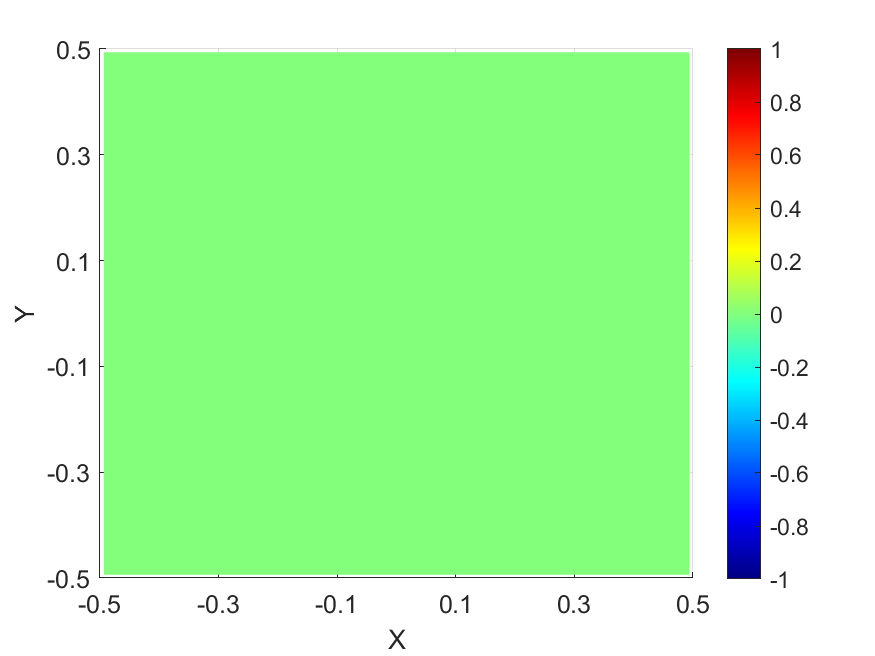}
			\includegraphics[width=0.24\linewidth]{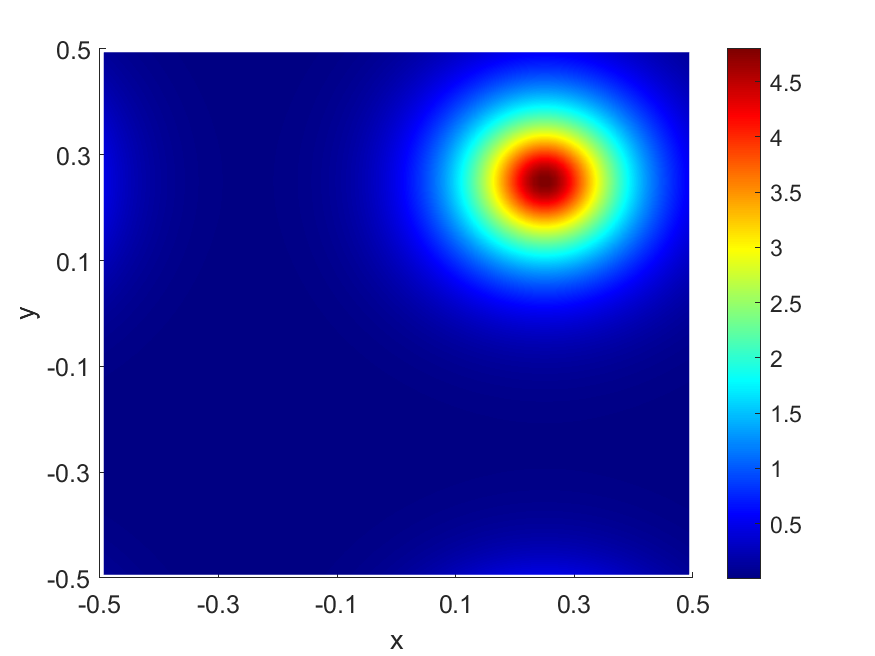}
			\includegraphics[width=0.24\linewidth]{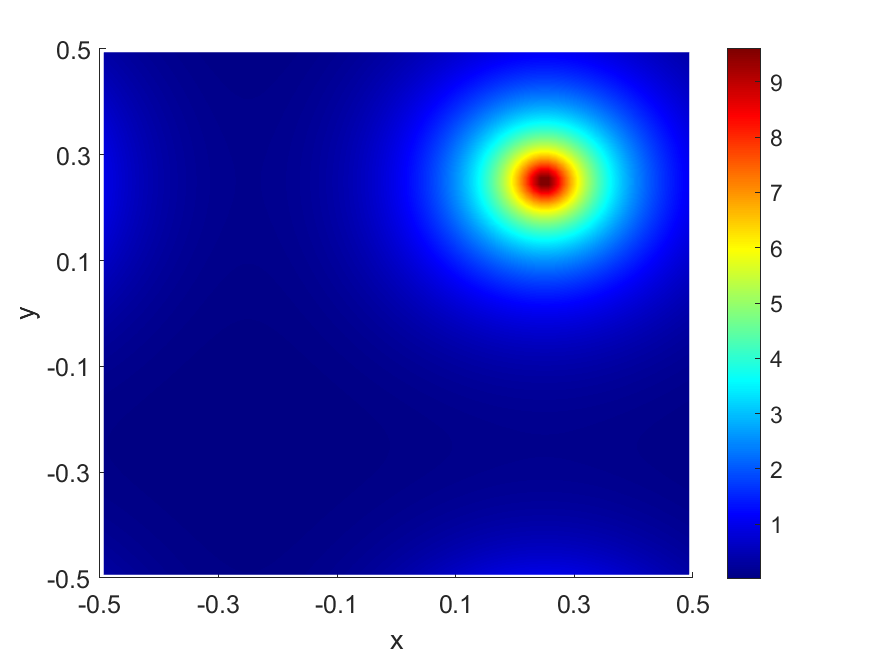}
			\includegraphics[width=0.24\linewidth]{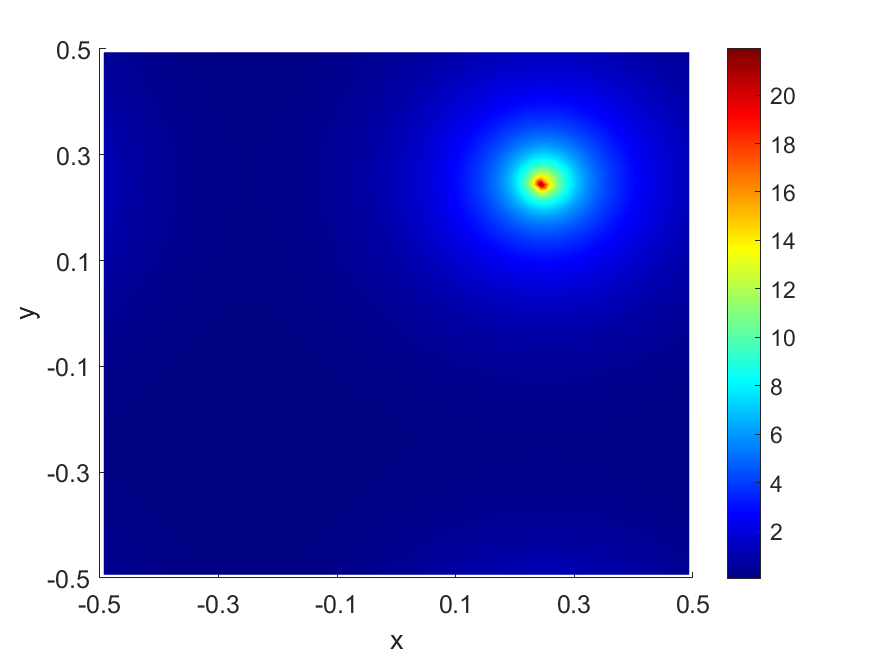}
			\caption{Contour plots of $u,v,c$ (from top to bottom) at time instants $t = 0,\ 8.0\times10^{-3},\ 1.5\times10^{-2},\ 1.8\times10^{-2}$ (from left to right)  for Example \ref{exam:eg4:pp}.}
			\label{fig:eg4:blow_up}
		\end{figure}

  In this simulation, we use $M=80$ grid points in each spatial direction, with a time stepsize $\tau=1.0\times10^{-5}$. Fig. \ref{fig:eg4:mass} illustrates the temporal evolution of the maximum and minimum values of $u$, $v$, and $c$, the total masses of the densities $u$ and $v$, and the original energy. The numerical results demonstrate that the densities $u$ and $v$ remain positive and conserve their respective masses, while the chemoattractant concentration $c$ also remains positive. Furthermore, the results in Fig. \ref{fig:eg4:mass} indicate that the original energy is dissipative. Finally, the temporal evolution of the numerical solutions $u$, $v$, and $c$ is presented in Fig. \ref{fig:eg4:blow_up}, where the cell densities $u$ and $v$ exhibit pronounced concentration and rapid growth as $t$ approaches $1.8\times 10^{-2}$, indicating the onset of finite-time blow-up behavior consistent with the expected dynamics of the Keller--Segel system.
        
	\begin{example}[Blow up in the parabolic-elliptic case]\label{exam:eg5:pe}
		In this example, we examine the finite-time blow-up behavior of the parabolic-elliptic system with two initially localized cell-density peaks on the computational domain $\Omega=(-0.5,0.5)^2$. The initial conditions are prescribed as 
		\begin{equation*}
			\begin{aligned}
				u^0(\bm x)&=600 \exp\left( -100 ((x-0.2)^2+y^2 )\right)+500\exp\left( -100 (x^2+(y-0.2)^2 )\right), \\
				v^0(\bm x)&=400 \exp\left( -100 ((x-0.2)^2+y^2 )\right)+300\exp\left( -100 (x^2+(y-0.2)^2 )\right).
			\end{aligned}
		\end{equation*}	
	For the initial distributions, the total masses of the two cell species satisfy $M[u^0]+M[v^0]\approx 34.56+ 21.99=56.55>8\pi$, and hence finite-time blow-up occurs.
	\end{example}

\begin{figure}[htbp]
			\centering
			\includegraphics[width=0.325\linewidth]{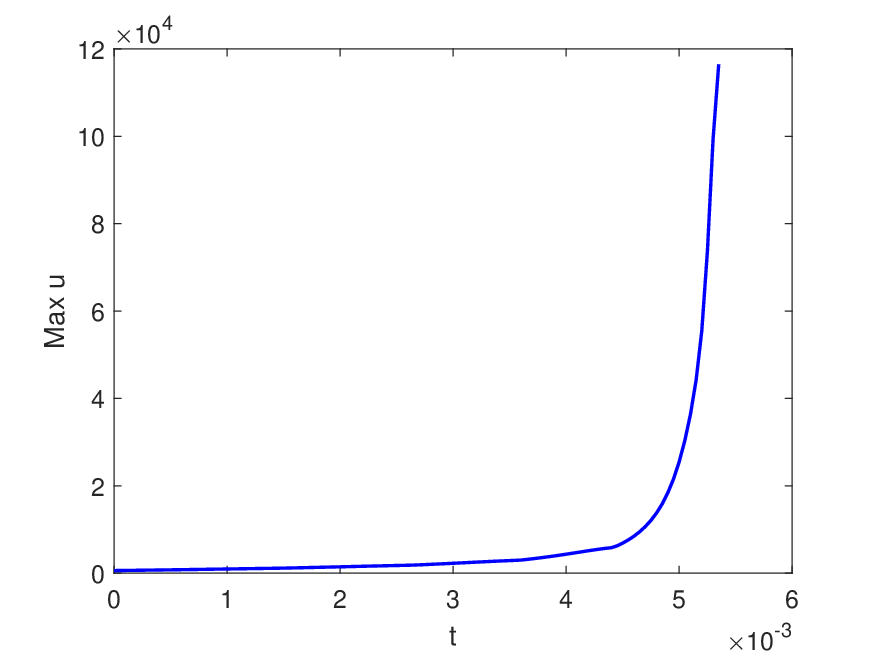}
			\includegraphics[width=0.325\linewidth]{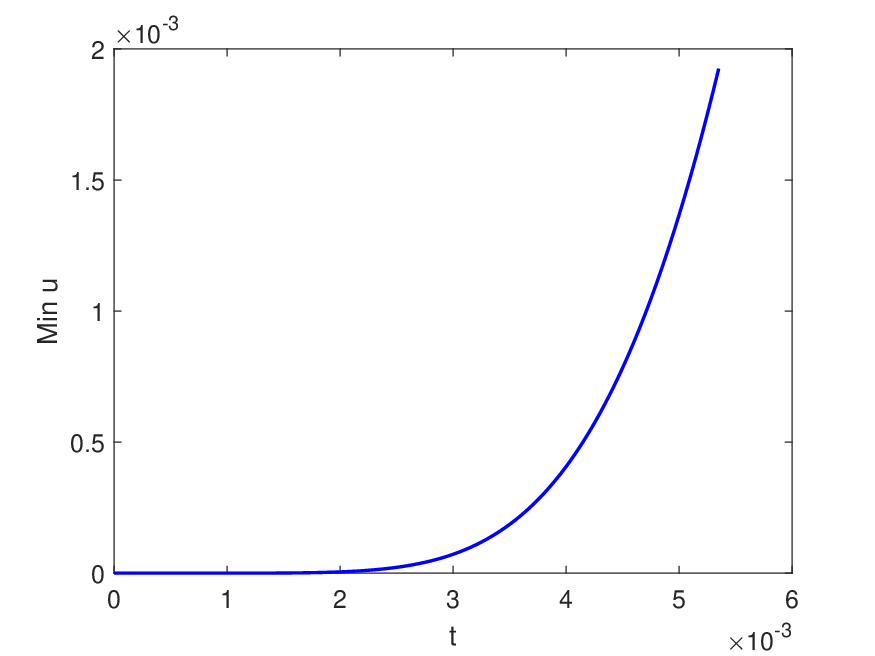}
			\includegraphics[width=0.325\linewidth]{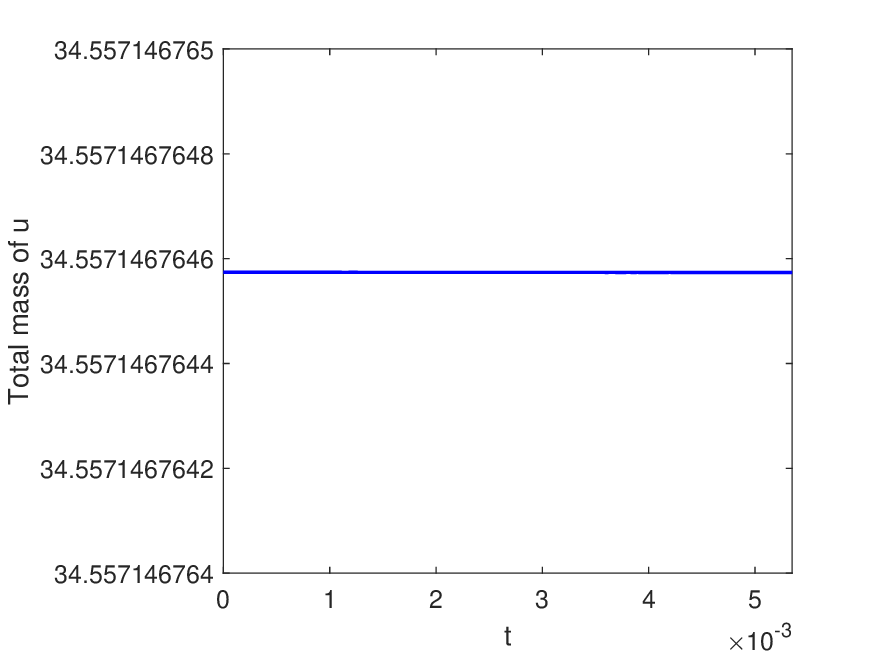}
			\includegraphics[width=0.325\linewidth]{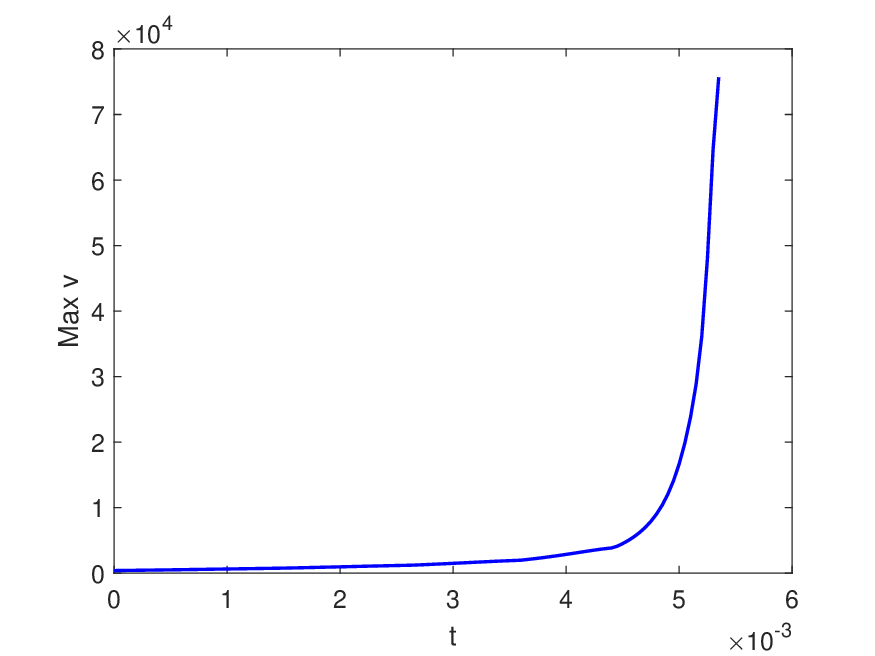}
			\includegraphics[width=0.325\linewidth]{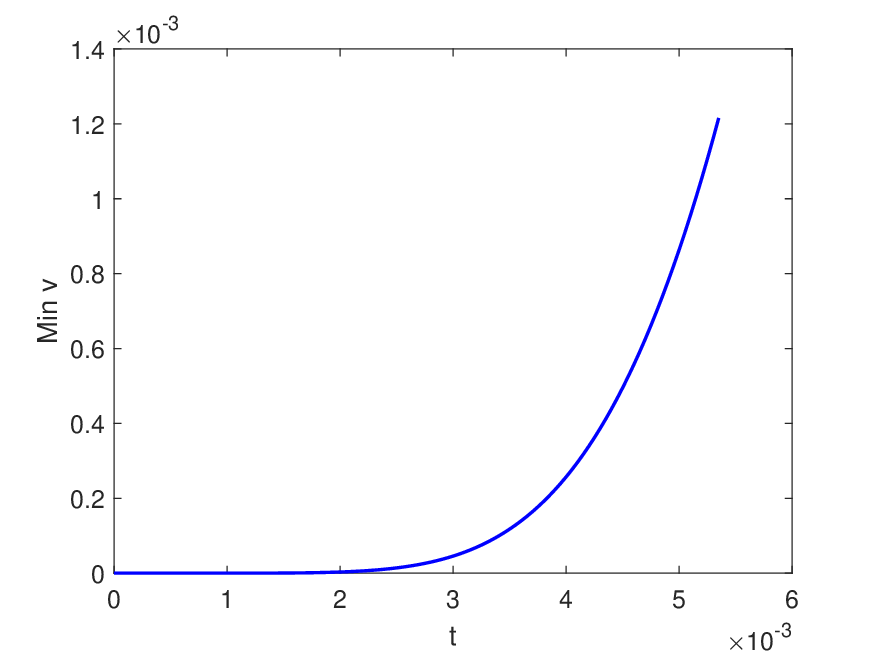}
			\includegraphics[width=0.325\linewidth]{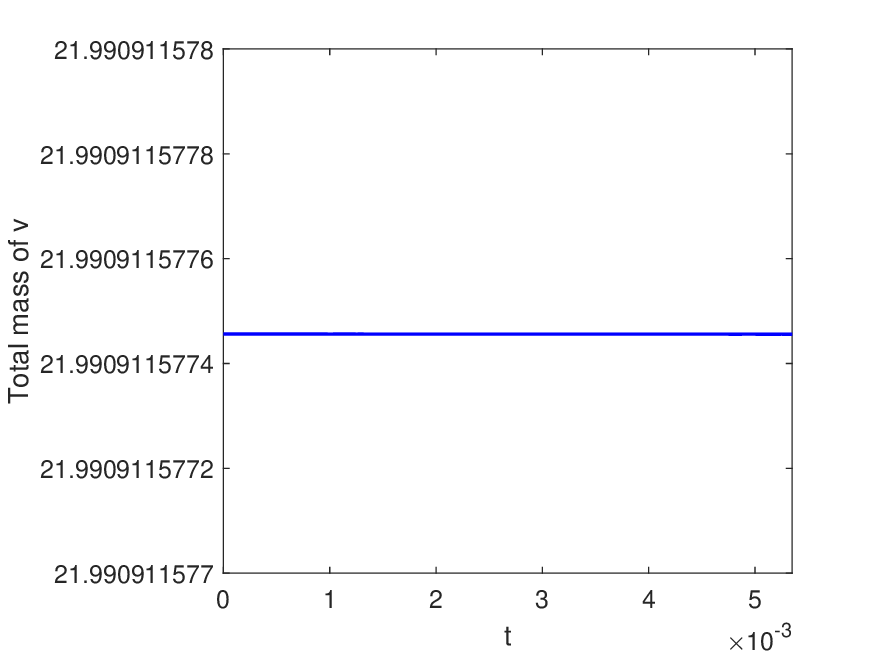}
			\includegraphics[width=0.325\linewidth]{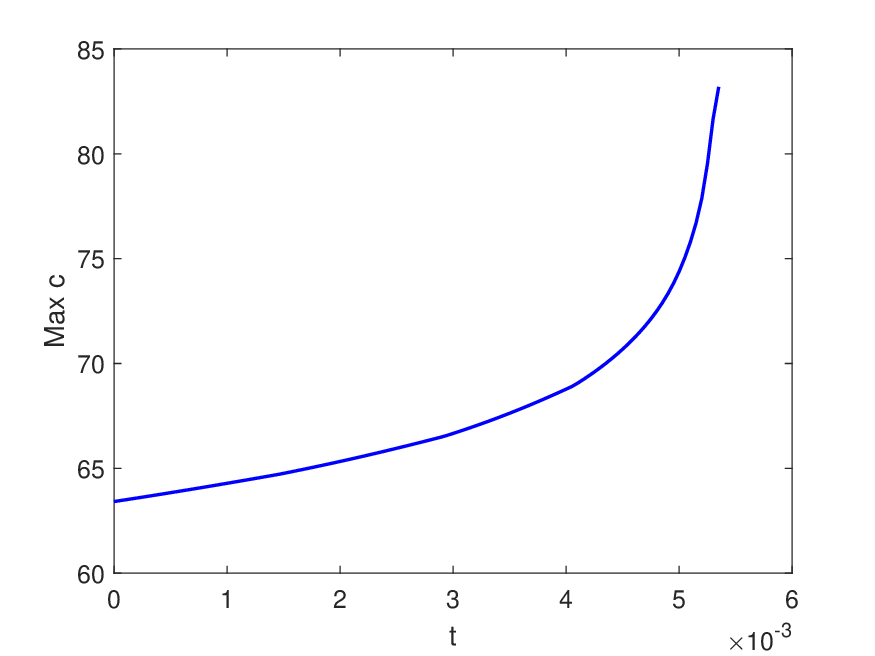}
			\includegraphics[width=0.325\linewidth]{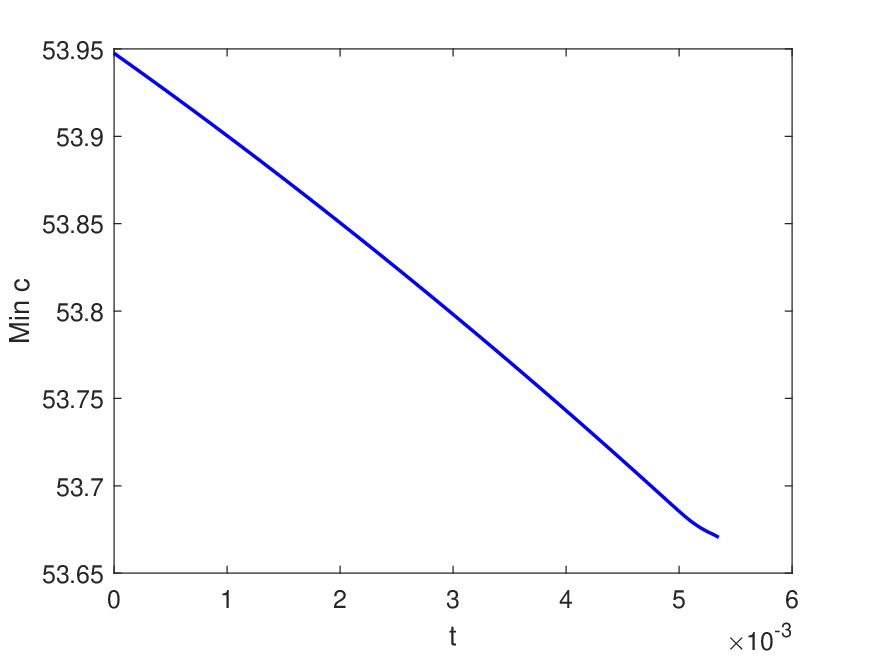}
			\includegraphics[width=0.325\linewidth]{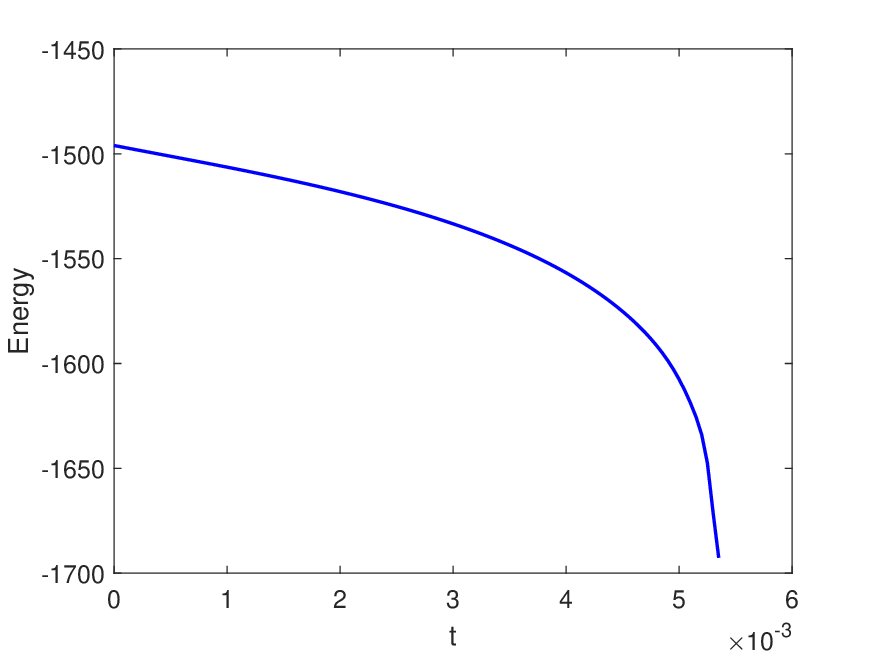}
			\caption{Evolutions of the extrema of $u$, $v$, $c$, and the total mass of $u,v$ and original energy (from left to right) for Example \ref{exam:eg5:pe}.}
			\label{fig:eg5:mass}
		\end{figure}

\begin{figure}[htbp]
			\centering
			\includegraphics[width=0.24\linewidth]{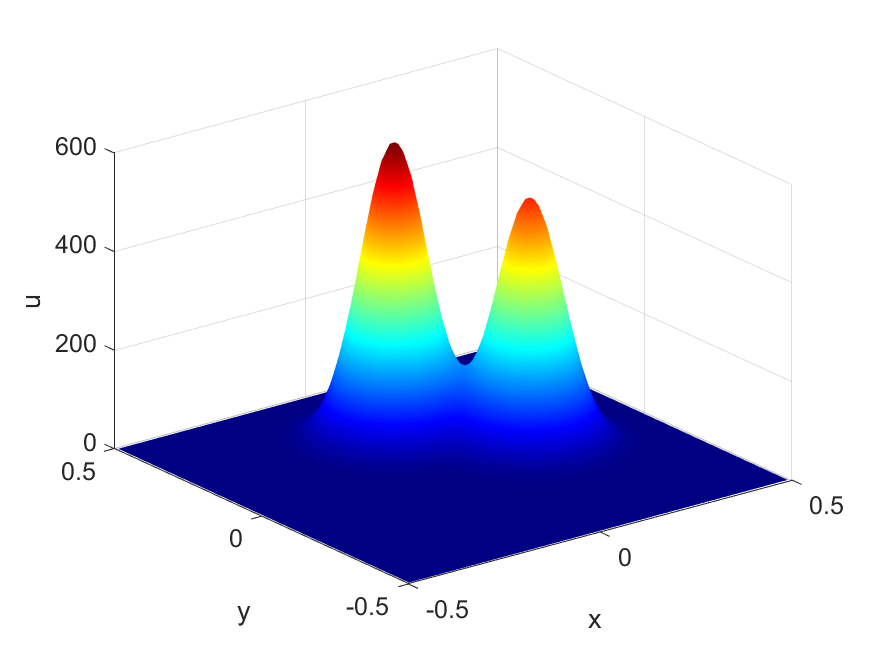}
			\includegraphics[width=0.24\linewidth]{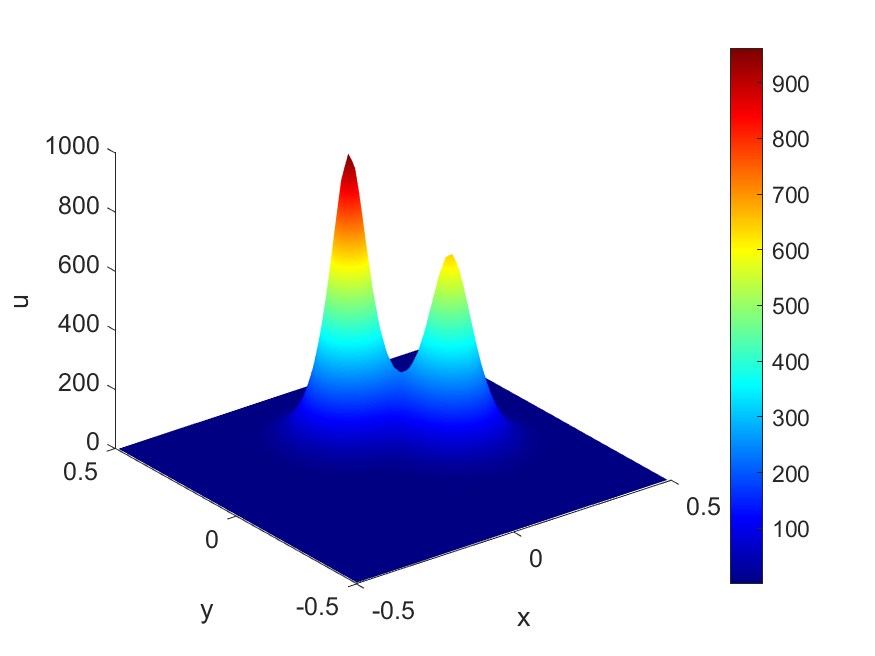}
			\includegraphics[width=0.24\linewidth]{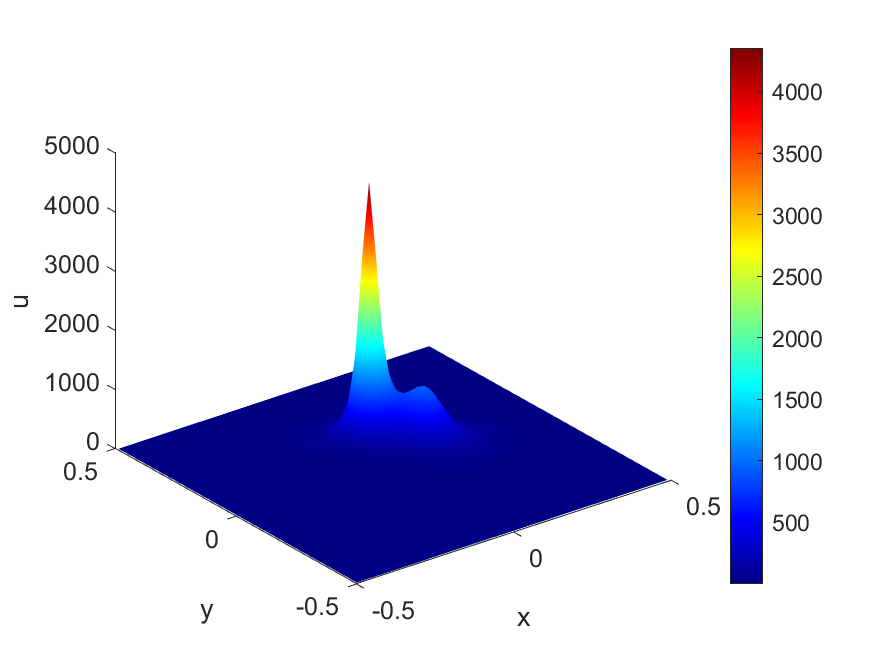}
			\includegraphics[width=0.24\linewidth]{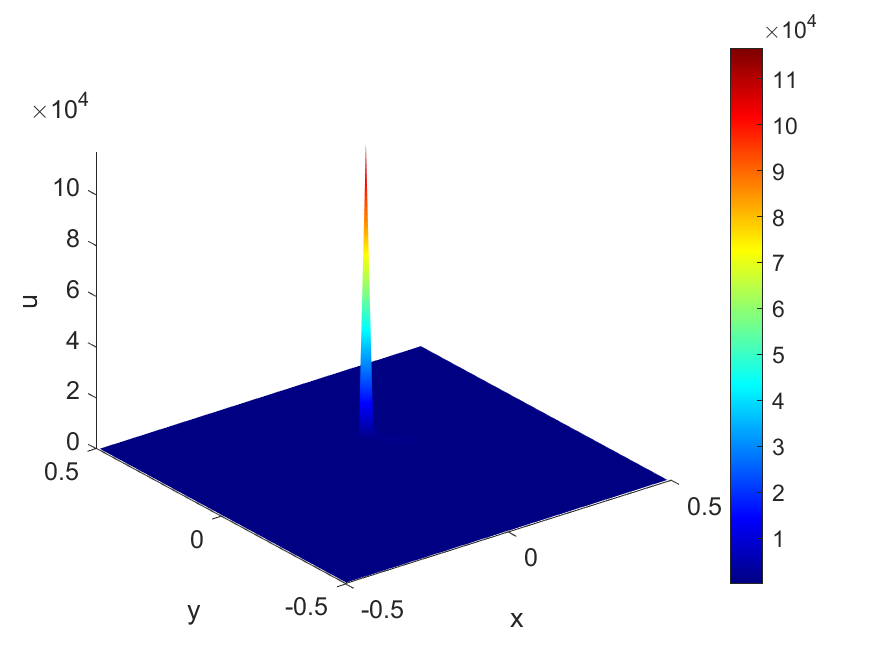}
			\includegraphics[width=0.24\linewidth]{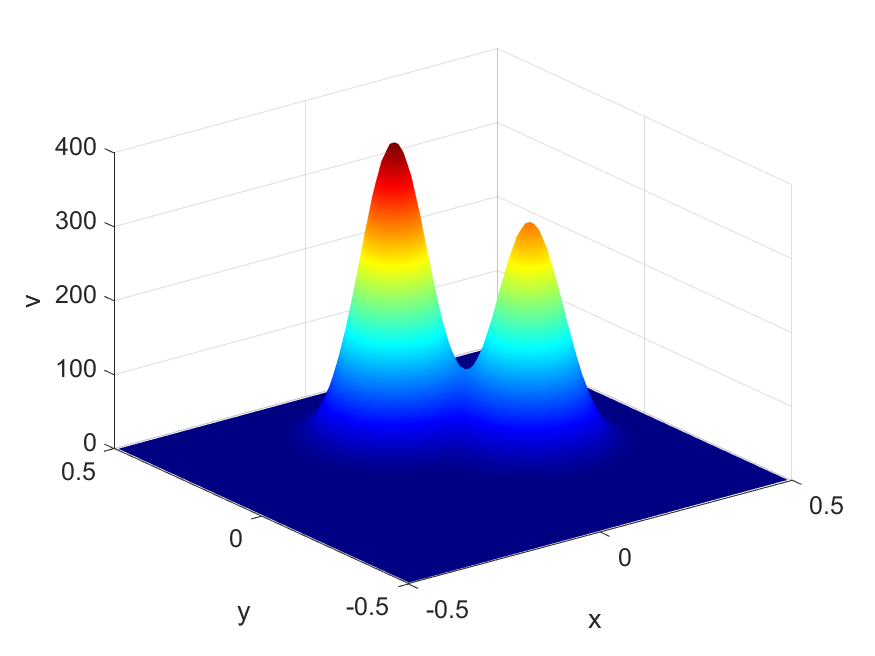}
			\includegraphics[width=0.24\linewidth]{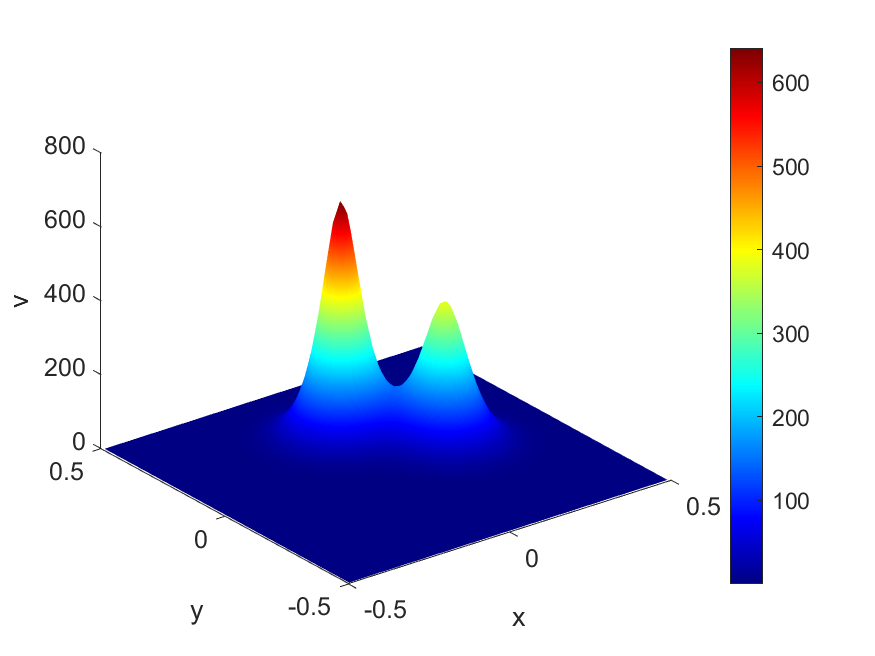}
			\includegraphics[width=0.24\linewidth]{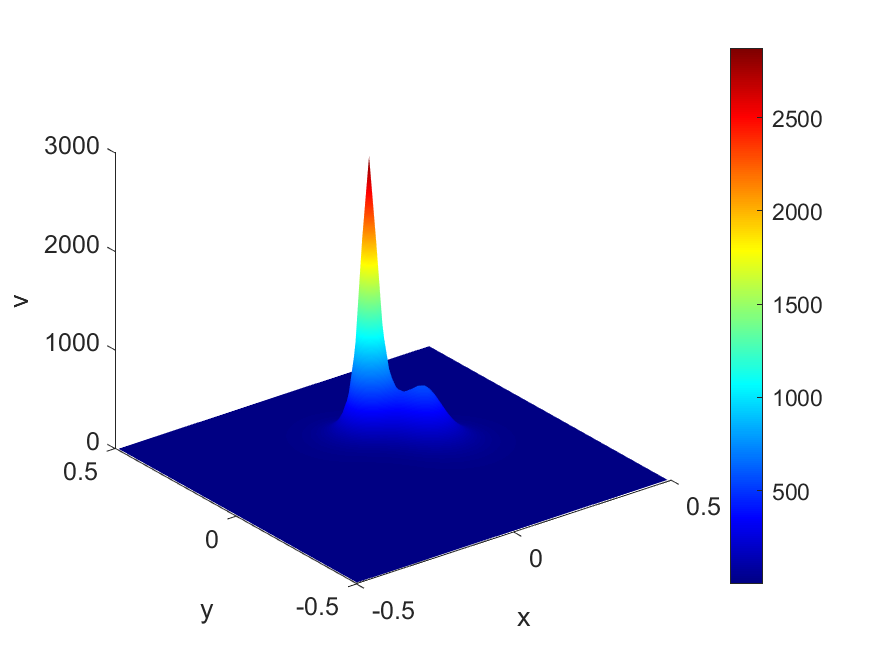}
			\includegraphics[width=0.24\linewidth]{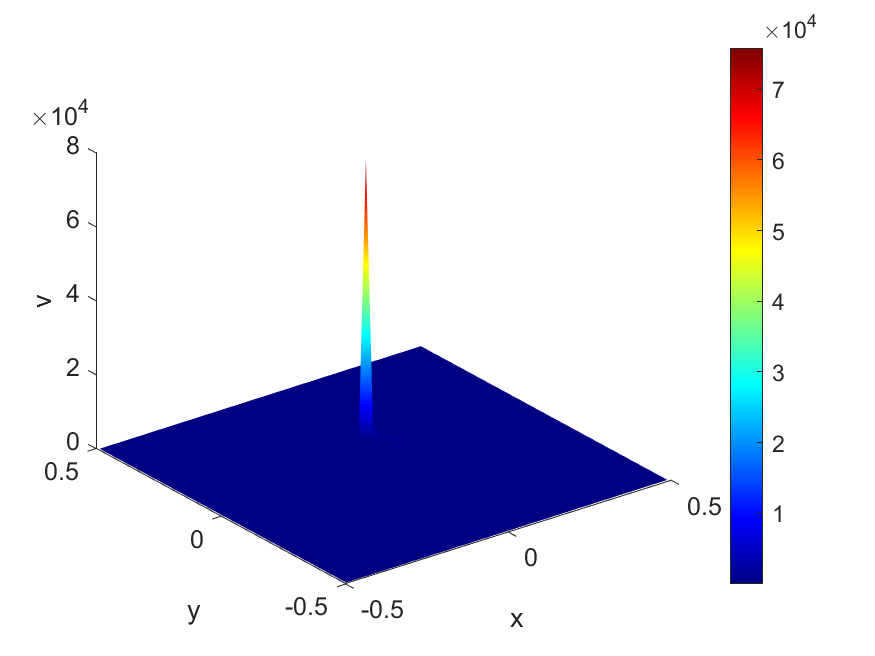}
			\includegraphics[width=0.24\linewidth]{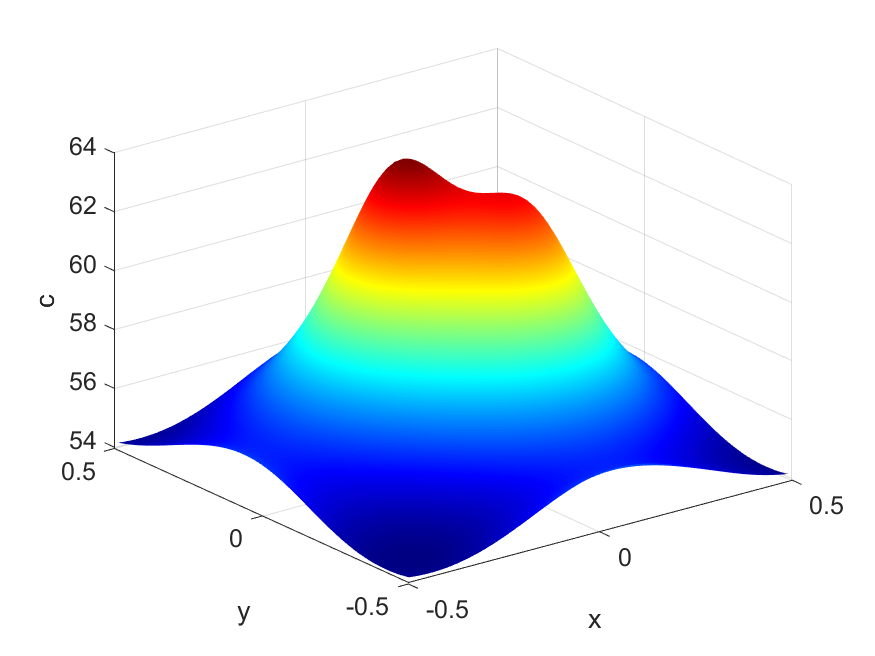}
			\includegraphics[width=0.24\linewidth]{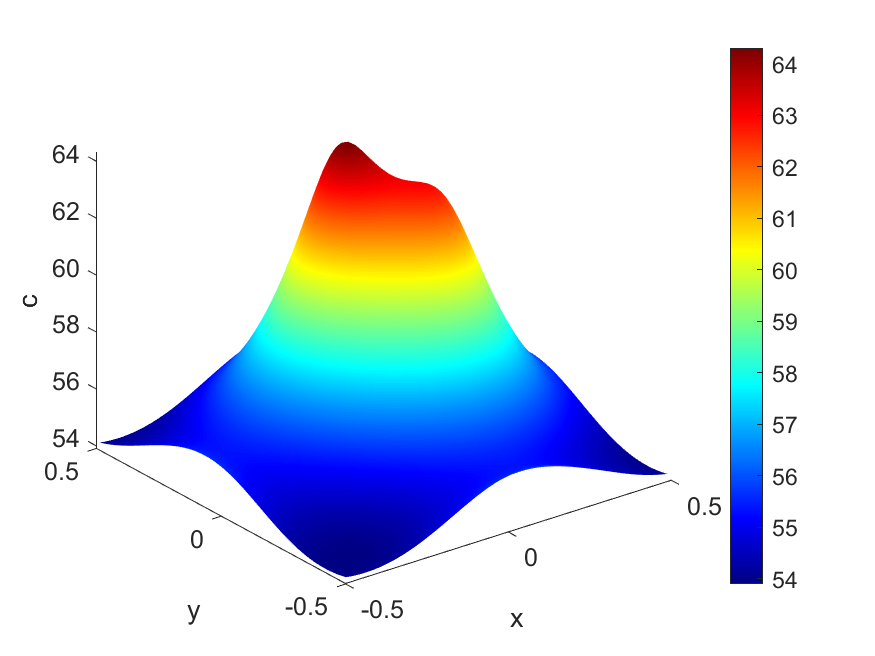}
			\includegraphics[width=0.24\linewidth]{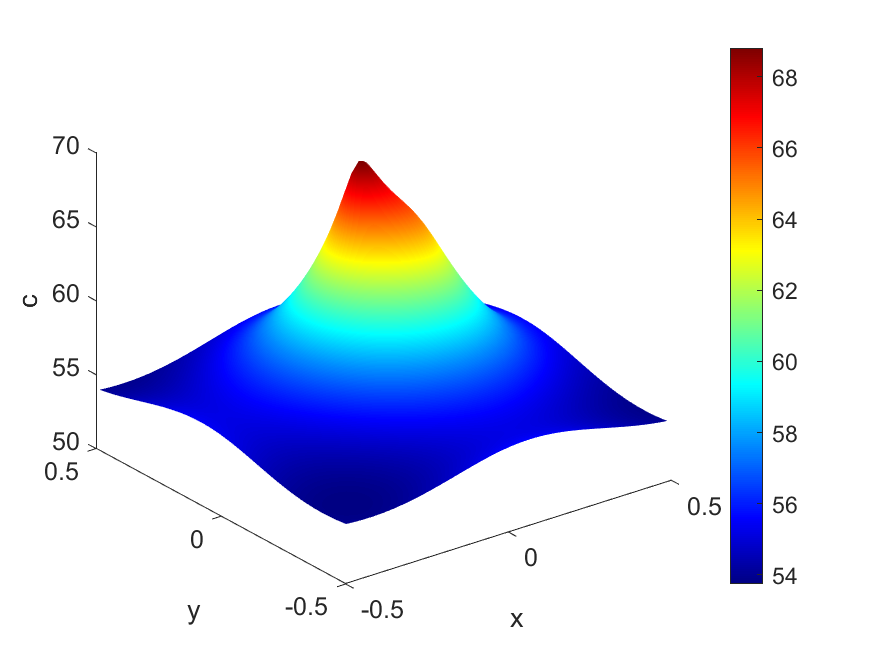}
			\includegraphics[width=0.24\linewidth]{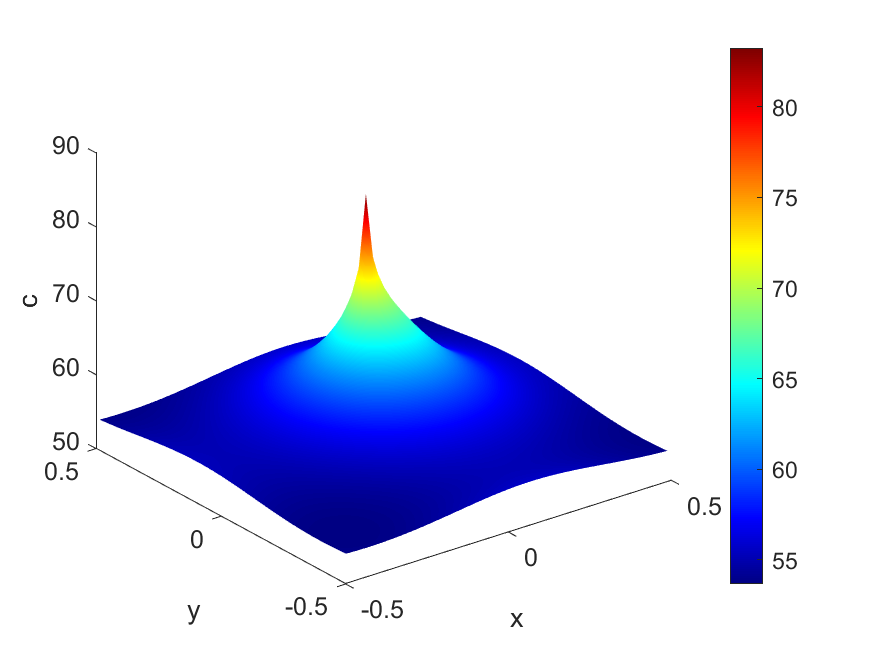}
			\caption{Snapshots of $u,v,c$ (from top to bottom) at time instants $t = 0,\ 1.0\times10^{-3},\ 4.0\times10^{-3},\ 5.35\times10^{-3}$ (from left to right)  for Example \ref{exam:eg5:pe}.}
			\label{fig:eg5:blow_up}
		\end{figure}

   To further assess the structure-preserving performance of the proposed scheme  prior to the onset of blow-up, numerical simulations are performed for the two-species Keller--Segel system with two initially localized cell-density peaks. The spatial discretization employs $M=80$ grids points in each spatial direction, with a time stepsize $\tau = 5.0 \times 10^{-5}$. As shown in Fig. \ref{fig:eg5:mass}, the numerical solutions retain the structure-preserving properties observed in Example \ref{exam:eg4:pp}. More notably, Fig. \ref{fig:eg5:blow_up} depicts the evolution of the two cell densities and the chemoattractant concentration, revealing a gradual merging of the two initially separated peaks into a highly concentrated region. This aggregation becomes increasingly pronounced as $T=5.35\times10^{-3}$, suggesting the emergence of finite-time blow-up behavior consistent with the expected aggregation dynamics of the Keller--Segel system.

%%%%%%%%%%%%%%%%%%%%%%%%%%%%%%%%%%%%%%%%%%%%%%%%%%%%%%%%%%%%%%%
%\subsection{3D simulation of blow up behavior}
	\begin{example}[Blow up in the 3D parabolic-parabolic case]\label{exam:bu:3D} As a final numerical experiment, we extent our simulations to the 3D two-species Keller--Segel system  on the cubic domain $\Omega=(0,1)^3$. To characterize the aggregation dynamics in 3D, we choose the following spatially localized initial distributions:
			\begin{equation*}
				\begin{aligned}
					u^0(\bm x)&=2000 \exp\left( -100 ((x-0.5)^2+(y-0.5)^2+(z-0.5)^2 )\right) , \\
					v^0(\bm x)&=1000 \exp\left( -100 ((x-0.5)^2+(y-0.5)^2+(z-0.5)^2 )\right) , \\
					c^0(\bm x)&= \exp\left( -50 ((x-0.5)^2+(y-0.5)^2+(z-0.5)^2 )\right) .
				\end{aligned}
			\end{equation*}
		These initial profiles are concentrated around the center $(0.5, 0.5, 0.5)$ of the domain, and the masses of the two cell species are $M[u^0] \approx 11.14$ and $M[v^0] \approx 5.57$, respectively.
	\end{example}

\begin{figure}[htbp]
			\centering
			\includegraphics[width=0.325\linewidth]{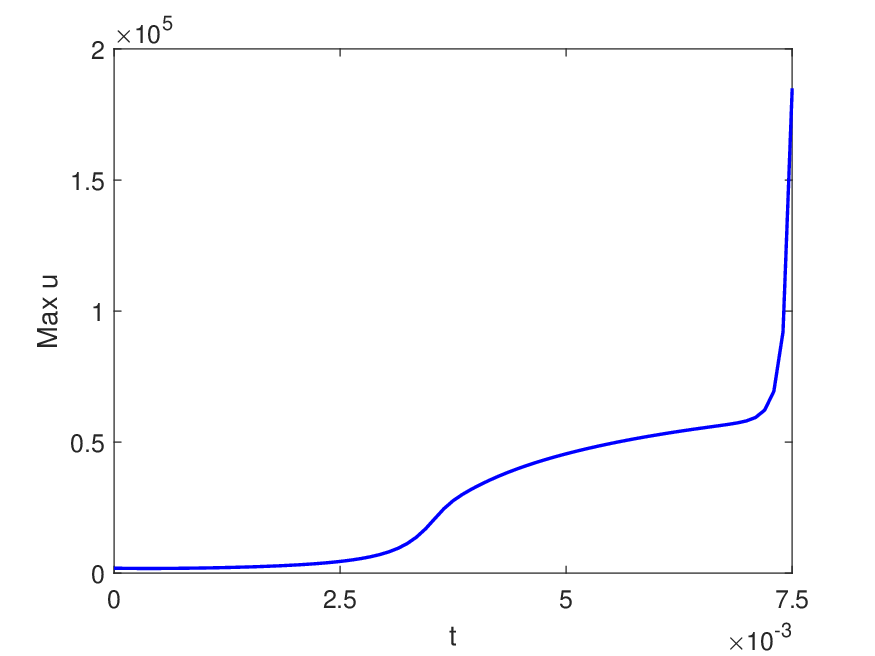}
			\includegraphics[width=0.325\linewidth]{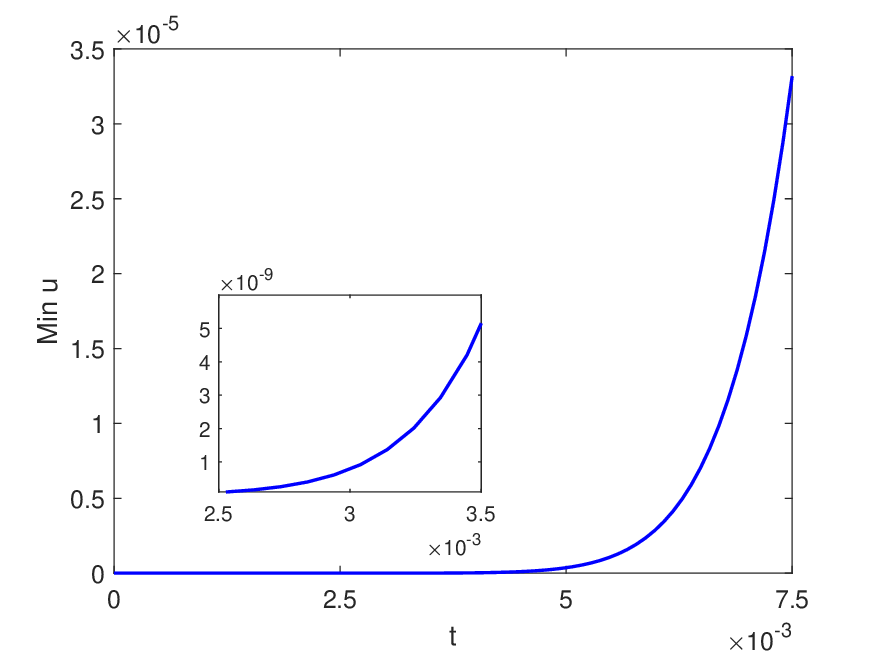}
			\includegraphics[width=0.325\linewidth]{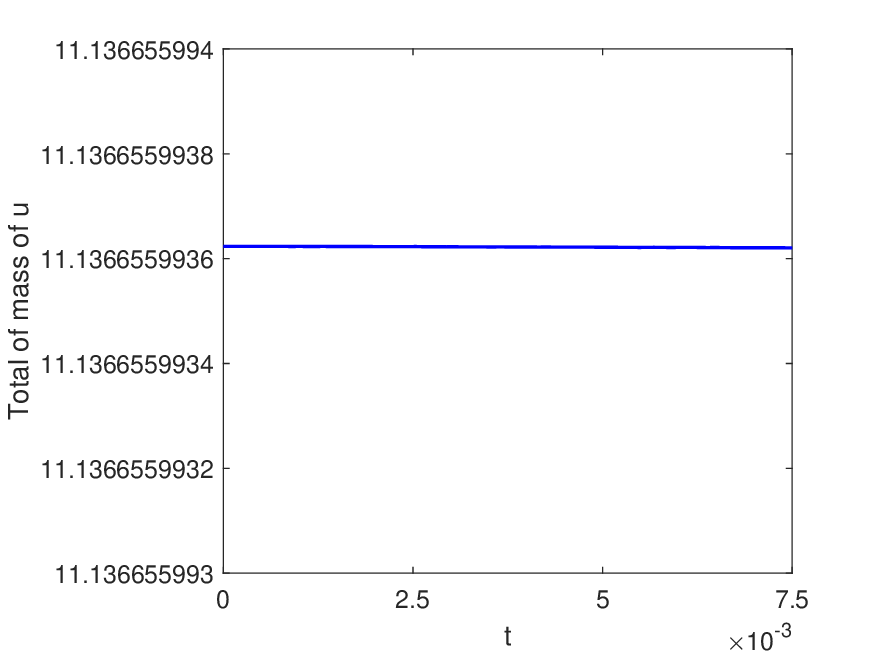}
			\includegraphics[width=0.325\linewidth]{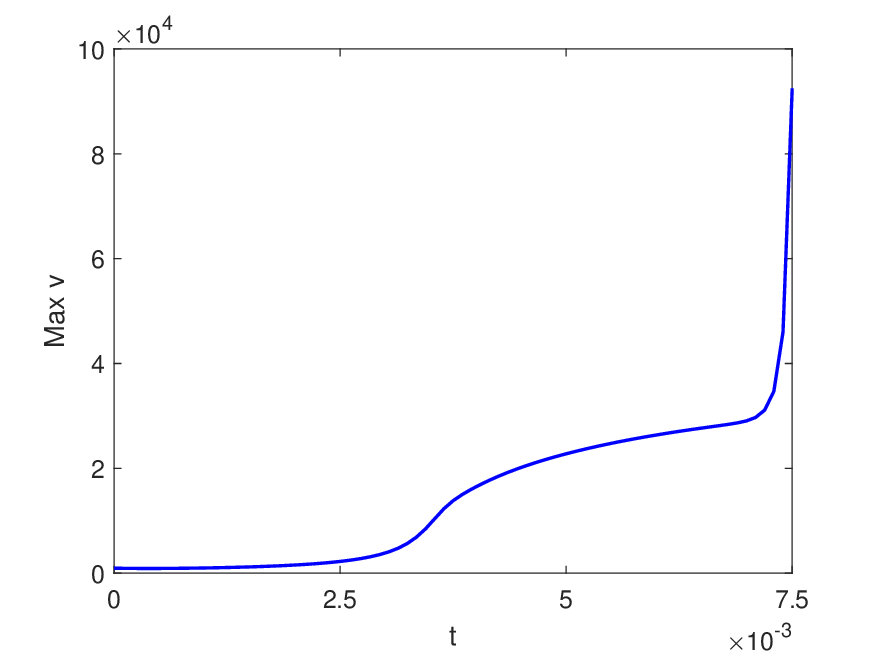}
			\includegraphics[width=0.325\linewidth]{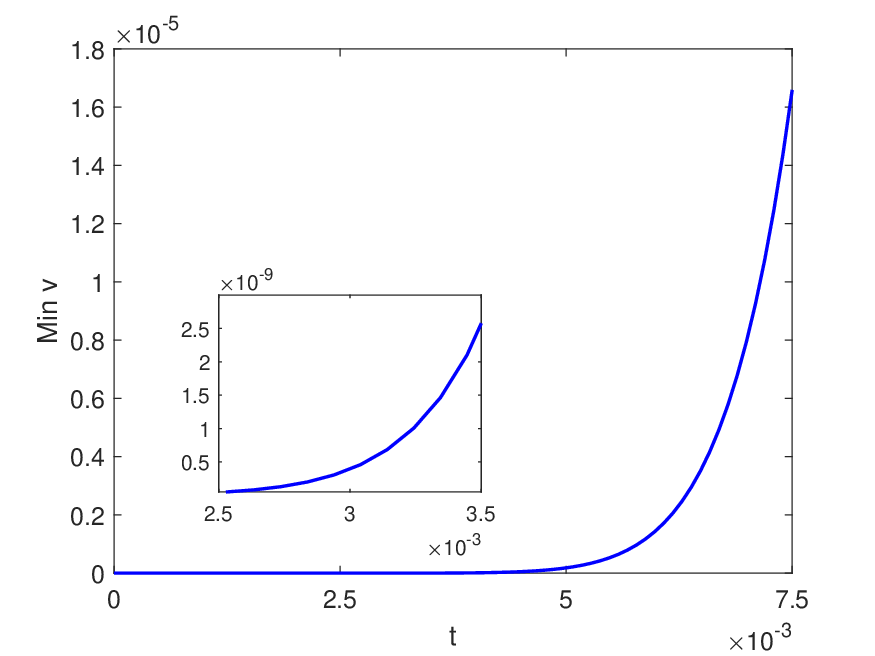}
			\includegraphics[width=0.325\linewidth]{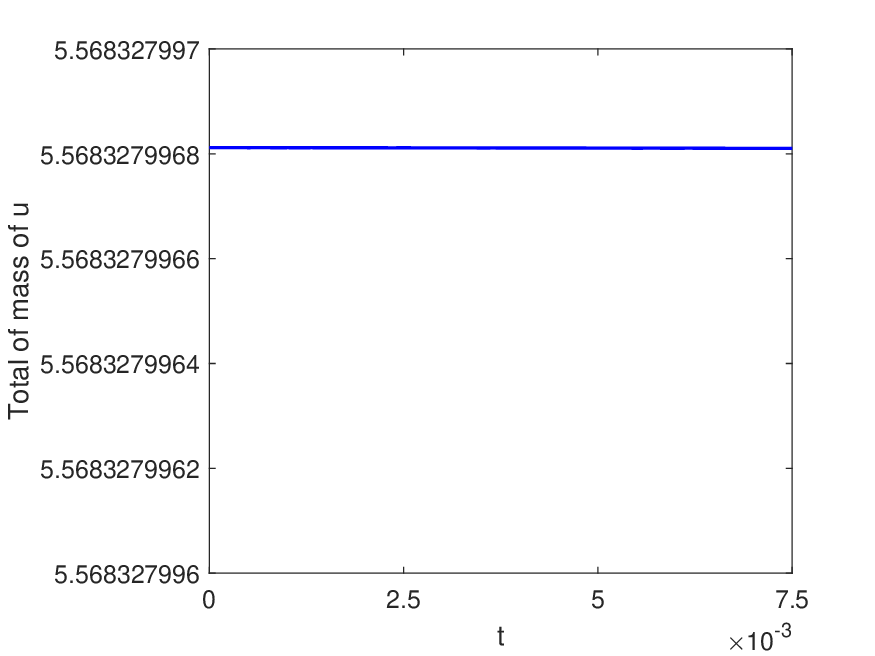}
			\includegraphics[width=0.325\linewidth]{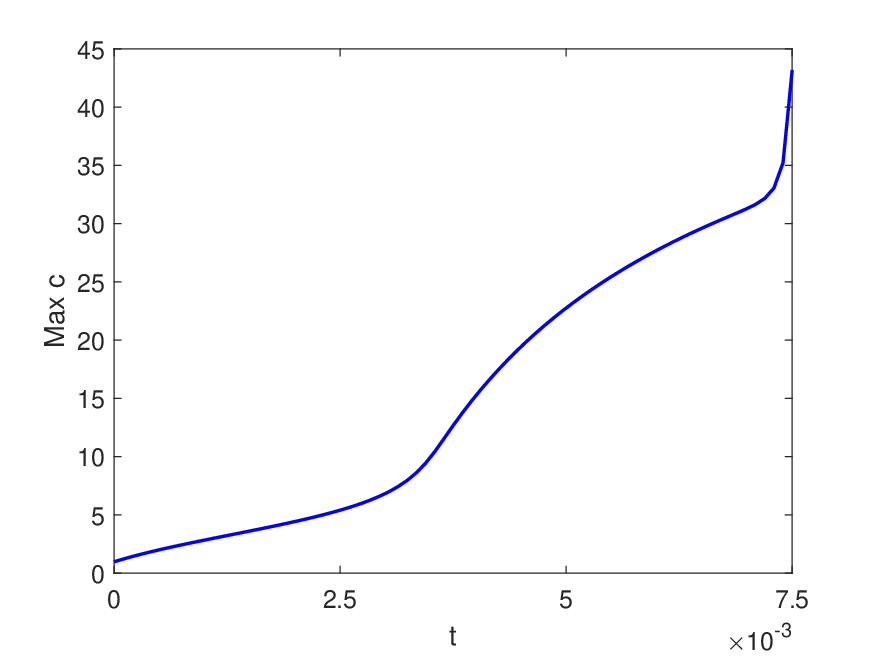}
			\includegraphics[width=0.325\linewidth]{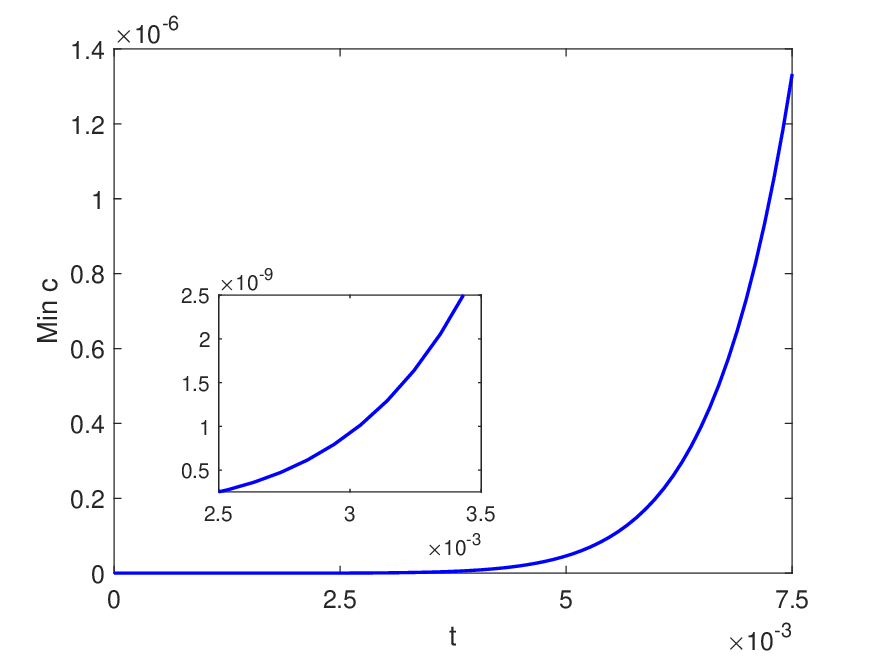}
			\includegraphics[width=0.325\linewidth]{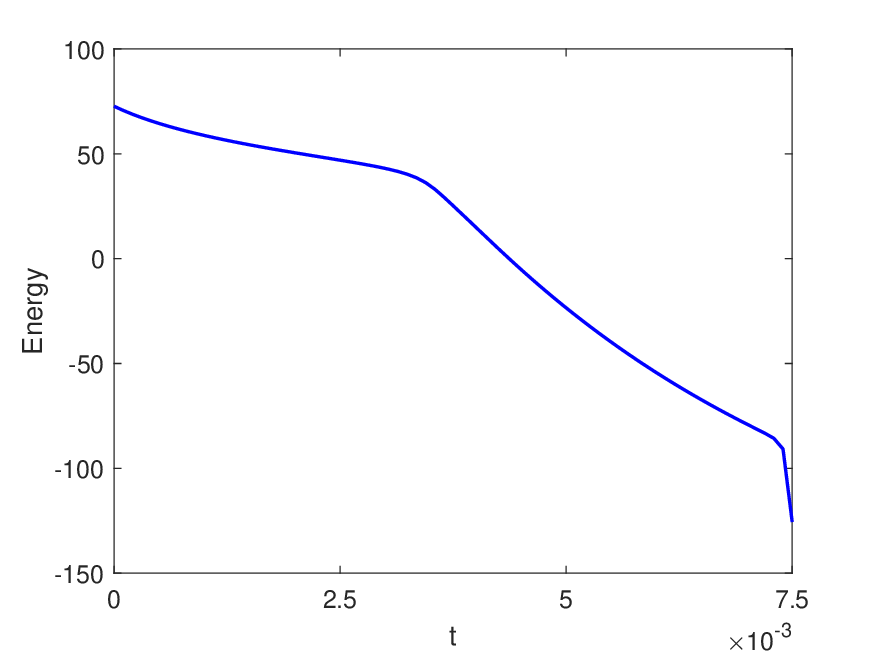}
			\caption{Evolutions of the extrema of $u$, $v$, $c$, and the total mass of $u,v$ and original energy (from left to right) for Example \ref{exam:bu:3D}.}
			\label{fig:eg6:mass}
	\end{figure}

\begin{figure}[htbp]
			\centering
			\includegraphics[width=0.24\linewidth]{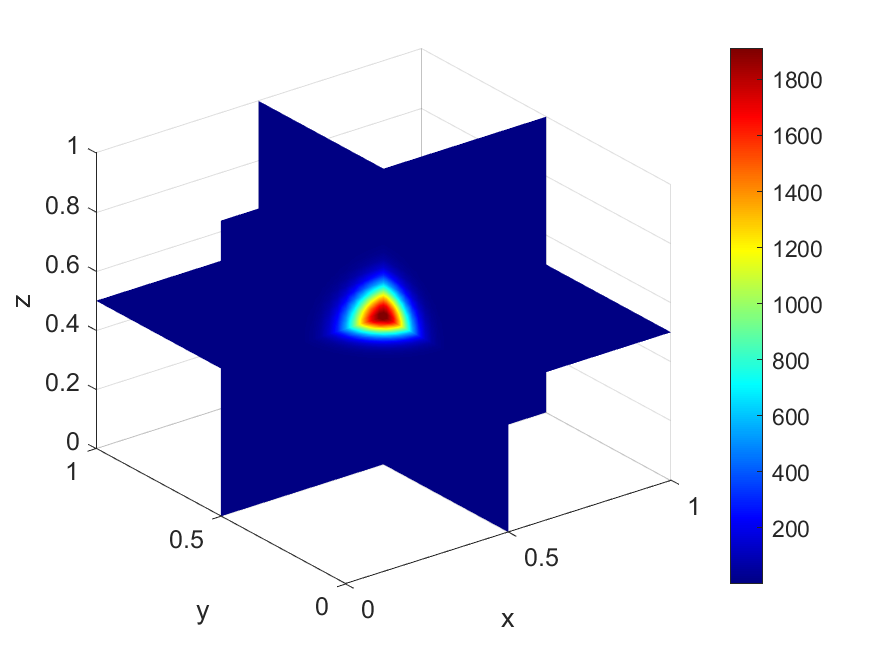}
            \includegraphics[width=0.24\linewidth]{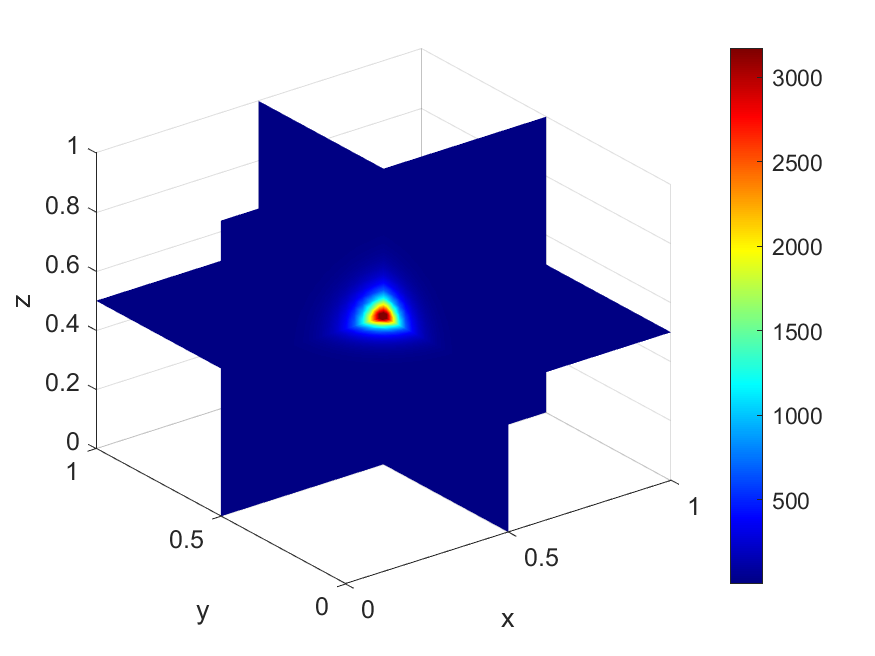}
			\includegraphics[width=0.24\linewidth]{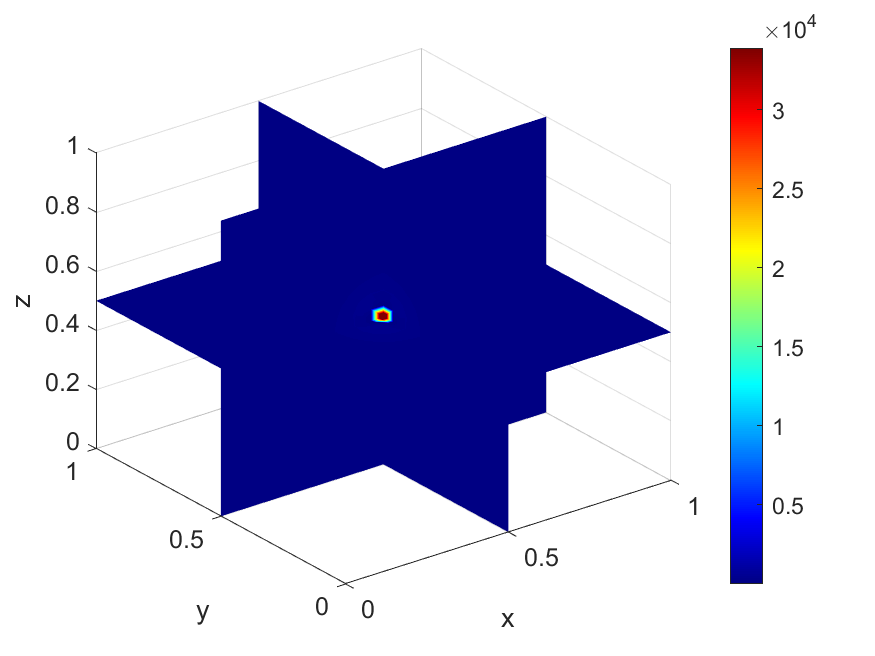}
			\includegraphics[width=0.24\linewidth]{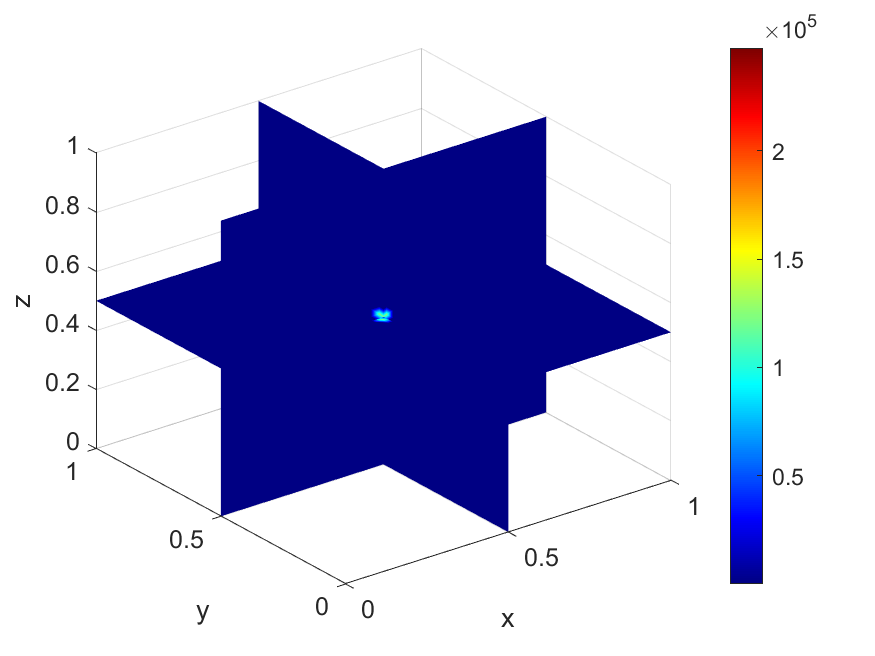}
			\includegraphics[width=0.24\linewidth]{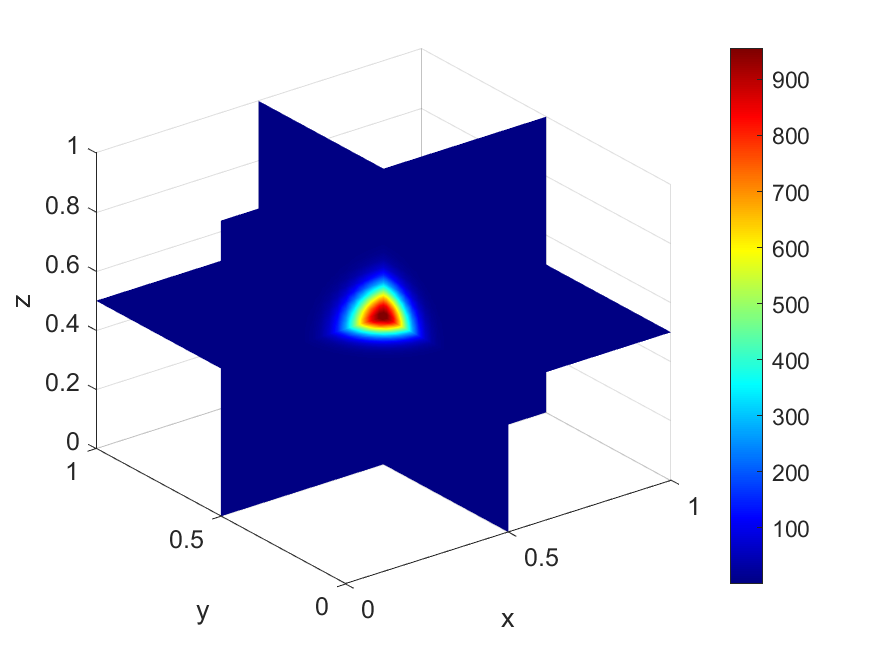}
            \includegraphics[width=0.24\linewidth]{t_0_v_eg4_s_new.eps}
			\includegraphics[width=0.24\linewidth]{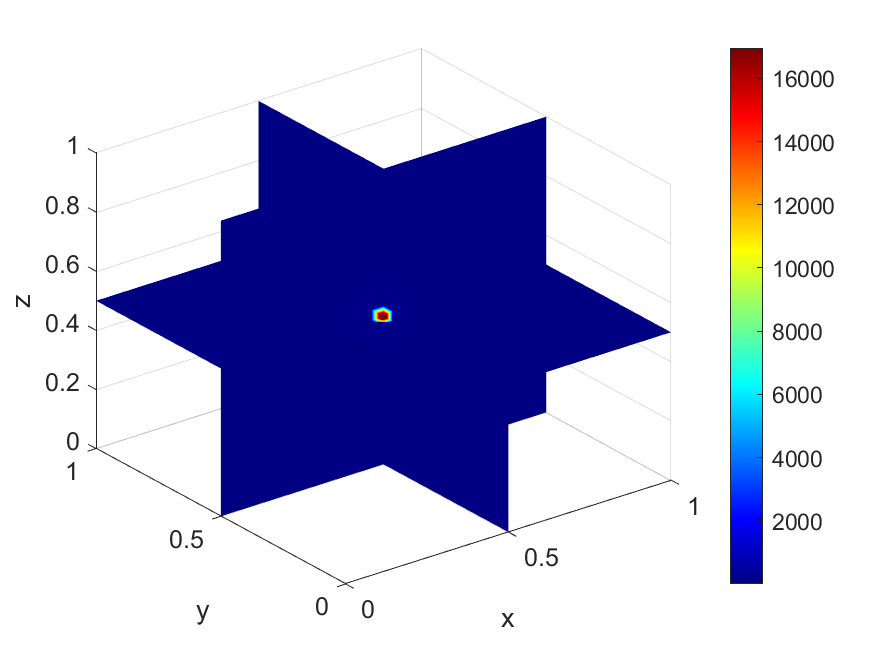}
			\includegraphics[width=0.24\linewidth]{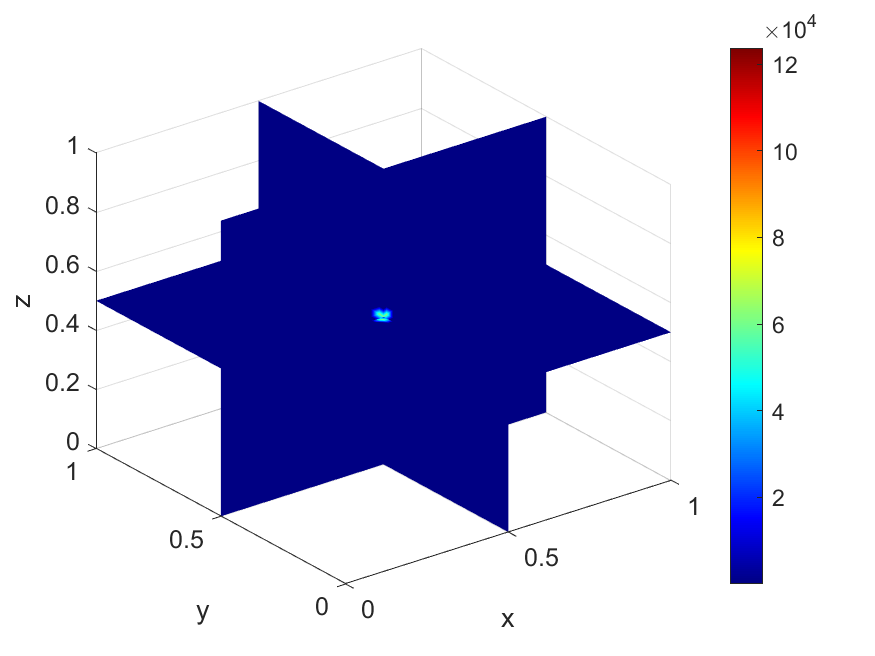}
			\includegraphics[width=0.24\linewidth]{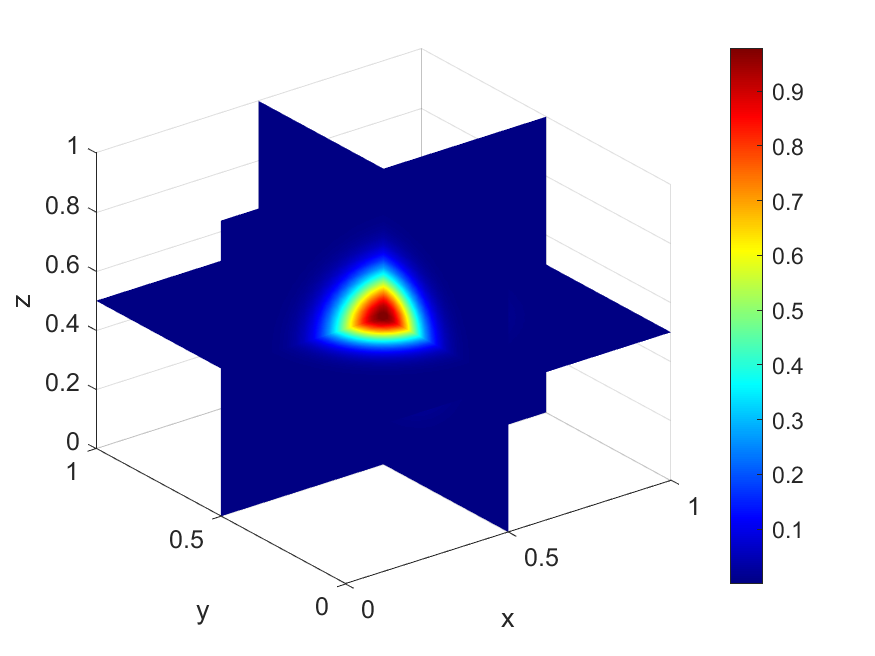}
            \includegraphics[width=0.24\linewidth]{t_0_c_eg4_s_new.eps}
			\includegraphics[width=0.24\linewidth]{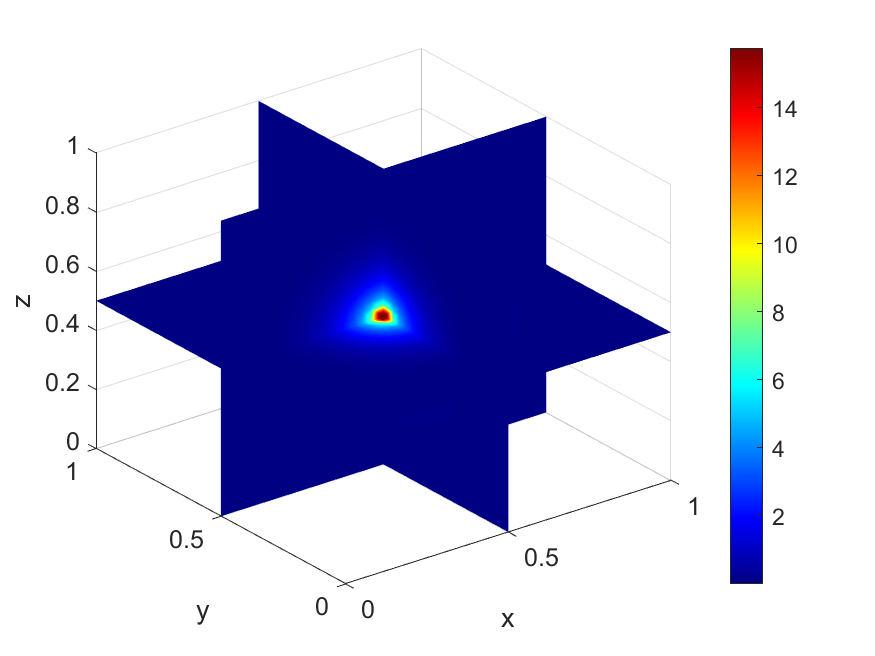}
			\includegraphics[width=0.24\linewidth]{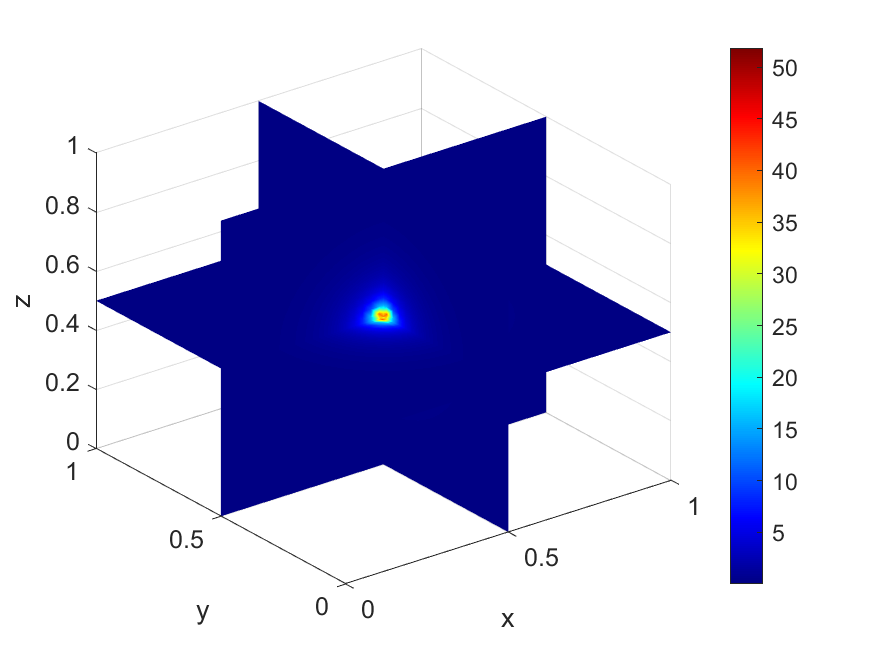}
			\caption{ Slices of $u,v,c$(from top to bottom)  at time instants $t=0,\ 2\times10^{-3},\ 4.0\times 10^{-3},\ 7.5\times 10^{-3}$ (from left to right) with $M = 40$ for Example \ref{exam:bu:3D}. }
			\label{fig:3D}
	\end{figure}

We simulate the structure-preserving property of the two-species Keller-Segel model in the 3D case, using $M=40$ grid points in each spatial direction, and a time stepsize of $\tau = 1.0 \times 10^{-4}$. Fig. \ref{fig:eg6:mass} presents the evolution of the maximum and minimum values of $u$, $v$, and $c$, together with the individual masses of two species and the discrete original energy. The numerical results confirm that the structure-preserving properties established in the preceding 2D experiments are retained in the 3D setting. To further visualize spatial dynamics, Fig. \ref{fig:3D} displays snapshots of the cell densities $u$ and $v$ and the chemoattractant concentration $c$ at times $t = 0,\, 4.0 \times 10^{-3},\, 7.5 \times 10^{-3}$, obtained from three orthogonal cross-sections at $x = 0.5$, $y=0.5$, and $z=0.5$. The cell densities exhibit increasingly pronounced aggregation near the center of the domain, with a sharp concentration developing as $t$ approaches $7.5\times10^{-3}$, suggesting finite-time blow-up behavior. These findings demonstrate that the proposed scheme captures the characteristic aggregation dynamics observed in the 2D case of Example \ref{exam:eg4:pp} while maintaining its structure-preserving properties in the 3D case.
	
%		
		
%%%%%%%%%%%%%%%%%%%%%%%%%%%%%%%%%%%%%%%%%%
	\section{Conclusion}\label{sec:conclusion}
%%%%%%%%%%%%%%%%%%%%%%%%%%%%%%%%%%%%%%%%%%
    In this paper, we have developed a fully discrete, linearly decoupled, structure-preserving BCFD scheme based on Slotboom transformation for the two-species Keller--Segel chemotaxis system. The proposed scheme effectively reduces the computational complexity associated with nonlinear coupling while preserving all fundamental physical properties of the continuous system, including the positivity of the cell densities and chemoattractant concentration, mass conservation for each cell species, and dissipation of the original energy. Moreover, the unique solvability of the numerical solutions at each time step is established. By employing high-order consistency analysis, we rigorously establish optimal-order error estimates for the densities in the discrete $L^2$-norm and for the chemoattractant concentration in the discrete $H^1$-norm, as stated in Theorem \ref{thm:covergence}. Extensive numerical experiments further confirm the accuracy and structure-preserving performance of the proposed scheme, including its ability to capture finite-time blow-up behavior in both 2D and 3D cases.
        
    Future work will focus on extending the present approach to high-order-in-time linearly decoupled structure-preserving numerical schemes for the Keller--Segel chemotaxis system, as well as more general chemotaxis–fluid coupled systems (see Refs. \cite{SN'19, HFXW'21, FHW'21,ZWP'26,WWQZ'26,BHP'25}).

        % In this paper, we have developed a fully discrete, linearly decoupled, and structure-preserving block-centered finite difference (BCFD) scheme based on the Slotboom transformation for the two-species Keller--Segel chemotaxis system \eqref{model:KS}. By reformulating the system into a symmetric form, the proposed method overcomes the computational challenges associated with nonlinear coupling and preserves the fundamental properties of the continuous model at the discrete level. The positivity of the cell densities and chemoattractant concentration, mass conservation of the cell densities, dissipation of the original energy, and unique solvability of the numerical solutions are rigorously established. Furthermore, optimal-order error estimates are obtained for the cell densities in the discrete $L^2$ norm and for the chemoattractant concentration in the discrete $H^1$ norm.
        
        % Extensive numerical experiments are conducted to validate the theoretical results and demonstrate the accuracy, stability, and efficiency of the proposed scheme. Future work will be devoted to the construction of high-order, linearly decoupled, and structure-preserving numerical methods for the two-species Keller--Segel chemotaxis system, as well as their extensions to more complex chemotaxis--fluid coupled systems; see Refs. \cite{SN'19, HFXW'21, FHW'21,ZWP'26,WWQZ'26,BHP'25}.
		
		%%%%%%%%%%%%%%%%%%%%%%%%%%%%%%%%%%%%%%%%%%%%%%%%%%%%%%		
		\section*{CRediT authorship contribution statement}
		\textbf{Ao Zhang}: Methodology, Formal analysis, Software, Writing-Original draft.
		\textbf{Bingyin Zhang}: Methodology, Formal analysis, Writing-Original draft.
		\textbf{Hongfei Fu}: Conceptualization, Supervision, Writing-Reviewing and Editing, Methodology, Funding acquisition.
		
		\section*{Data availability} 
        Data are available upon reasonable request.
        
		\section*{Acknowledgements}
		This work was partially supported by the National Natural Science Foundation of China under grant 12131014, by the Shandong Provincial Natural Science Foundation under grant ZR2024MA023, and by the Shandong Provincial Key Laboratory of Stochastic System Control and Scientific Computing.

        \section*{Declarations} 
        \textbf{Conflict of interest} The authors declare that they have no competing interests.
		
		%%%%%%%%%%%%%%% References %%%%%%%%%%%%%%%%%%%%%%%
		\bibliographystyle{elsarticle-num}
		\bibliography{reference.bib}
		
	\end{document}